\documentclass{article}
\usepackage{arxiv}
\usepackage[utf8]{inputenc} 
\usepackage[T1]{fontenc}    
\usepackage{hyperref}       
\hypersetup{
	pdftitle={Eigenvalue Stabilization Technique},
	pdfsubject={Computational Mechanics},
	pdfauthor={S.~Eisentraeger},
	pdfkeywords={Computational Mechanics},
	colorlinks=false,
	hidelinks=true,
}
\usepackage{url}            
\usepackage{booktabs}       
\usepackage{longtable}
\usepackage{colortbl}
\usepackage{overpic}
\usepackage{amsfonts}       
\usepackage{amsmath}
\usepackage{amssymb}
\usepackage{mathtools}
\usepackage{cancel}
\usepackage{xfrac}			
\usepackage{microtype}      
\usepackage{cleveref}       
\usepackage{graphicx}
\usepackage{doi}
\usepackage{subfig}
\usepackage{placeins}
\usepackage{scalerel}
\usepackage{diagbox}		
\usepackage{arydshln}		
\usepackage{multirow}
\usepackage[bottom,symbol]{footmisc}	
\usepackage[noadjust]{cite}
\usepackage{xcolor}
\definecolor{Matlab1}{rgb}{0.0,0.0,0.0}
\definecolor{Matlab2}{rgb}{0.878431372549020,0.015686274509804,0.015686274509804}
\definecolor{Matlab3}{rgb}{0.101960784313725,0.352941176470588,0.627450980392157}
\definecolor{Matlab4}{rgb}{0.6,0.6,0.6}
\definecolor{Matlab5}{rgb}{1.0,0.443137254901961,0.0}
\definecolor{Matlab6}{rgb}{0.274509803921569,0.560784313725490,1.0}
\definecolor{Matlab7}{rgb}{1.0,0.0,0.0}
\definecolor{Matlab8}{rgb}{0.0,0.0,1.0}
\definecolor{Matlab9}{rgb}{0.0,1.0,0.0}
\definecolor{Grey1}{rgb}{0.5,0.5,0.5}
\definecolor{Grey2}{rgb}{0.9,0.9,0.9}
\usepackage{bbding}
\usepackage{pifont}
\usepackage{relsize}

\usepackage{contour}
\newcommand{\fett}[1]{\mbox{\boldmath$#1$}} %
\newcommand{\fatgreek}[1]{{\contourlength{0.001em}\contour{black}{\fett{#1}}}}
\usepackage{lineno}
\modulolinenumbers[2]
\let\originalleft\left
\let\originalright\right
\def\left#1{\mathopen{}\originalleft#1}
\def\right#1{\originalright#1\mathclose{}}
\usepackage{tikz}

\usepackage{enumitem}
\graphicspath{{./figs/}} 
\title{An EigenValue Stabilization Technique for Immersed Boundary Finite Element Methods in Explicit Dynamics}
\author{
\href{https://orcid.org/0000-0001-8774-9732}{\includegraphics[scale=0.06]{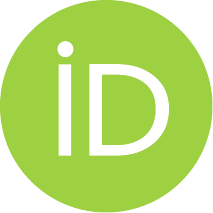}\hspace{1mm}S.~Eisentr\"ager$^{1,}$\footnote[1]{Corresponding author}} \hspace*{1ex}and D.~Juhre$^{6}$ \\
Institute of Mechanics\\
Otto von Guericke University Magdeburg\\
Magdeburg, 39106, Germany\\
\And
L.~Radtke$^{2}$, W.~Garhuom$^{3}$ and A.~D\"uster$^{5}$ \\
Numerical Structural Analysis with \\
Application in Ship Technology (M-10)\\
Hamburg University of Technology\\
Hamburg, 21073, Germany
\And
S.~L\"ohnert$^{4}$\\
Institute of Mechanics and Shell Structures\\
Technische Universit\"at Dresden\\
Dresden, 01062, Germany
\And
D.~Schillinger$^{7}$\\
Institute for Mechanics\\
Technical University of Darmstadt\\
Darmstadt, 64287, Germany
}
\renewcommand{\headeright}{Preprint submitted to \textit{CAMWA}}
\renewcommand{\undertitle}{Preprint submitted to \textit{CAMWA}}
\renewcommand{\shorttitle}{Eigenvalue stabilization technique}
\begin{document}
%
%
\maketitle
\begin{abstract}
The application of immersed boundary methods in static analyses is often impeded by poorly cut elements (small cut elements problem), leading to ill-conditioned linear systems of equations and stability problems. While these concerns may not be paramount in explicit dynamics, a substantial reduction in the critical time step size based on the smallest volume fraction $\chi$ of a cut element is observed. This reduction can be so drastic that it renders explicit time integration schemes impractical. To tackle this challenge, we propose the use of a dedicated eigenvalue stabilization (EVS) technique.\\
The EVS-technique serves a dual purpose. Beyond merely improving the condition number of system matrices, it plays a pivotal role in extending the critical time increment, effectively broadening the stability region in explicit dynamics. As a result, our approach enables robust and efficient analyses of high-frequency transient problems using immersed boundary methods. A key advantage of the stabilization method lies in the fact that only element-level operations are required.\\
This is accomplished by computing all eigenvalues of the element matrices and subsequently introducing a stabilization term that mitigates the adverse effects of cutting. Notably, the stabilization of the mass matrix $\mathbf{M}_\mathrm{c}$ of cut elements -- especially for high polynomial orders $p$ of the shape functions -- leads to a significant raise in the critical time step size $\Delta t_\mathrm{cr}$.\\
To demonstrate the efficacy of our technique, we present two specifically selected dynamic benchmark examples related to wave propagation analysis, where an explicit time integration scheme must be employed to leverage the increase in the critical time step size.
\end{abstract}
%
\footnote[0]{$^*$ Corresponding author: \texttt{sascha.eisentraeger@ovgu.de}; the order of the authors is indicated by a numeric (arabic numbers) superscript, e.g., Surname$^i$.}
\renewcommand*{\thefootnote}{\arabic{footnote}}
\setcounter{footnote}{0}
\section{Introduction}
\label{sec:Intro}
Transient analyses play a crucial role in various scientific and engineering fields. However, despite the significant computational resources available today, solving high-frequency dynamics problems in the time domain remains challenging. The demanding requirements for fine spatial and temporal resolutions call for highly efficient algorithms, which are not readily accessible at present. Although the conventional finite element method (FEM) is widely used for numerical analysis across diverse problem domains, it does have limitations. For instance, it lacks an automated mesh generation pipeline for complex structures, high-order convergence, and advanced mass lumping techniques. 

To address these shortcomings, alternative numerical methods have gained considerable traction over the past two decades. Notable approaches include isogeometric analysis (IGA)~\cite{BookCottrell2009} and immersed boundary methods\footnote{Remark: The article uses the term ``immersed boundary methods'' interchangeably with ``fictitious domain methods'' or ``embedded domain methods'', aligning with their widespread usage in other relevant literature.}. In the context of IGA, a mathematical theory of mass lumping has recently been developed in Ref.~\cite{ArticleVoet2022}. Furthermore, ongoing efforts are focused on devising high-order convergent mass lumping techniques. These methods rely on approximate dual shape functions within a Petrov-Galerkin framework, as detailed in Refs.~\cite{ArticleNguyen2024} and \cite{ArticleHeld2024}. Nonetheless, several challenges remain to be addressed, including the implementation of outlier removal techniques \cite{ArticleHiemstra2021}, multi-patch analysis, and trimming. The integration of these elements into the framework is essential to render it suitable for real-world problems.

In this article, we primarily focus on immersed boundary methods, specifically examining the finite cell method (FCM)~\cite{ArticleParvizian2007, ArticleDuester2008, ArticleElhaddad2015} and spectral cell method (SCM)~\cite{ArticleDuczek2014, ArticleJoulaian2014a}, its extension to explicit dynamic problems. It is worth mentioning other notable immersed approaches, including CutFEM~\cite{ArticleBurman2022}, Cartesian grid FEM (cgFEM)~\cite{ArticleJimenez2020}, and Aggregated FEM~\cite{ArticleBadia2021} to name just a few. While our discussions focus on the FCM/SCM, the techniques presented in this paper are applicable to all types of immersed boundary methods being based on finite element technologies.

Immersed boundary methods provide a means of discretizing structures using Cartesian grids, employing elements that do not conform to the boundary of the geometry of interest. However, this non-conforming spatial discretization introduces three challenges that need to be addressed. First, accurate numerical evaluation of integrals over cut elements requires sophisticated techniques. Second, the imposition of Dirichlet (essential) boundary conditions becomes more complex. Third, stability and conditioning issues arise for cut elements that intersect the physical domain boundary.

In this paper, we will focus on the third problem and briefly explore potential remedies. It is well-known that ill-conditioning and stability problems occur when cut elements are sparsely filled with material, lacking sufficient support in the physical domain \cite{ArticleDePrenter2017}. Therefore, stabilizing the fictitious domain is a common requirement across all immersed boundary methods to prevent severe numerical issues.

Various stabilization techniques have been developed in recent years to address this challenge. For example, the ghost penalty method is utilized in CutFEM applications \cite{ArticleBurman2010a, ArticleBurman2012}. This technique involves introducing an additional term to the weak form, penalizing jumps in the normal derivatives of shape functions between neighboring cut and non-cut elements. It effectively supports the cut elements using the interior non-cut elements, improving the condition number without significantly altering the underlying mathematical problem.

Another stabilization approach, known as the fictitious material or $\alpha$-method, is commonly employed in the FCM~\cite{ArticleDuester2008, InbookDuester2018}. Here, a very soft material is introduced in the fictitious domain of every cut element, yielding favorable results across various applications. However, this method has the drawback of adding the same artificial stiffness to all points in the fictitious domain, potentially modifying the solution and decreasing the overall accuracy. Additionally, its performance is limited for nonlinear problems, imposing restrictions on the stability in numerical analysis and thus, on the achievable deformation at finite strains \cite{ArticleGarhuom2022}.

The basis function removal (BFR) strategy has been proposed to eliminate shape functions with minimal contributions to the global stiffness matrix \cite{ArticleVerhoosel2015}. While simple in concept, this technique is not reliable in curing the small cut element problem and lacks robustness. However, when combined with a dedicated remeshing strategy, favorable outcomes can be achieved by creating a new mesh when the old mesh can no longer accommodate further deformations (highly distorted elements) \cite{ArticleGarhuom2020}.

Tailor-made preconditioning techniques, such as the symmetric incomplete permuted inverse Cholesky preconditioner \cite{ArticleDePrenter2017}, offer another avenue to mitigate conditioning issues. Furthermore, the use of an additive Schwarz preconditioner in conjunction with multigrid techniques has been proposed, demonstrating robustness and efficiency \cite{ArticleDePrenter2019, ArticleJomo2021}, while the feasibility of applying these preconditioners to practical problems within high-performance computing environments (parallel-computing) has been investigated in Ref.~\cite{ArticleJomo2019}. The study showcased their promising scalability properties, making them valuable tools for real-world applications.

For a more in-depth exploration of strategies for effectively addressing the challenges associated with small cut elements, we recommend consulting the recent review article by de~Prenter et al. \cite{ArticlePrenter2022}. This comprehensive review offers valuable insights and a detailed analysis of various techniques and methodologies used in handling poorly cut elements within immersed boundary methods. The article also covers important topics, including numerical integration over cut elements, the imposition of boundary conditions, stability and conditioning issue, as well as dedicated stabilization techniques. It serves as a valuable resource for researchers and practitioners seeking a deeper understanding of the complexities associated with small cut elements, providing guidance on selecting suitable approaches for different situations.

Within the context of the extended finite element method (XFEM), which can also be seen as a fictitious domain technique, particularly when analyzing void regions, a set of techniques has been established that are theoretically applicable to other immersed boundary methods such as the FCM or CutFEM. One option is to utilize the node-moving technique \cite{ArticleLoehnert2014}. When an element has minimal intersection with the physical domain, the volume fraction can be increased by moving nodes within the physical domain. However, this approach sacrifices the advantages of Cartesian meshes. To address poorly conditioned stiffness matrices, Bechet et al. developed a specialized preconditioner for enriched finite elements based on Cholesky decompositions of submatrices \cite{ArticleBechet2015}. Menk and Bordas proposed another effective preconditioner similar to the FETI (finite element tearing and interconnecting \cite{ArticleFarhat1991}) domain decomposition technique \cite{ArticleMenk2011}. Babu\v{s}ka and Banerjee introduced a stabilization technique that tackles convergence issues arising from high condition numbers of the global stiffness matrix \cite{ArticleBabuska2012}. This technique is solver-independent and particularly useful for large-scale 3D simulations that frequently employ library-based parallel iterative equation solvers with various preconditioning techniques.

In Refs.~\cite{ArticleLoehnert2014, ArticleBeese2017}, a simple -- yet efficient -- stabilization technique based on an eigenvalue decomposition of the elemental stiffness matrix was proposed. Referred to as the EVS-technique (eigenvalue stabilization) throughout the article, this method operates at the element level, ensuring minimal additional numerical costs and high parallelizability. The stabilization matrices are calculated during the element assembly process, which not only preserves method flexibility, but also facilitates its implementation in high-performance computing environments.

Previously, the EVS-scheme has proven successful in handling quasi-static and dynamic crack propagation problems. Note that in these applications, the focus was primarily on solving ill-conditioning problems related to enriched elements. Its extension to immersed boundary methods and the FCM, in particular, was achieved in Ref.~\cite{ArticleGarhuom2022}. Here, the EVS-technique was implemented to reduce the condition number of cut elements without significantly affecting solution quality. To ensure accurate results in nonlinear analyses (e.g., hyperelastic material model at finite strains), an iteratively updated force correction term was incorporated into the solution procedure. The computational overhead incurred by applying the EVS-technique in nonlinear analyses is limited for two reasons: First, the required eigenvalue decomposition is performed at the element level, considering eigenvalues and mode shapes only for cut elements. Second, nonlinear analyses inherently require an incremental/iterative solver (e.g., Newton-Raphson algorithm), naturally incorporating the iterative correction scheme. In this process, modes with small or zero eigenvalues are grouped and stabilized. It is essential to note that without a stabilization technique critical modes can render the FCM less robust, especially when dealing with badly cut elements and high-order shape functions. Hence, a dedicated stabilization technique, targeting specific parts of the stiffness  and/or mass matrices, becomes crucial to effectively address these issues.

This contribution presents a further extension of the EVS-technique to encompass dynamics, specifically explicit time stepping. The novel approach pursues two primary objectives: First, increasing the critical time step size in explicit dynamics and second, reducing the condition number of the system matrices. The first objective is crucial for efficiently analyzing wave propagation or impact problems (e.g., crash tests) using explicit time stepping methods. On the other hand, the second objective only holds significant importance for large-scale simulations utilizing implicit time stepping methods, where iterative solvers are commonly employed and conditioning issues can significantly impede convergence. It is important to note that these two objectives are interconnected to a certain extent and not mutually exclusive. However, our findings, as presented in this contribution, demonstrate that for explicit analyses, we can leave the stiffness matrix unstabilized, as its impact on the critical time step size is minimal. In contrast, this would not be a viable option for implicit schemes.

In contrast to prior implementations of the EVS-technique, which primarily focused on stabilizing the stiffness matrix to mitigate ill-conditioning, our novel approach advocates for the stabilization of the mass matrix. As previously noted, in the context of (linear) explicit dynamics, the stabilization of the stiffness matrix is of lesser significance since simple matrix--vector products are sufficient to advance in time. Consequently, there is no need to solve linear systems of equations, where concerns related to conditioning typically manifest.

To address time-dependent problems in the context of immersed boundary methods, the authors have previously developed the SCM \cite{ArticleDuczek2014,ArticleJoulaian2014a}, a specialized variant of the FCM. Unlike the FCM, the SCM utilizes nodal shape functions based on Lagrangian interpolation polynomials defined on non-equidistant nodal distributions. This unique approach enables the application of mass lumping techniques, which are crucial for highly efficient explicit time stepping algorithms. Consequently, the integration of the SCM with the EVS-technique and explicit time stepping offers a compelling framework for complex dynamic simulations. This approach provides notable benefits, particularly in terms of increased computational efficiency. Thus, by leveraging the EVS-technique, a robust framework, enabling researchers and practitioners to confidently tackle challenging dynamic simulations across diverse domains, including structural dynamics, crash simulations, and acoustic analyses, is proposed.
\section{Governing equations of elastodynamics and finite element discretization}
\label{sec:Dyn}
This article focuses on problems within the field of linear elastodynamics \cite{BookChopra2011}. Specifically, we address these problems by utilizing a variational formulation of the following form:
\begin{equation}
	\mathcal{A}(\ddot{\mathbf{u}},\mathbf{v}) + \mathcal{C}(\dot{\mathbf{u}},\mathbf{v}) + \mathcal{B}(\mathbf{u},\mathbf{v}) = \mathcal{F}(\mathbf{v})\,,
	\label{eq:DynVar}
\end{equation}
with
\allowdisplaybreaks
\begin{alignat}{3}
		&\mathcal{A}(\ddot{\mathbf{u}},\mathbf{v}) && = && \int\limits_\Omega \rho \mathbf{v}^\mathrm{T}\ddot{\mathbf{u}}\,\mathrm{d}\Omega\,, \\
		&\mathcal{C}(\dot{\mathbf{u}},\mathbf{v}) && = && \int\limits_\Omega \mathbf{v}^\mathrm{T}\fatgreek{\kappa}\dot{\mathbf{u}}\,\mathrm{d}\Omega\,, \\
		&\mathcal{B}(\mathbf{u},\mathbf{v}) && = && \int\limits_\Omega \left[\mathbf{L}\mathbf{v}\right]^\mathrm{T}\mathbb{C}\left[\mathbf{L}\mathbf{u}\right]\,\mathrm{d}\Omega\,, \\
		\shortintertext{and} & && \nonumber \\[-18pt]
		&\mathcal{F}(\mathbf{v}) && = && \int\limits_\Omega \mathbf{v}^\mathrm{T}\mathbf{f}\,\mathrm{d}\Omega + \int\limits_{\Gamma_\mathrm{N}} \mathbf{v}^\mathrm{T}\bar{\mathbf{t}}\,\mathrm{d}\Gamma + \mathbf{v}^\mathrm{T}\mathbf{f}_\mathrm{p}\,.
\end{alignat}
The variational formulation includes the bi-linear forms $\mathcal{A}$, $\mathcal{C}$, and $\mathcal{B}$, along with the linear form $\mathcal{F}$, representing terms associated with inertia, damping, stiffness, and external forces, respectively. In this context, $\rho$ denotes the mass density, $\fatgreek{\kappa}$ represents the matrix of damping parameters, $\mathbb{C}$ is the elasticity/constitutive matrix, and $\mathbf{f}$ stands for the vector of volume loads acting on the domain $\Omega$. The trial and test functions, denoted as $\mathbf{u}$ and $\mathbf{v}$, respectively, typically correspond to the displacement vector. The first and second temporal derivatives are represented by $\dot{\square}$ and $\ddot{\square}$. Furthermore, the linear strain-displacement operator is denoted as $\mathbf{L}$. Neumann boundary conditions, such as surface tractions $\bar{\mathbf{t}}$ and point forces $\mathbf{f}_\mathrm{p}$, are prescribed along the boundary $\Gamma_\mathrm{N}$ or at individual points. To complete the set of equations, Dirichlet boundary conditions must be imposed along the boundary $\Gamma_\mathrm{D}$
\begin{equation}
	\mathbf{u} = \bar{\mathbf{u}} \quad \text{on} \quad \Gamma_\mathrm{D}\,.
\end{equation}
In the case of dynamical problems, it is also necessary to consider initial conditions for the displacement and velocity fields, which are given by
\begin{alignat}{2}
	\mathbf{u}(t=0) = \mathbf{u}_0\,, \\
	\shortintertext{and}  \nonumber \\[-18pt]
	\dot{\mathbf{u}}(t=0) = \dot{\mathbf{u}}_0\,.
\end{alignat}
This set of equations is defined within the physical domain, assuming a geometry-conforming discretization. However, in the next section (Sect.~\ref{sec:FCM}), we will briefly outline an extension to non-conforming meshes using immersed boundary methods.

Following the standard FEM-procedure \cite{BookZienkiewicz2000a}, the displacement field within each finite element (with element domain $\Omega_e$) is approximated using simple polynomial shape functions, given by
\begin{equation}
	\mathbf{u}(\mathbf{x},t) = \mathbf{N}(\mathbf{x})\mathbf{U}_e(t) \qquad \forall \; \,\mathbf{x} \in \Omega_e\,,
\end{equation}
where $\mathbf{N}$ represents the matrix of shape functions and $\mathbf{U}$ denotes the vector of nodal displacements. It is worth noting that, in a Bubnov-Galerkin approach, the same shape functions are also employed for the test functions
\begin{equation}
	\mathbf{v}(\mathbf{x},t) = \mathbf{N}(\mathbf{x})\mathbf{V}_e(t)\quad \forall \; \mathbf{x} \in \Omega_e\,.
\end{equation}
By substituting the discretized versions of the displacement field $\mathbf{u}$ and the test function $\mathbf{v}$ into the weak form, and performing some algebraic manipulations, we can assemble the elemental contributions, resulting in the well-known semi-discrete equations of motion
\begin{equation}
	\mathbf{M}\ddot{\mathbf{U}} + \mathbf{C}\dot{\mathbf{U}} + \mathbf{K}\mathbf{U} = \mathbf{F}_\mathrm{ext}\,,
	\label{eq:EquationMotion}
\end{equation}
where $\mathbf{M}$, $\mathbf{C}$, and $\mathbf{K}$ denote the mass, damping, and stiffness matrices, respectively. The external load vector is represented by $\mathbf{F}_\mathrm{ext}$. At this point we want to point out that the damping matrix $\mathbf{C}$ is obtained by using Rayleigh's hypothesis, resulting in a linear combination of the mass and stiffness matrices
\begin{equation}
	\mathbf{C} = \alpha_\mathrm{R}\mathbf{M} + \beta_\mathrm{R}\mathbf{K}\,.
	\label{eq:RayleighDamping}
\end{equation}
The coefficients $\alpha_\mathrm{R}$ and $\beta_\mathrm{R}$ determine the effective range of damping. For further insight into determining these coefficients, especially in more complex applications, refer to Ref.~\cite{ArticleSong2017}. It should be noted that in explicit analyses, $\beta_\mathrm{R}$ is typically not considered. We have to realize that including a stiffness-proportional damping term would render the effective (dynamic) stiffness matrix  $\tilde{\mathbf{K}}$ non-diagonal, assuming the application of the central difference method (CDM) to advance in time. More details regarding the temporal discretization and time stepping are introduced later in Sect.~\ref{sec:StabTimeInt}.
\section{Finite cell method}
\label{sec:FCM}
\begin{figure}[t!]
	\centering
	\begin{overpic}[width=1.0\textwidth]{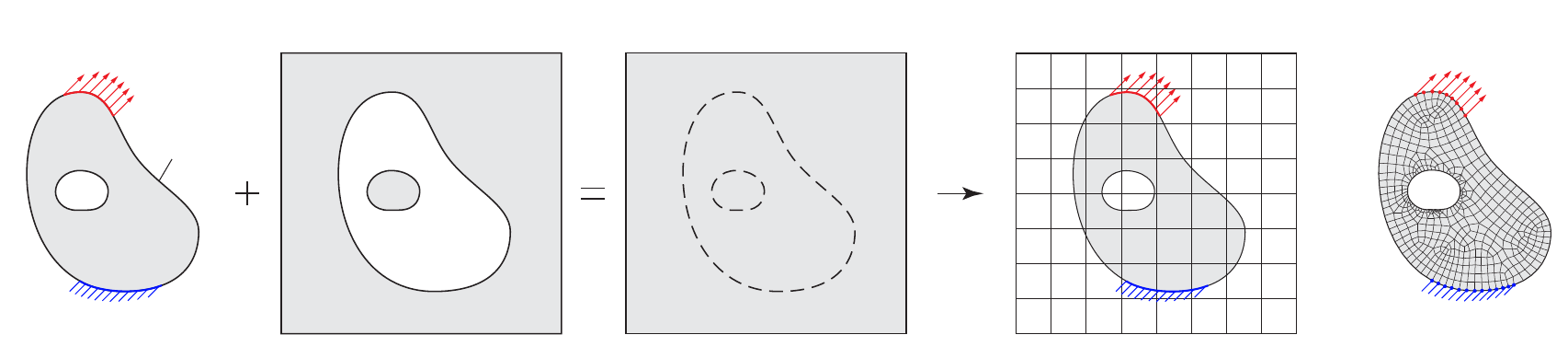}
		\put(0,20){Physical domain}
		\put(18,20){Fictitious domain}
		\put(40,20){Extended domain}
		\put(69.5,20){FC-mesh}
		\put(87,20){FE-mesh}
		\put(84,10){vs.}
		\put(3,2){$\hat{\mathbf{u}}$}
		\put(10,2){$\Gamma_\mathrm{D}$}
		\put(5,7){$\Omega_\mathrm{phys}$}
		\put(11,13){$\partial\Omega$}
		\put(9,17){$\Gamma_\mathrm{N}$}
		\put(2.5,17){$\hat{\mathbf{t}}$}
		\put(30,15){$\Omega_\mathrm{fict}$}
		\put(52,15){$\Omega^\mathrm{ex}$}
	\end{overpic}
	\caption{Fundamental idea of immersed boundary methods (Cartesian mesh) in comparison to a typical finite element discretization (body-fitted mesh).}
	\label{fig:FCM_concept_1}
\end{figure}

The FCM is a notable example of immersed boundary methods, which have undergone significant development in the past decade \cite{ArticleParvizian2007, ArticleDuester2008, ArticleSchillinger2014}. These numerical methods are designed to solve the governing partial differential equations (PDEs) on an extended domain $\Omega_\mathrm{ex}$ instead of the complex physical domain $\Omega_\mathrm{phys}$. This simplifies the meshing process, while necessitating advanced techniques for numerical integration of the system matrices \cite{ArticleKudela2016, ArticleJoulaian2016, ArticlePetoe2020, ArticleLegrain2021} and the imposition of Dirichlet boundary conditions \cite{ArticleBurman2012, ArticleSchillinger2016}. In this section, we focus on discussing the FCM in the context of void regions with arbitrary shapes. For a comprehensive analysis of multi-material problems that require an enrichment of the ansatz space, we refer interested readers to Refs.~\cite{ArticleJoulaian2013, ArticlePetoe2023a}. When only void regions are considered, the extended domain consists of two disjoint regions: the physical domain $\Omega_\mathrm{phys}$ and the fictitious domain $\Omega_\mathrm{fict}$, as illustrated in Fig.~\ref{fig:FCM_concept_1}. We want to stress at this point that the distinction between physical and fictitious domains is not necessary for geometry-aligning discretization methods such as the FEM and therefore, the notation slightly differs from that introduced in Sect.~\ref{sec:Dyn}.

To distinguish between physical and fictitious domains, an indicator function $\alpha_\mathrm{FCM}$, facilitating point-membership tests, is introduced
\begin{equation}
	\alpha_\mathrm{FCM}(\mathbf{x})= 
	\begin{cases}
		1.0,  & \forall\; \mathbf{x} \in \Omega_{\mathrm{phys}}\\
		\alpha_0,  & \forall\; \mathbf{x} \in \Omega_{\mathrm{fict}}
	\end{cases}.
	\label{eq:indicator_function}
\end{equation}%
The theoretically ideal value for the indicator function within the fictitious domain, denoted as $\alpha_0$, is 0, preserving the original form of the PDE. However, to prevent severe ill-conditioning problems, it is common practice to select a small positive value for $\alpha_0$, typically in the range of $10^{-12}$ to $10^{-5}$. The specific choice of $\alpha_0$ can also depend on the material properties, as demonstrated in Refs.~\cite{ArticleRuess2013, ArticleEisentraeger2022a}. Hence, the introduction of the $\alpha$-method, as it will be referred to throughout this article, involves adding a small stiffness to the fictitious domain, which serves to stabilize the equations. It should be noted that the use of the indicator function introduces a discontinuous function in the volume integrals of the weak form. This leads to non-smooth integrands, i.e., discontinuities in the integrals, which is a characteristic feature of finite cell-based numerical methods
\begin{equation}
	\int \limits_{\Omega_{\mathrm{phys}}} \mathcal{P} (\mathbf{x})~\mathrm{d}\Omega \approxeq 
	\int \limits_{\Omega^\mathrm{ex}} \alpha_\mathrm{FCM}(\mathbf{x}) \mathcal{P}(\mathbf{x})~\mathrm{d} \Omega\,.
	\label{eq:volume_integrals}
\end{equation}%

\begin{figure}[b!]
\centering
\begin{overpic}[width=1.0\textwidth]{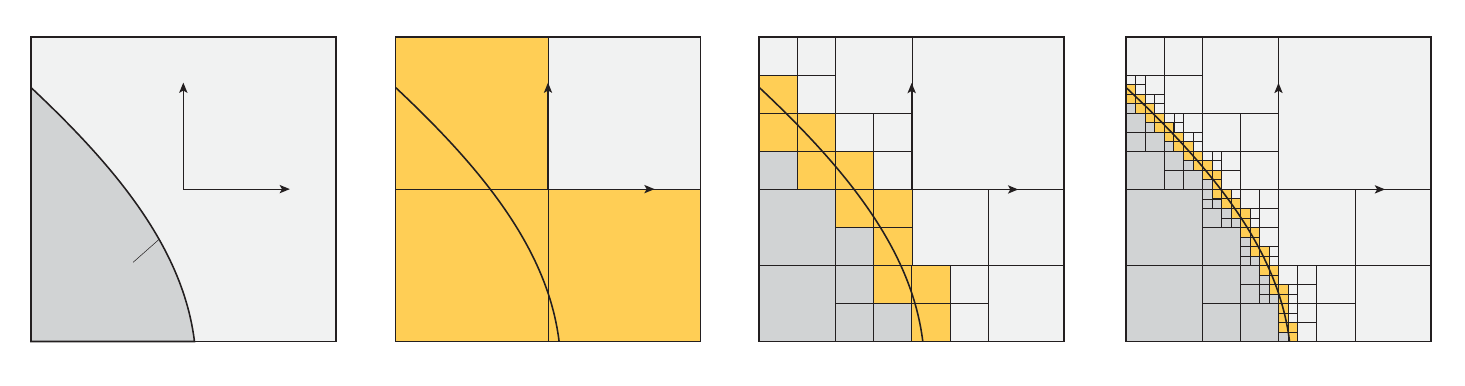}
	\put(20.5,11){$\xi$}
	\put(10.5,21.5){$\eta$}
	\put(45.25,11){$\xi$}
	\put(35.25,21.5){$\eta$}
	\put(70.5,11){$\xi$}
	\put(60,21.5){$\eta$}
	\put(94.75,11){$\xi$}
	\put(85,21.5){$\eta$}
	\put(17,19){$\Omega_\mathrm{fict}$}
	\put(3,4.5){$\Omega_\mathrm{phys}$}
	\put(5.5,7){$\partial\Omega$}
	\put(10.5,0){$k\,{=}\,0$}
	\put(35.5,0){$k\,{=}\,1$}
	\put(60,0){$k\,{=}\,3$}
	\put(85.5,0){$k\,{=}\,5$}
\end{overpic}
\caption{Construction of the element-level integration grid based on the quadtree decomposition of a cut element. Three different refinement levels $k = 1,3,5$ are depicted. Dark gray subcells indicate that the integration domain belongs to the physical domain, while a light grey color refers to the fictitious domain, and cut integration subcells are marked in yellow.}
\label{fig:Quadtree_levels}
\end{figure}
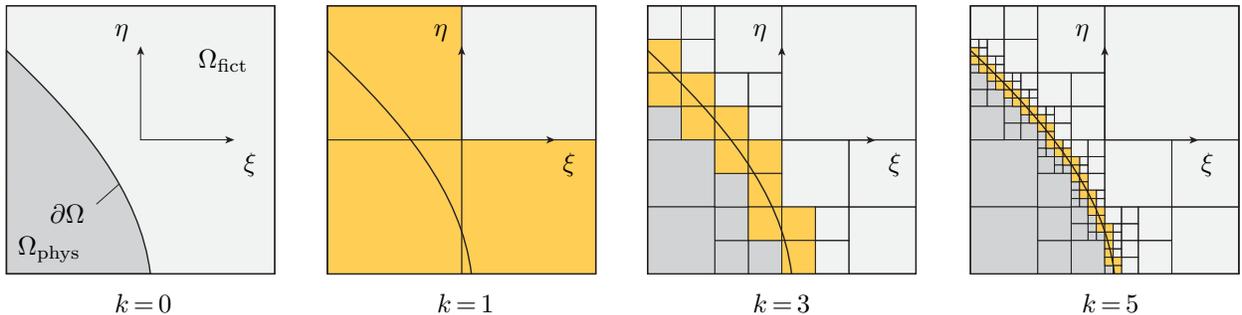

As a consequence, standard Gaussian quadrature rules commonly used in the FEM for numerical integration of the system matrices do not yield accurate results when applied to the discontinuous integrands in immersed boundary methods. To address this issue, tailored solutions have been developed and extensively discussed in Ref.~\cite{ArticlePetoe2020}. Among these solutions, quadtree/octree based approaches have gained popularity due to their robustness and flexibility (as shown in Fig.~\ref{fig:Quadtree_levels}). Additionally, methods based on moment fitting \cite{ArticleMueller2013, ArticleDuester2019, ArticleLegrain2021, ArticleGarhuom2022b} or the divergence theorem \cite{ArticleDuczek2015b} have proven to be more efficient. These specialized integration techniques enable accurate and efficient numerical computations within the fictitious domain, ensuring reliable results in simulations involving complex geometries. Furthermore, the treatment of Neumann and Dirichlet boundary conditions becomes more intricate due to the non-conforming nature of the spatial discretization \cite{InbookDuester2018}.

It is important to highlight that all simulations discussed in this article are conducted using the \emph{SCM}, a specialized version of the FCM specifically designed for explicit dynamics \cite{ArticleDuczek2014}. The SCM is based on the spectral element method (SEM), which utilizes a nodal basis instead of the hierarchical one commonly employed in \emph{p}-FEM. A notable advantage of SEM is its direct generation of a diagonal mass matrix through the nodal quadrature technique, eliminating the need for heuristic methods like row-summing \cite{BookCook1989} or diagonal scaling (HRZ-method) \cite{ArticleHinton1976}. The key element lies in the utilization of Lagrange shape functions defined on a Gau\ss{}-Lobatto-Legendre grid, which mitigates the issues associated with Runge's phenomenon and enables highly accurate results. For more in-depth information on spectral shape functions and the SEM, interested readers are encouraged to refer to the monographs by Pozrikidis \cite{BookPozrikidis2014} and Karniadakis and Sherwin \cite{BookKarniadakis2005}, which offer comprehensive discussions on the topic. Further improvements in accuracy for dynamic problems can be achieved by adopting a higher-order mass formulation, a concept initially pioneered by Goudreau~\cite{PhDGoudreau1970, ArticleGoudreau1973} and later endorsed by Hughes~\cite{BookHughes1987}. This approach involves computing a weighted average of the consistent and lumped mass matrices, which is subsequently employed in the simulation process. Additionally, Ainsworth and Wajid have demonstrated an alternative route to achieve the same accuracy boost -- by developing a tailored integration rule for the mass matrix~\cite{ArticleAinsworth2010}. This innovative approach leads to the creation of optimally blended spectral elements, which have been comprehensively investigated in Ref.~\cite{ArticleRadtke2021}.

It is important to recognize that the advantages of the SEM do not directly carry over to spectral cells. Thus, let us revisit the methodology for achieving a diagonal mass matrix within the SCM framework. In the wide body of literature, two principal approaches have emerged for addressing cut elements:
\begin{enumerate}
	\item The HRZ-method as detailed in Ref.~\cite{ArticleJoulaian2014a}.
	\item Non-negative moment fitting combined with the nodal quadrature technique as detailed in Refs.~\cite{ArticleNicoli2022,ArticleNicoli2023}.
\end{enumerate}
Regardless of the approach chosen, uncut elements (standard spectral elements) are consistently lumped using nodal quadrature, ensuring optimal convergence, as stated in Ref.~\cite{ArticleDuczek2019b}., while the aforementioned approaches are exclusively applied to cut elements. However, it is worth noting that HRZ-lumping of cut elements negatively affects achievable convergence rates. On the other hand, the nodal quadrature technique in conjunction with moment fitting leads to significant under-integration issues of the mass matrix. In the remainder of this section, we will provide a concise overview of the nodal quadrature technique for spectral elements and the HRZ-lumping technique for cut elements, which is the preferred approach in this paper.

At this point, let us briefly outline the fundamentals of mass lumping for both spectral elements and cut elements in a one-dimensional context. It is important to note that this methodology can be readily extended to multi-dimensional problems using a tensor product formulation. Spectral shape functions essentially involve Lagrange interpolation polynomials defined on a non-uniform grid of nodal points, with Gauss-Lobatto-Legendre (GLL) points being a commonly used choice for this purpose, as documented in Refs.~\cite{BookKarniadakis2005, BookPozrikidis2014}. GLL points are determined as the roots of completed Lobatto polynomials, expressed as:
\begin{equation}
	\big( 1- \xi^2 \big) L_{p-1}(\xi) = 0\,.
\end{equation}
The inclusion of the term $\big( 1- \xi^2 \big)$ guarantees that nodes are located at the element boundaries ($\xi = \pm1$), making it suitable for continuous Galerkin formulations. The term $L_{p-1}(\xi)$ represents the first derivative of a Legendre polynomial of order $p$, often referred to as a Lobatto polynomial. The element shape functions are defined as Lagrangian interpolation polynomials supported at GLL nodes, with their mathematical expression as follows:
\begin{equation}
	N_i^p(\xi) = \prod\limits_{j=1,\, j\ne i}^{p+1} \cfrac{\xi - \xi_j^p}{\xi_i^p - \xi_j^p}\,,
\end{equation}
where $\xi_k^p$ denotes the $k^\mathrm{th}$ GLL-point of order $p$. The individual shape functions are then assembled in the matrix of shape functions
\begin{equation}
	\mathbf{N}(\xi) = 
	\begin{bmatrix}
		N_1^p & N_2^p & N_3^p & \ldots & N_{p+1}^p
	\end{bmatrix},
\end{equation}
which is employed to define the consistent mass matrix (CMM) of a spectral element:
\begin{equation}
	\mathbf{M}_\mathrm{e} = \int\limits_{\Omega_\mathrm{e}}^{} \rho \left(\mathbf{N}^\mathrm{T} \mathbf{N}\right) \, \mathrm{d}\Omega\,,
	\label{eq:CMM}
\end{equation}
with a single component being computed as:
\begin{equation}
	{}_\mathrm{CMM}M^{ij}_\mathrm{e} = \int\limits_{\Omega_\mathrm{e}}^{} \rho(x) N_i^p(\xi)N_j^p(\xi)\, \mathrm{d}\Omega = \sum\limits_{k=1}^{p+1}\rho(x(\xi_k^p)) N_i^p(\xi_k^p)N_j^p(\xi_k^p) w_k^p \det{(\mathbf{J})}|_{\xi_k^p} \,.
\end{equation}
Here, $(\xi_k^p,\,w_k^p)$ denote the integration points and weights and $\mathbf{J}$ represents the Jacobi-matrix of the geometric mapping. Given that all shape functions based on Lagrange polynomials satisfy the Kronecker-delta property, it is feasible to diagonalize the mass matrix, as defined in Eq.~\eqref{eq:CMM}, by means of the nodal quadrature technique. That is to say, GLL-points are not only employed to define the shape functions, but also to numerically integrate the mass matrix. This results, by definition, in a lumped (diagonal) mass matrix (LMM), which is high-order convergent \cite{ArticleFried1975, ArticleMalkus1986, BookCook1989}. Due to the simplicity and elegance of this approach, the components of the mass matrix can be efficiently computed using the following expression:
\begin{equation}
	{}_\mathrm{LMM}M^{ii}_\mathrm{e} = \rho \, w_i^p \det{(\mathbf{J})}|_{\xi_i^p}\,.
	\label{eq:LMM_NQ}
\end{equation}
However, for cut elements that require a more sophisticated numerical integration technique to address the discontinuous integrand, achieving lumping solely through nodal quadrature is feasible only if a certain level of under-integration is accepted. In such cases, a moment fitting approach can be adopted to adjust the integration weights, as outlined in Refs.~\cite{ArticleNicoli2022, ArticleNicoli2023}.

In this contribution, an alternative path is pursued, wherein we first compute the consistent mass matrix of a cut element
\begin{equation}
	\mathbf{M}_\mathrm{c} = \int\limits_{\Omega_\mathrm{e}}^{} \alpha_\mathrm{FCM}\rho \left(\mathbf{N}^\mathrm{T} \mathbf{N}\right) \, \mathrm{d}\Omega\,.
	\label{eq:CMM_cut}
\end{equation}
Here, we accurately evaluate the integral form~\eqref{eq:CMM_cut} by means of a quadtree-based numerical integration technique \cite{ArticlePetoe2020}
\begin{equation}
	\begin{split}
		{}_\mathrm{CMM}M^{ij}_\mathrm{c} & = \sum\limits_{s=1}^{n_\mathrm{s}}\,\int\limits_{\Omega_{s}}^{} \alpha_\mathrm{FCM}(x)\rho(x) N_i^p(\xi)N_j^p(\xi)\, \mathrm{d}\Omega \\
		& = \sum\limits_{s=1}^{n_\mathrm{s}}\sum\limits_{k=1}^{p+1}\alpha_\mathrm{FCM}(x(\xi(r_k^p)))\rho(x(\xi(r_k^p))) N_i^p(\xi(r_k^p))N_j^p(\xi(r_k^p)) w_k^p \det{(\mathbf{J}_\mathrm{s})}|_{\xi(r_k^p)} \det{(\mathbf{J})}|_{\xi(r_k^p)} \,.
	\end{split}
	\label{eq:Mij_CMM_cut}
\end{equation}
In Eq.~\eqref{eq:Mij_CMM_cut}, $n_\mathrm{s}$ denotes the number of integration subdomains and $\mathbf{J}_\mathrm{s}$ is the second Jacobi matrix utilized to accommodate the geometric transformation from the subdomain with local coordinate $r$ to the element reference frame with local coordinate $\xi$.

Moreover, we note that this approach yields a fully-populated mass matrix, which is subsequently subjected to diagonalization via the HRZ-lumping scheme, as described in Ref.~\cite{ArticleHinton1976}. In mathematical terms, this method is represented as:
\begin{equation}
	{}_\mathrm{LMM}M^{ii}_c = \beta {}_\mathrm{CMM}M^{ii}_c\,,
	\label{eq:LMM_HRZ}
\end{equation}
where $\beta$ serves as a scaling factor that ensures mass conservation. This scaling factor is defined as the ratio of the total mass $m_\mathrm{c}$ of the cut element divided by the sum of the components on the main diagonal:
\begin{equation}
	\beta = \cfrac{m_\mathrm{c}}{\sum\limits_{i=1}^{p+1}{}_\mathrm{CMM}M^{ii}_c} = \cfrac{\sum\limits_{i=1}^{p+1} \sum\limits_{j=1}^{p+1} {}_\mathrm{CMM}M^{ij}_c}{\sum\limits_{i=1}^{p+1}{}_\mathrm{CMM}M^{ii}_c}\,.
\end{equation}
By employing this two-step approach, we are able to obtain a diagonal mass matrix even within the context of immersed boundary methods utilizing formulations based on nodal shape functions.

Despite the acknowledged limitations regarding mass lumping, the SCM has consistently demonstrated a good performance across various applications, including smart structure analysis \cite{ArticleDuczek2015a}, guided wave propagation analysis for structural health monitoring \cite{ArticleNicoli2022, ArticleNicoli2023}, seismic wave propagation analysis incorporating nonlinear effects \cite{ArticleGiraldo2017}, and the study of heterogeneous materials like sandwich panels with foam cores, utilizing multiple GPUs and CPUs \cite{ArticleMossaiby2019}. These notable examples illustrate the versatility of the SCM in dynamic simulations.
\section{Eigenvalue stabilization technique}
\label{sec:EVST}
In this study, we adopt a stabilization approach based on an eigenvalue decomposition of the system matrices at the element level. This ensures that all operations are performed on individual elements, resulting in minimal computational overhead for most practical applications, particularly in nonlinear and transient analyses.

This section provides a comprehensive discussion of the EVS-technique, which is subdivided into three parts: Firstly, we present its extension to dynamic problems, introducing a novel mass matrix stabilization technique. This approach represents the main innovation in this contribution, as it has not been previously explored. Secondly, building upon the insights from previous studies, which have revealed substantial improvements in terms of the condition number and robustness, particularly in nonlinear problems, we explore the application of the EVS-technique to the stiffness matrix of a cut element, a task that proves to be more intricate when compared to mass stabilization. Thirdly, we put forward an innovative scaling approach that establishes a relationship between the stabilization matrices and the corresponding finite element matrices. This crucial step ensures that the stabilization remains unaffected by the choice of units within the analysis framework.

Following the outlined methodology, we assert that the stabilization of the elemental mass and stiffness matrices holds the potential to extend the critical time step size in explicit dynamics. Additionally, the proposed procedure, being based on an eigenvalue decomposition, offers a more targeted approach to stabilization by only adjusting those modes that directly contribute to ill-conditioning, in contrast to the $\alpha$-method, which applies stabilization across the entire fictitious domain. This refined strategy contributes to a more robust implementation of immersed boundary methods, promising substantial enhancements in terms of stability, performance, and reliability.
\subsection{Mass matrix}
\label{subsec:EVST_mass}
In the context of transient analyses, it is crucial to assess the effectiveness of stabilizing the mass matrix $\mathbf{M}_\mathrm{c}$ of a cut element using the EVS-technique. Since the critical time step size is inversely proportional to the largest eigenvalue of the matrix product $\mathbf{M}_\mathrm{c}^{-1}\mathbf{K}_\mathrm{c}$. a key question arises: Can stabilizing the mass matrix, or both the mass and stiffness matrices, significantly increase the critical time step size $\Delta t_\mathrm{cr}$ for explicit time integration schemes? A positive answer to this question would enhance the efficiency and potentially the accuracy of explicit simulations.

To identify the modes that require stabilization, we calculate all eigenvalues and mode shapes (eigenvectors) of the mass matrix for a cut element by solving the following eigenvalue problem:
\begin{equation}
	\mathbf{M}_\mathrm{c}\, \fatgreek{\psi}_i = \omega_i\, \fatgreek{\psi}_i\quad \;\forall i = 1,\,2,\,\ldots,\,n_\mathrm{DOF}\,.
	\label{eq:EVP_Mc}
\end{equation}

Hence, the derivation of the stabilization technique for the mass matrix of a cut element begins with an eigenvalue decomposition
\begin{equation}
	\mathbf{M}_\mathrm{c} = \fatgreek{\Psi}\, \fatgreek{\Omega}\, \fatgreek{\Psi}^\mathrm{T}\,.
	\label{eq:EigenvalueDecomposition_M}
\end{equation}
Here, the mode shape matrix $\fatgreek{\Psi}$ and the diagonal eigenvalue matrix $\fatgreek{\Omega}$ are introduced. These matrices have dimensions of $[n_\mathrm{DOF} \times n_\mathrm{DOF}]$, with $n_\mathrm{DOF}$ denoting the number of degrees of freedom per finite element. Thus, in the absence of constraints on the mass matrix (no Dirichlet boundary conditions), we have $n_\mathrm{DOF}$ eigenvalue/eigenvector pairs. The individual mode shapes or eigenvectors $\fatgreek{\psi}_i$ corresponding to the eigenvalue $\omega_i$ are stored column-wise in $\fatgreek{\Psi}$.

To ensure the uniqueness of the stabilization technique, it is essential to normalize all eigenvectors such that their Euclidean norms are equal to one:
\begin{equation}
	\lVert\fatgreek{\psi}_i\rVert_{\scriptsize\raisebox{-1.5ex}{$2$}} = 1\,. \label{eq:ModeNormalization}
\end{equation}
Therefore, the following normalization step is performed
\begin{equation}
	\fatgreek{\psi}_i = \cfrac{\fatgreek{\psi}_i}{\lVert\fatgreek{\psi}_i\rVert_{\scriptsize\raisebox{-1.5ex}{$2$}}}\,.
	\label{eq:Normalization}
\end{equation}%

This normalization step guarantees consistency and facilitates meaningful interpretations of the mode shapes during the stabilization process. It is important to note that the individual modes are partitioned into two subspaces:
\begin{equation}
	\mathbf{M}_\mathrm{c} =
	\begin{bmatrix}
		\bar{\fatgreek{\Psi}} & \hat{\fatgreek{\Psi}}
	\end{bmatrix}
	\begin{bmatrix}
		\bar{\fatgreek{\Omega}} & \mathbf{0}\\
		\mathbf{0} & \hat{\fatgreek{\Omega}}
	\end{bmatrix}
	\begin{bmatrix}
		\bar{\fatgreek{\Psi}} & \hat{\fatgreek{\Psi}}
	\end{bmatrix}^{\scriptsize\raisebox{+1.0ex}{$\mathrm{T}$}}.
	\label{eq:EigenvalueDecomposition_M_Partitioned}
\end{equation}%
The quantities with an overbar (e.g., $\bar{\square}$) represent meaningful results from the non-zero eigenspace, while those with a hat (e.g., $\hat{\square}$) pertain to the zero eigenspace. It is important to note that due to the effect of inertia, \emph{rigid body modes (RBMs) are not part of the zero eigenspace of the mass matrix} \cite{ArticleLoehnert2015}; instead, (nearly) zero eigenvalues arise from the intersection of the physical boundary with an element. Consequently, all mode shapes contained in the matrix $\hat{\fatgreek{\Psi}}$ are physically meaningless and require stabilization.

However, in cases involving enriched formulations, an additional extraction procedure (as discussed in the next section for the stiffness matrix) should be used to extract all physically meaningful eigenvectors corresponding to (nearly) zero eigenvalues. In contrast, in the context of immersed boundary methods that solely deal with voids and do not require an enrichment of the ansatz space, no orthogonalization procedure is necessary to extract specific modes from the zero eigenspace of the mass matrix, simplifying the implementation of the mass stabilization technique and reducing computational costs compared to the stiffness stabilization procedure (see Sect.~\ref{subsec:EVST_stiff}).

At this point, it is crucial to define what qualifies as a (nearly) zero eigenvalue, which is vital for matrix partitioning. Such a definition is vital for the purpose of matrix partitioning. In this study, a practical approach based on the ratio $r_i$ between the $i$\textsuperscript{th} eigenvalue and the largest eigenvalue is used to determine whether a mode requires stabilization:
\begin{equation}
	\epsilon_i = \epsilon_\mathrm{S}f^-(r_i - \epsilon_\lambda)\,, \qquad\text{with}\quad f^-(x) = 
	\begin{cases}
		1 & \forall\; x \le 0 \\
		0 & \forall\; x > 0
	\end{cases}
	\qquad\text{and}\quad r_i = \cfrac{\lambda_i}{\lambda_\mathrm{max}}\,.
	\label{eq:StabFac1}
\end{equation}%
Here, the stabilization parameter $\epsilon_\mathrm{S}$ and the threshold value $\epsilon_\lambda$ are user-defined input parameters\footnote{Remark: Based on the recommendations provided in previous studies on the EVS-technique \cite{ArticleLoehnert2014, ArticleBeese2017, ArticleGarhuom2022}, the values for $\epsilon_\lambda$ and $\epsilon_\mathrm{S}$ are selected as follows: the threshold value for determining whether a mode should be stabilized or not is typically chosen in the interval $\epsilon_\lambda \in [10^{-7},\,10^{-3}]$, while the stabilization parameter is chosen in the interval $\epsilon_\mathrm{S} \in [10^{-7},\,10^{-2}]$. These values have been found to ensure reasonable accuracy and good performance in improving the conditioning of the system matrices.}. By taking Eq.~\eqref{eq:StabFac1} into account, the relation

\begin{equation}
	r_i = \cfrac{\lambda_i}{\lambda_\mathrm{max}} < \epsilon_\lambda
	\label{eq:EVratio}
\end{equation}%
is automatically fulfilled. To achieve this, the function $f^-(x)$ is utilized, which takes a value of zero for positive arguments and a value of one for negative arguments. As a result, a stabilization factor of either $\epsilon_\mathrm{S}$ or zero is determined for all modes. Consequently, all $n_\mathrm{s}$ (nearly) zero eigenvalues and their corresponding eigenvectors can be easily collected in $\hat{\fatgreek{\Omega}}$ and $\hat{\fatgreek{\Psi}}$, respectively. It is important to note that the value of the stabilization parameter $\epsilon_i$ determines whether the corresponding mode requires stabilization. Therefore, if $r_i$ is smaller than $\epsilon_\lambda$, Eq.~\eqref{eq:EVratio} ensures that the mode shape $\fatgreek{\psi}_i$ is included in $\hat{\fatgreek{\Psi}}$, and $\omega_i$ is contained in $\hat{\fatgreek{\Omega}}$.

For simplicity, the user-defined threshold value $\epsilon_\lambda$ is identical for both mass and stiffness stabilization approaches. Finally, the mass stabilization matrix $\tilde{\mathbf{M}}_\mathrm{c}^\mathrm{S}$ for a cut element is computed as follows:
\begin{equation}
	\tilde{\mathbf{M}}_\mathrm{c}^\mathrm{S} = \epsilon_\mathrm{S} \sum\limits_{i=1}^{n_\mathrm{s}} \left(\hat{\fatgreek{\psi}}_{i}\,\hat{\fatgreek{\psi}}_{i}^\mathrm{T}\right)\,.
	\label{eq:Stabilization Matrix_M}
\end{equation}
It should be noted that the value of the stabilization factor $\epsilon_i$ provided in Eq.~\eqref{eq:StabFac1} is mode-independent (constant) and, therefore, all mode shapes requiring stabilization are multiplied by the same factor $\epsilon_\mathrm{S}$\footnote{Remark: Previous studies have shown that the stabilization factor $\epsilon_\mathrm{S}$ can be chosen within a relatively large interval without significantly affecting the results. However, a value around $10^{-3}$ has been suggested as a rule of thumb in Ref.~\cite{ArticleGarhuom2022}.}. Other possibilities, where the stabilization parameter is a function of the eigenvalue and/or the volume fraction of the cut element are discussed in Appendix~\ref{App:Eps_S}.

Finally, the contributions from all badly cut elements are assembled in a conventional finite element manner to obtain the overall mass stabilization matrix. This is achieved by summing up the individual stabilization matrices of each cut element, resulting in the expression
\begin{equation}
	\mathbf{M}^\mathrm{S} = \bigcup\limits_{c=1}^{n_\mathrm{C}} \tilde{\mathbf{M}}_c^\mathrm{S}\,.
\end{equation}
Here, $n_\mathrm{C}$ represents the number of cut elements that require stabilization. In the subsequent discussions, all stabilized quantities -- also referred to as modified quantities in the remainder of this article -- will be denoted by a superscript $\square^\mathrm{Mod}$ to indicate that they have undergone the EVS-technique. For example, the modified mass matrix is denoted as $\mathbf{M}^\mathrm{Mod}$
\begin{equation}
	\mathbf{M}^\mathrm{Mod} = \mathbf{M} + \mathbf{M}^\mathrm{S}\,.
\end{equation}
For an alternative and more compact expression, please refer to Appendix~\ref{App:AlternativeExp}.

When considering the stabilization technique for the mass matrix, it is important to address additional factors, particularly in the context of transient analyses and wave propagation simulations. To achieve efficient simulations using explicit time integration schemes, mass lumping plays a crucial role. This raises the question of whether it is advisable to directly compute the mass stabilization matrix of a cut element from the lumped or consistent elemental mass matrix. Considering the consistent mass matrix formulation, the resulting stabilization matrix is also fully populated, necessitating a subsequent mass lumping step to maintain a diagonal mass matrix. Within the framework of the FEM, two commonly employed options for mass lumping are: (i) the row-summing technique, and (ii) the HRZ-method, also known as diagonal scaling \cite{ArticleHinton1976, BookCook1989}. These methods will also be subject to performance testing for the mass stabilization matrix.
\subsection{Stiffness matrix}
\label{subsec:EVST_stiff}
When considering the stiffness matrix, it is worth noting that badly cut elements, characterized by a small volume fraction ($\chi$) within the physical domain, exhibit additional eigenvalues close to zero. These eigenvalues, distinct from those associated with the rigid body modes (RBMs) of the structure, can lead to ill-conditioning and stability issues. This presents challenges for both direct and iterative equation solvers, necessitating the need for stabilization.

However, a critical aspect is to avoid stabilizing the RBMs to preserve physically meaningful results. This adds an extra layer of complexity compared to mass stabilization. To identify the modes that genuinely require stabilization, we calculate all eigenvalues and mode shapes (eigenvectors) of the stiffness matrix for a cut element by solving the following eigenvalue problem~\cite{ArticleLoehnert2014}:
\begin{equation}
	\mathbf{K}_\mathrm{c}\, \fatgreek{\phi}_i = \lambda_i\, \fatgreek{\phi}_i\quad \;\forall i = 1,\,2,\,\ldots,\,n_\mathrm{DOF}\,.
	\label{eq:EVP_Kc}
\end{equation}
Hence, the spectral decomposition of the elemental stiffness matrix can be expressed as
\begin{equation}
	\mathbf{K}_\mathrm{c} = \fatgreek{\Phi}\, \fatgreek{\Lambda}\, \fatgreek{\Phi}^\mathrm{T}\,.
	\label{eq:EigenvalueDecomposition_K}
\end{equation}%
Here, $\fatgreek{\Phi}$ represents the matrix of mode shapes, and $\fatgreek{\Lambda}$ is the diagonal eigenvalue matrix. In the absence of Dirichlet boundary conditions, we have $n_\mathrm{DOF}$ eigenvalue/eigenvector pairs, which include $n_0$ physically meaningful zero eigenvalues.

In the context of immersed boundary methods (without enrichment), these $n_0$ physically meaningful zero eigenvalues correspond to the RBMs of a finite element. For two-dimensional applications, we typically have three RBMs (two translational and one rotational), resulting in $n_0\,{=}\,3$. In three-dimensional problems, there are six RBMs (three translational and three rotational), resulting in $n_0\,{=}\,6$.

However, in the case of badly cut elements, there are additional modes within the zero eigenspace, solely caused by a small value of $\chi$. Therefore, when applying the EVS-technique, it is crucial to distinguish between physically meaningful (nearly) zero eigenvalues and those arising from cutting a finite element. The goal is to stabilize all singular modes while preserving the RBMs to avoid unphysical outcomes. To achieve this, we introduce the following partitioning of the matrix of mode shapes and the eigenvalue matrix:
\begin{equation}
	\mathbf{K}_\mathrm{c} =
	\begin{bmatrix}
		\bar{\fatgreek{\Phi}} & \hat{\fatgreek{\Phi}}
	\end{bmatrix}
	\begin{bmatrix}
		\bar{\fatgreek{\Lambda}} & \mathbf{0}\\
		\mathbf{0} & \hat{\fatgreek{\Lambda}}
	\end{bmatrix}
	\begin{bmatrix}
		\bar{\fatgreek{\Phi}} & \hat{\fatgreek{\Phi}}
	\end{bmatrix}^{\scriptsize\raisebox{+1.0ex}{$\mathrm{T}$}}.
	\label{eq:EigenvalueDecomposition_K_Partitioned}
\end{equation}%
Here, variables denoted with an overbar ($\bar{\square}$) refer to nonzero quantities, while a hat above a variable ($\hat{\square}$) signifies (nearly) zero eigenvalues and their corresponding eigenvectors. It is important to emphasize that, without loss of generality, we maintain the normalization of all mode shapes, as previously discussed in Eq.~\eqref{eq:ModeNormalization}.

To identify the modes requiring stabilization, which are collected in the matrix $\hat{\fatgreek{\Phi}}$, the same condition as given in Eq.~\eqref{eq:EVratio}
is utilized. For simplicity, the user-defined threshold value $\epsilon_\lambda$ is identical for both the stiffness and mass stabilization approaches. However, it is vital to keep in mind that the collected modes still include the physically meaningful zero eigenspace associated with the RBMs of the structure, which must remain unchanged. In the next step of the analysis, it becomes necessary to extract these mode shapes. Consequently, we further partition the matrices $\hat{\fatgreek{\Phi}}$ and $\hat{\fatgreek{\Lambda}}$
\begin{equation}
	\hat{\fatgreek{\Phi}} = 
	\begin{bmatrix}
		\hat{\fatgreek{\Phi}}_0 & \hat{\fatgreek{\Phi}}_\mathrm{u}
	\end{bmatrix}
	\quad\text{and}\quad
	\hat{\fatgreek{\Lambda}} = 
	\begin{bmatrix}
		\hat{\fatgreek{\Lambda}}_0 & \mathbf{0}\\
		\mathbf{0} & \hat{\fatgreek{\Lambda}}_\mathrm{u}
	\end{bmatrix}\,,
	\qquad\text{with}\quad
	\hat{\fatgreek{\Lambda}}_0 = \mathbf{0}\,.
	\label{eq:SplitPhiLambda}
\end{equation}%
In this context, the matrix $\hat{\fatgreek{\Phi}}_\mathrm{u}$ contains all $n_\mathrm{u}$ mode shapes corresponding to (nearly) zero eigenvalues, while $\hat{\fatgreek{\Phi}}_0$ contains the RBMs. It is worth noting that in Eq.~\eqref{eq:SplitPhiLambda}, the eigenvalue matrix $\hat{\fatgreek{\Lambda}}_0$ is  is theoretically defined as a zero matrix. However, in practical numerical computations, round-off errors introduce very small but non-zero eigenvalues.

Moreover, it is important to consider that in enriched methods like XFEM, additional modes beyond the RBMs must be included in $\hat{\fatgreek{\Phi}}_0$, as discussed in Refs.~\cite{ArticleLoehnert2014, ArticleBeese2017}. In the context of nonlinear problems, zero eigenvalues associated with stability issues such as buckling or material instabilities must also be accounted for, adding an additional layer of complexity to the identification of the valid zero eigenspace.

Depending on the clustering of (nearly) zero eigenvalues, an extraction procedure is required to isolate the RBMs from the unphysical modes. To achieve this, we employ a Gram-Schmidt orthogonalization procedure
\begin{equation}
	\hat{\fatgreek{\phi}}_i^\mathrm{GS} = \hat{\fatgreek{\phi}}_i - \sum\limits_{j=1}^{n_0} \left( \hat{\fatgreek{\phi}}_i^\mathrm{T}\,\fatgreek{\phi}_{\mathrm{RBM},j} \right) \fatgreek{\phi}_{\mathrm{RBM},j}\qquad \forall \; i \in [1,n_\mathrm{s}]\qquad\text{with}\quad n_\mathrm{s} = n_0 + n_\mathrm{u}\,.
	\label{eq:GramSchmidt}
\end{equation}%
Here, the RBMs -- $\fatgreek{\phi}_{\mathrm{RBM},j}$ -- are extracted from the singular set. By applying Eq.~\eqref{eq:GramSchmidt}, we obtain a new set of eigenvectors denoted as $\hat{\fatgreek{\Phi}}^\mathrm{GS}$, which consists of the orthogonalized unphysical modes. 

In order to identify the physically meaningful zero eigenvalue modes, the vector-norm of each eigenvector $\hat{\fatgreek{\phi}}_i^\mathrm{GS}$ can be computed. Typically, RBMs can be distinguished by their significantly smaller vector-norms compared to the other modes:
\begin{equation}
	\left\lVert\hat{\fatgreek{\phi}}_i^\mathrm{GS}\right\rVert_{\scriptsize\raisebox{-2.5ex}{$2$}} < 1\times 10^{-3}\qquad \forall \; i \in [1,n_\mathrm{s}]\,.
	\label{eq:ExtractRBM}
\end{equation}%
The norm should also be significantly smaller than one since mode shapes are orthogonal to each other. Therefore, only modes that are primarily composed of a linear combination of RBMs are affected by the orthogonalization process. If condition~\eqref{eq:ExtractRBM} is satisfied, indicating that the norm of the eigenvector is sufficiently small, the corresponding eigenvector is deleted from the set $\hat{\fatgreek{\Phi}}^\mathrm{GS}$.

However, for badly cut elements with only a small volume fraction of the physical domain $\chi$, it is sometimes observed that there are no pure RBMs that can be deleted after the orthogonalization step. In these cases, all modes corresponding to (nearly) zero eigenvalues can be described as a linear combination of RBMs and spurious mode shapes\footnote{Remark: The term ``spurious mode shapes'' refers to modes that arise due to the small support of certain degrees of freedom. These modes negatively impact the conditioning of the problem and result in inaccurate approximations, hence requiring stabilization.}, which negatively impact problem conditioning and accuracy, thus requiring stabilization. The orthogonalization procedure extracts the RBM components from the original modes, leading to a notable reduction in the eigenvector's norm. In such cases, condition~\eqref{eq:ExtractRBM} cannot be satisfied, as the components of the mode shapes are still non-zero due to significant contributions from the spurious eigenvectors.

An alternative criterion for checking the presence of RBMs is based on the absolute value of the smallest eigenvalue, $\lambda_\mathrm{min}$. If the following condition is met, the corresponding mode can be deleted:
\begin{equation}
	\left(\hat{\fatgreek{\phi}}_i^\mathrm{GS}\right)^{\scriptsize\raisebox{+1.5ex}{$\mathrm{T}$}} \; \hat{\fatgreek{\phi}}_i^\mathrm{GS} < \lambda_\mathrm{min}\,.
	\label{eq:ExtractRBM_v2}
\end{equation}
The second criterion, based on the numerical properties of the system, appears to be more versatile and will be employed in all simulations throughout the rest of this article. It should be noted that in cases where well-separated RBMs exist, the norm of all eigenvectors in $\hat{\fatgreek{\Phi}}_0$ theoretically becomes zero after the orthogonalization step. The resulting modified set, denoted as $\hat{\fatgreek{\Phi}}^\mathrm{GS}_\mathrm{u}$, contains only the modes used in the stabilization technique. Following the orthogonalization procedure, a normalization step according to Eq.~\eqref{eq:Normalization} is performed. This step ensures that the vector-norm remains equal to one.

To apply the Gram-Schmidt algorithm, as described in Eq.~\eqref{eq:GramSchmidt}, we need analytical expressions describing all RBMs $\fatgreek{\phi}_{\mathrm{RBM},j}$ with $j\in\,[1,n_0]$. Fortunately, such expressions are well-documented in the literature \cite{ArticleJoensthoevel_2013, ArticleBaggio2017}. For the sake of completeness, we provide the expressions for the three RBMs in two-dimensional problems:
\begin{alignat}{2}
	\fatgreek{\phi}_{\mathrm{RBM},1} & = &&
	\begin{bmatrix}
		1 & 0 & 1 & 0 & 1 & \ldots & 0
	\end{bmatrix}^\mathrm{T}\,, \\
	\fatgreek{\phi}_{\mathrm{RBM},2} & = && 
	\begin{bmatrix}
		0 & 1 & 0 & 1 & 0 & \ldots & 1
	\end{bmatrix}^\mathrm{T}\,,\\
	\fatgreek{\phi}_{\mathrm{RBM},3} & = && 
	\begin{bmatrix}
		+R_1 \cos(\theta_1) & -R_1 \sin(\theta_1) & +R_2 \cos(\theta_2) & -R_2 \sin(\theta_2) & +R_3 \cos(\theta_3) & \ldots & -R_{n_\mathrm{N}} \sin(\theta_{n_\mathrm{N}})
	\end{bmatrix}^\mathrm{T}\,.
	\label{eq:RBM}
\end{alignat}%
Here, $R_i$ represents the distance from a specific node to the origin of a element-specific coordinate system ($\tilde{x}\tilde{y}$), and $\theta_i$ denotes the angle with respect to the local $\tilde{x}$-axis. This formulation utilizes a polar coordinate system to derive the expression for the rotational RBM. In three-dimensional cases, a spherical coordinate system is employed \cite{ArticleJoensthoevel_2013, ArticleBaggio2017}. It is assumed that the coordinate vector of an element is defined as follows:
\begin{equation}
	\mathbf{x} = \{x_1\;y_1\;x_2\;y_2\;\cdots\;x_{n_\mathrm{N}}\;y_{n_\mathrm{N}}\}^\mathrm{T}\,,
\end{equation}%
and thus, the radius and the angle are defined as
\begin{equation}
	R_i = \sqrt{\tilde{x}_i^2 + \tilde{y}_i^2}\quad\text{and}\quad \theta_i = \arctan{\left(\cfrac{\tilde{y}_i}{\tilde{x}_i}\right)}\,,
\end{equation}%
with $\tilde{x}_i$ and $\tilde{y}_i$ representing the nodal coordinates of the $i$\textsuperscript{th} node in the elemental coordinate system. In this coordinate system, the origin coincides with the centroid of the element with straight edges. It is important to note that only the corner nodes are considered to determine the location of the centroid
\begin{equation}
	x_c = \cfrac{1}{4}\sum\limits_{i=1}^{4}x_i \quad \text{and} \quad y_c = \cfrac{1}{4}\sum\limits_{i=1}^{4}y_i\,.
\end{equation}
Therefore, the new coordinates in the element-specific coordinate system are
\begin{equation}
	\begin{bmatrix}
		\tilde{x} \\
		\tilde{y}
	\end{bmatrix}
	=
	\begin{bmatrix}
		x \\
		y
	\end{bmatrix}
	-
	\begin{bmatrix}
		x_c \\
		y_c
	\end{bmatrix}
	\qquad \text{or} \qquad
	\tilde{\mathbf{x}} = \mathbf{x} - \mathbf{x}_c\,.
\end{equation}
That means, by knowing the global coordinates of all nodes associated with a finite element, we can easily compute the eigenvectors for all RBMs. The number of nodes involved in the computation depends on the chosen polynomial degree $p$ of the shape functions
\begin{equation}
	n_\mathrm{N} = (p + 1)^d\,,
\end{equation}%
where $d$ denotes the dimensionality of the problem. The expression provided assumes a tensor product formulation of the finite element shape functions, which restricts the discussion to quadrilateral and hexahedral elements in this article. This choice is consistent with most publications in the context of immersed boundary methods, where regular Cartesian meshes are commonly used. For applications involving unstructured discretizations, we refer the reader to Refs.~\cite{ArticleVarduhn2015, ArticleDuczek2016a, ArticleDuczek2016b}. Additionally, we need to consider hierarchic or modal shape functions typically used in the $p$-version of FEM and the FCM~\cite{BookSzabo1991, InbookDuester2018}. In such cases, the higher-order degrees of freedom (DOFs) do not have a direct physical interpretation and can be regarded as unknowns of the high-order polynomial ansatz space. Consequently, when applying rigid body displacements to an element, only the nodal shape functions associated with the corner nodes of the element are activated. This implies that all DOFs associated with high-order shape functions are set to zero in Eq.~\eqref{eq:RBM}. For quadrilateral elements, this results in only eight non-zero components (four nodes with two DOFs each), while for hexahedral elements, 24 components may have non-zero values (eight nodes with three DOFs each). According to Ref.~\cite{ArticleLoehnert2014}, a second orthogonalization step might be necessary as the extraction of the RBMs from $\hat{\fatgreek{\Phi}}$ can lead to an $n_0$-fold linear dependence. Therefore, the mode shapes used to construct the stabilization matrix are subject to additional orthogonalization
\begin{equation}
	\hat{\fatgreek{\phi}}_{\mathrm{u},i}^\mathrm{2GS} = \hat{\fatgreek{\phi}}_{\mathrm{u},i}^\mathrm{GS} - \sum\limits_{j=1}^{i-1} \left( \hat{\fatgreek{\phi}}_{\mathrm{u},i}^\mathrm{GS,T}\,\hat{\fatgreek{\phi}}_{\mathrm{u},j}^\mathrm{GS} \right) \hat{\fatgreek{\phi}}_{\mathrm{u},j}^\mathrm{GS}\qquad \forall\; i \in [1,n_\mathrm{u}]\,.
\end{equation}%
After the second orthogonalization step, an additional normalization step is required. Equation~\eqref{eq:Normalization} can be used once again for this purpose.

Finally, the stiffness stabilization matrix $\tilde{\mathbf{K}}_\mathrm{c}^\mathrm{S}$ for a cut element is computed as follows:
\begin{equation}
	\tilde{\mathbf{K}}_\mathrm{c}^\mathrm{S} = \epsilon_\mathrm{S}\sum\limits_{i=1}^{n_\mathrm{u}} \left(\hat{\fatgreek{\phi}}_{\mathrm{u},i}^\mathrm{2GS}\,\hat{\fatgreek{\phi}}_{\mathrm{u},i}^\mathrm{2GS,T}\right)\,.
	\label{eq:Stabilization Matrix1}
\end{equation}

The global stiffness stabilization matrix is again obtained by aggregating the contributions from all poorly cut elements
\begin{equation}
	\mathbf{K}^\mathrm{S} = \bigcup\limits_{c=1}^{n_\mathrm{C}} \tilde{\mathbf{K}}_c^\mathrm{S}\,,
\end{equation}
and the stabilized (or modified) stiffness matrix is defined as
\begin{equation}
	\mathbf{K}^\mathrm{Mod} = \mathbf{K} + \mathbf{K}^\mathrm{S}\,.
\end{equation}
To verify the correctness of the implementation of the stabilization procedure, it is advisable to perform an eigenvalue analysis of the stabilized stiffness matrix. In this analysis, only $n_0$ zero eigenvalues corresponding to the RBMs of the element should be present. Additionally, the condition number of the stabilized element subject to minimal Dirichlet boundary conditions can be computed and compared to the original value. If the condition number is significantly reduced, it indicates that the implementation is functioning correctly.
\subsection{Scaling approach}
\label{subsec:Scaling}
A severe shortcoming of the stabilization introduced in Sects.~\ref{subsec:EVST_mass} and \ref{subsec:EVST_stiff} is that its parameter $\epsilon_\mathrm{S}$ is independent of the system of units being employed in the analysis. In other words, while the mode shapes remain unaffected by a change in units, the absolute values of the matrix components do vary with them. Because of this reason, we propose to augment the stabilization factor by an additional scaling procedure, which should be a function of the largest component of the stiffness/mass matrix of an uncut (finite) element. It is important to note that, for the remainder of this section, all explanations are related to the stiffness term, but the same technique is also applied to the mass term.

To achieve our goal, the following approach is suggested:
\begin{enumerate}
	\item Determine the maximum absolute value of all components in the elemental stiffness matrix $\mathbf{K}_\mathrm{e}$\footnote{Remark: We explicitly use the finite element matrix $\mathbf{K}_\mathrm{e}$, which is entirely located in the physical domain, i.e., the element is uncut. The rationale behind that decision is that the contribution of one finite element in a sufficiently refined mesh is already small, such that any additional dependence of the stabilization magnitude on the volume fraction $\chi$ of a cut element should be avoided a priori.}:
	\begin{equation}
		k_\mathrm{e,max} = \max[(\mathbf{K}_\mathrm{e})_{ij}]  \qquad \forall \; i,j = 1,2,\ldots n_\mathrm{DOF}
	\end{equation}
	\item Determine the maximum absolute value of all components in the elemental stiffness stabilization matrix $\tilde{\mathbf{K}}_\mathrm{c}^\mathrm{S}$:\\
	\begin{equation}
		k_\mathrm{c,max}^\mathrm{S} = \max[(\tilde{\mathbf{K}}_\mathrm{c}^\mathrm{S})_{kl}] \qquad \forall \; k,l = 1,2,\ldots n_\mathrm{DOF}
	\end{equation}
	\item Compute the so-called scaling parameter $n_\alpha$:
	\begin{equation}
		n^\mathrm{k}_\alpha = 10^{\gamma}\,, \qquad\text{with}\quad \gamma = \left \lfloor \log_{10}\left(\cfrac{k_\mathrm{e,max}}{k_\mathrm{c,max}^\mathrm{S}} \; \epsilon_\mathrm{S}\right) \right \rceil \,.
		\label{eq:Alpha}
	\end{equation}%
	\item Scale the stiffness stabilization matrix:
	\begin{equation}
		\mathbf{K}_\mathrm{c}^\mathrm{S} = n^\mathrm{k}_\alpha \tilde{\mathbf{K}}_\mathrm{c}^\mathrm{S}\,.
	\end{equation}
\end{enumerate}
This approach ensures that for a specific stabilization parameter $\epsilon_\mathrm{S}$, a similar degree of stabilization is achieved regardless of material properties. Furthermore, Eq.~\eqref{eq:Alpha} ensures that the largest component in the stabilization matrix is approximately $\epsilon_\mathrm{S}^{-1}$ times smaller than the largest entry in the corresponding finite element matrix, where $\lfloor \square \rceil$ denotes a rounding operation. At this point, it is essential to emphasize once more that we handle the mass stabilization matrix in an analogous fashion. To this end, we simply substitute all stiffness-related terms in the previously outlined process with their mass-related counterparts. Consequently, we will utilize distinct scaling parameters denoted as $n^\mathrm{k}_\alpha$ and $n^\mathrm{m}_\alpha$ for stiffness and mass stabilization, respectively.

It is worth emphasizing that our strategy for determining the actual stabilization parameter differs significantly from previous works \cite{ArticleLoehnert2014, ArticleBeese2017, ArticleGarhuom2022}. Instead of applying a nonlinear function to determine the actual stabilization factor, we opt for a constant one, which is scaled based on the maximum components of the stabilization matrix and the corresponding finite element matrix. Therefore, all recommendations regarding a suitable value for $\epsilon_\mathrm{S}$ are not directly transferable to our case. They offer guidance to initiate our investigations, but we expect that a different range of values will emerge as the optimal choice.

While we can intuitively understand the logic behind using a nonlinear function to stabilize modes with smaller eigenvalues by using larger stabilization parameters, there is no concise theoretical justification or numerical evidence available to support the notion that this is the best choice. Hence, we choose the simplest and most elegant approach, i.e., the one that has been established in previous sections. However, it is certainly worthwhile to refine the formulation of the stabilization parameter based on a simple benchmark model and multi-objective optimization. For transient problems, which are the primary focus of this contribution, the main goal would be to maximize the critical time increment, while maintaining a certain level of accuracy.
\subsection{Variants of the eigenvalue stabilization scheme}
\label{subsec:Stab_variants}
In the original implementation of the FCM or SCM, the $\alpha$-method was utilized to stabilize the results in the fictitious domain. Therefore, the $\alpha$-method is chosen as a reference to evaluate the performance of different variants of the EVS-technique. The EVS-technique can be broadly classified into three main categories based on which matrices are stabilized:
\begin{enumerate}
\item Stiffness-Stabilized Systems,
\item Mass-Stabilized Systems,
\end{enumerate}
and
\begin{enumerate}
	\item[3.] Mass-Stiffness-Stabilized Systems.
\end{enumerate}
\begin{table}[b!]
	\caption{Variants of the eigenvalue stabilization technique.\label{tab:Variants}}
	\centering
	\begin{tabular}{rcccccccccc}
		\toprule[2pt]
		& No. & Var. & \phantom{xx} & $\mathbf{K}_c^\mathrm{Mod}$ & $\mathbf{M}_c^\mathrm{Mod}$ & Compute $\fatgreek{\Psi}$ and $\fatgreek{\Omega}$ & Lumping: $\mathbf{M}^\mathrm{S}_c$ & Lumping: $\mathbf{M}_c$ & $\alpha_0$ & \phantom{x}\\
		\midrule[1pt]
		& 1  & 0a && \XSolidBrush & \XSolidBrush & n/a & n/a          & \Checkmark   & $0$ &\\
		& 2  & 0b && \XSolidBrush & \XSolidBrush & n/a & n/a          & \XSolidBrush & $0$ &\\
		& 3  & 0c && \XSolidBrush & \XSolidBrush & n/a & n/a          & \Checkmark   & $10^{-10}$ &\\
		& 4  & 0d && \XSolidBrush & \XSolidBrush & n/a & n/a          & \XSolidBrush & $10^{-10}$ &\\
		\rowcolor{Matlab6}
		Reference: & 4 & 0e && \XSolidBrush & \XSolidBrush & n/a & n/a          & \Checkmark   & $10^{-5}$ &\\
		& 6  & 0f && \XSolidBrush & \XSolidBrush & n/a & n/a          & \XSolidBrush & $10^{-5}$ &\\\hline\\[-1.75ex]
		& 7  & 1a && \Checkmark   & \XSolidBrush & n/a & n/a          & \Checkmark   & 0         &\\
		& 8  & 1b && \Checkmark   & \XSolidBrush & n/a & n/a          & \XSolidBrush & 0         &\\\hline\\[-1.75ex]
		& 9  & 2a && \XSolidBrush & \Checkmark   & CMM & \XSolidBrush & \Checkmark   & 0         &\\
		\rowcolor{Matlab2}
		Optimal: & 10 & 2b && \XSolidBrush & \Checkmark   & CMM & HRZ          & \Checkmark   & 0         &\\
		& 11 & 2c && \XSolidBrush & \Checkmark   & CMM & Row-Sum      & \Checkmark   & 0         &\\
		& 12 & 2d && \XSolidBrush & \Checkmark   & CMM & \XSolidBrush & \XSolidBrush & 0         &\\
		& 13 & 2e && \XSolidBrush & \Checkmark   & CMM & HRZ          & \XSolidBrush & 0         &\\
		& 14 & 2f && \XSolidBrush & \Checkmark   & CMM & Row-Sum      & \XSolidBrush & 0         &\\\hdashline\\[-1.75ex]
		& 15 & 2g && \XSolidBrush & \Checkmark   & LMM & \XSolidBrush & \Checkmark   & 0         &\\
		& 16 & 2h && \XSolidBrush & \Checkmark   & LMM & HRZ          & \Checkmark   & 0         &\\
		& 17 & 2i && \XSolidBrush & \Checkmark   & LMM & Row-Sum      & \Checkmark   & 0         &\\
		\rowcolor{Matlab5}
		Example: & 18 & 2j && \XSolidBrush & \Checkmark   & LMM & \XSolidBrush & \XSolidBrush & 0         &\\
		& 19 & 2k && \XSolidBrush & \Checkmark   & LMM & HRZ          & \XSolidBrush & 0         &\\
		& 20 & 2l && \XSolidBrush & \Checkmark   & LMM & Row-Sum      & \XSolidBrush & 0         &\\\hline\\[-1.75ex]
		& 21 & 3a && \Checkmark   & \Checkmark   & CMM & \XSolidBrush & \Checkmark   & 0         &\\
		& 22 & 3b && \Checkmark   & \Checkmark   & CMM & HRZ          & \Checkmark   & 0         &\\
		& 23 & 3c && \Checkmark   & \Checkmark   & CMM & Row-Sum      & \Checkmark   & 0         &\\
		& 24 & 3d && \Checkmark   & \Checkmark   & CMM & \XSolidBrush & \XSolidBrush & 0         &\\
		& 25 & 3e && \Checkmark   & \Checkmark   & CMM & HRZ          & \XSolidBrush & 0         &\\
		& 26 & 3f && \Checkmark   & \Checkmark   & CMM & Row-Sum      & \XSolidBrush & 0         &\\\hdashline\\[-1.75ex]
		& 27 & 3g && \Checkmark   & \Checkmark   & LMM & \XSolidBrush & \Checkmark   & 0         &\\
		& 28 & 3h && \Checkmark   & \Checkmark   & LMM & HRZ          & \Checkmark   & 0         &\\
		& 29 & 3i && \Checkmark   & \Checkmark   & LMM & Row-Sum      & \Checkmark   & 0         &\\
		& 30 & 3j && \Checkmark   & \Checkmark   & LMM & \XSolidBrush & \XSolidBrush & 0         &\\
		& 31 & 3k && \Checkmark   & \Checkmark   & LMM & HRZ          & \XSolidBrush & 0         &\\
		& 32 & 3l && \Checkmark   & \Checkmark   & LMM & Row-Sum      & \XSolidBrush & 0         &\\
		\bottomrule[2pt]
	\end{tabular}
\end{table}%
Moreover, various subcategories arise from different choices, such as the order of lumping the mass matrix before or after computing the eigenvalues, the specific lumping technique used, and more. For ease of reference, all these choices are summarized in Table~\ref{tab:Variants}, which also includes an extra column for the value of $\alpha_0$. Setting the indicator function to zero within the fictitious domain allows us to focus solely on the impact of the EVS-technique. However, it is worth noting that combining the $\alpha$-method with the EVS-technique (targeted stabilization) is indeed a valid option, and this approach will be further investigated in this study. On the other-hand side, it is obvious that the application of $\alpha$-stabilization to cut cells diminishes the effectiveness of the EVS-technique since the entire fictitious domain is already stabilized. As a result, fewer modes of a cut element need stabilization, leading to a partial loss of the actual benefits of the EVS-technique. Consequently, combining both stabilization schemes does not guarantee an improved performance.

The presence of numerous submethods is related to the mass stabilization matrix, which can be computed using either the consistent mass matrix (CMM) or the lumped mass matrix (LMM). In this article, the HRZ method is employed for lumping the mass matrix within the framework of the SCM, as discussed in Ref.~\cite{ArticleJoulaian2014a}. Other options such as (constrained) moment fitting approaches \cite{ArticleNicoli2022} are out of the scope of the current contribution. However, for lumping the mass stabilization matrix two different options exist: (i) the HRZ-method and (ii) row-summing. Clearly, another option is to use the consistent matrix for stabilization. However, this approach has drawbacks, especially when dealing with explicit time integration schemes. For the sake of completeness and despite the aforementioned issues, we have still included all variants that employ either a consistent mass matrix or a consistent mass stabilization matrix in our list. Nevertheless, it is important to note that these approaches are not competitive in terms of the overall numerical costs associated with an explicit time integration scheme. In the following, our investigations first focus on the potential increase in the critical time increment. However, it is essential to consider the overall efficiency of the developed approach in the later stages of our analysis (see Sect.~\ref{sec:Examples}).

To explain the notation put forward in Table~\ref{tab:Variants}, let us consider \textcolor{Matlab5}{\emph{variant 2j}} as an example. In this case, only the mass matrix is stabilized (\ding{51}), while the original matrix is used for the stiffness matrix (\ding{55}). The mass stabilization matrix, denoted as $\mathbf{M}^\mathrm{S}$, is computed based on the lumped (diagonalized) mass matrix (LMM), i.e., $\fatgreek{\Psi}$ and $\fatgreek{\Omega}$ are derived from the LMM formulation. However, neither the cell's mass matrix $\mathbf{M}_c$ (\ding{55}) nor the mass stabilization matrix $\mathbf{M}^\mathrm{S}_c$ (\ding{55}) are lumped in this variant of the EVS-technique. The other rows in this table marking other variants of the stabilization scheme are to be understood/interpreted in a similar manner.
\section{Central difference method}
\label{sec:StabTimeInt}
This section explores the ramifications of stabilizing the mass matrix and/or stiffness matrix on the chosen time integration scheme. Specifically, we focus on the central difference method (CDM) as a representative example of explicit time integration methods \cite{BookChopra2011}. The observed characteristics and performance can shed light on the broader implications of stabilizing the matrices in explicit time stepping algorithms.
\subsection{Theoretical derivation}
The derivation of temporal discretization schemes begins with the semi-discrete equations of motion, as expressed in Eq.~\eqref{eq:EquationMotion}. In order to solve these second-order ordinary differential equations in time, initial conditions must be specified. These initial conditions can be written as:
\begin{equation}
	\mathbf{U}(t=0) = \mathbf{U}_0\quad\text{and}\quad \dot{\mathbf{U}}(t=0) = \dot{\mathbf{U}}_0\,,
	\label{eq:IC}
\end{equation}
where $\mathbf{U}(t=0)$ and $\dot{\mathbf{U}}(t=0)$ represent the nodal displacement and velocity vectors at the start of the computation, typically at $t\,{=}\,0$. These conditions are necessary to form a well-posed problem for the subsequent derivation of explicit time integration methods. Only by incorporating the given initial conditions, the solution of the equations of motion becomes feasible.

To derive \emph{explicit} time integrators, the balance of momentum equation is expressed at a specific time step $i$, which can be written as follows:
\begin{equation}
	\mathbf{M}\ddot{\mathbf{U}}_i + \mathbf{C}\dot{\mathbf{U}}_i + \mathbf{K}\mathbf{U}_i = \mathbf{F}^i_\mathrm{ex}\,.
	\label{eq:EquationMotion_i}
\end{equation}
In the next step, the temporal derivatives are approximated by central difference formulae
\begin{alignat}{2}
	\dot{\mathbf{U}}_i & = && \cfrac{\mathbf{U}_{i+1} - \mathbf{U}_{i-1}}{2\Delta t}\,, \label{eq:FirstDerivative}\\
	\shortintertext{and} & && \nonumber \\[-18pt]
	\ddot{\mathbf{U}}_i & = && \cfrac{\mathbf{U}_{i+1} - 2\mathbf{U}_{i} + \mathbf{U}_{i-1}}{\Delta t^2}\,, \label{eq:SecondDerivative}
\end{alignat}
which are substituted into Eq.~\eqref{eq:EquationMotion_i}. Subsequently, the terms associated with quantities at time steps $i\,{-}\,1$, $i$, and $i\,{+}\,1$ (i.e., the values we intend to calculate) are gathered independently, resulting in
\begin{equation}
	\left(\cfrac{1}{\Delta t^2}\mathbf{M} + \cfrac{1}{2\Delta t}\mathbf{C}\right)\mathbf{U}_{i+1} = \mathbf{F}^i_\mathrm{ex} - \left[\left(\mathbf{K} - \cfrac{2}{\Delta t^2}\mathbf{M}\right)\mathbf{U}_i + \left(\cfrac{1}{\Delta t^2}\mathbf{M} - \cfrac{1}{2\Delta t}\mathbf{C}\right)\mathbf{U}_{i-1}\right]\,.
	\label{eq:CDM}
\end{equation}
For a comprehensive derivation of the expressions and an extensive analysis of their numerical characteristics, we recommend consulting the monographs by Chopra~\cite{BookChopra2011} or Bathe~\cite{BookBathe2002}.

It is worth noting that the derived expression holds true for every time step and can be reformulated to resemble the static equilibrium by introducing an effective stiffness matrix $\tilde{\mathbf{K}}$ and an effective load vector $\tilde{\mathbf{F}}$
\begin{equation}
	\tilde{\mathbf{K}}\mathbf{U}_{i+1} = \tilde{\mathbf{F}}_i\,.
	\label{eq:StaticEquilibrium_i}
\end{equation}
When employing explicit time stepping schemes in linear analyses, it is crucial to consider two key factors: 
\begin{enumerate}
	\item Temporal progression relies exclusively on the numerical results from preceding time steps, such as $i$ and $i\,{-}\,1$.
	\item Efficient matrix-vector products suffice for advancing in time.
\end{enumerate}
It is worth noting that the second point is only realized when both the mass and damping matrices are diagonal, eliminating the need to solve a system of equations. This aspect underscores the significance of mass lumping within the realm of explicit time stepping investigations.

In explicit methods, the assembly of a global stiffness matrix is typically unnecessary, as only the vector of internal forces is required. Consequently, it is common practice to employ element-by-element techniques. Given that explicit methods are typically conditionally stable, solving a system of equations at each time step would be prohibitively costly.

While it is a viable option in linear elastodynamics to factorize the effective stiffness matrix to reduce the computational burden, the numerical costs are often still dominated by the stability limit. For a comprehensive exploration of mass lumping across different finite element families, we refer the reader to Refs.~\cite{ArticleDuczek2019a, ArticleDuczek2019b} and the related references therein.

Keep in mind that the CDM is not self-starting, necessitating some remarks regarding the start-up procedure. By substituting the initial conditions stated in Eq.~\eqref{eq:IC}, into the semi-discrete equation of motion~\eqref{eq:EquationMotion_i}, the initial acceleration can be derived as follows:
\begin{equation}
	\ddot{\mathbf{U}}(t=0) = \ddot{\mathbf{U}}_0 = \mathbf{M}^{-1}\left(\mathbf{F}^0_\mathrm{ex} - \mathbf{K}\mathbf{U}_0 - \mathbf{C}\dot{\mathbf{U}}_0\right)\,.
	\label{eq:InitialAcc}
\end{equation}
Furthermore, during the initialization phase of the time integration method (i.e., for $t\,{=}\,0$ or equivalently $i\,{=}\,0$), it is necessary to know the displacement vector at $t\,{=}\,-\Delta t$ (or $i\,{=}\,-1$) as well \cite{BookChopra2011} -- cf. Eq.~\eqref{eq:CDM}. This value can be obtained by
\begin{equation}
	\mathbf{U}_{-1} = \mathbf{U}_0 - \Delta t\dot{\mathbf{U}}_0 - \cfrac{\Delta t^2}{2}\ddot{\mathbf{U}}_0\,.
\end{equation}
\subsection{Stiffness- and mass-stabilization}
To integrate the EVS-technique into time stepping schemes, it is only necessary to replace the original system matrices in Eq.~\eqref{eq:EquationMotion_i} with their stabilized/modified counterparts, resulting in
\begin{equation}
	\textcolor{red}{\mathbf{M}^\mathrm{Mod}}\ddot{\mathbf{U}}_i + \textcolor{red}{\mathbf{C}^\mathrm{Mod}}\dot{\mathbf{U}}_i + \textcolor{red}{\mathbf{K}^\mathrm{Mod}}\mathbf{U}_i = \mathbf{F}^i_\mathrm{ex}\,.
	\label{eq:EquationMotion_i_mod}
\end{equation}
Please note that the structure of Eq.~\eqref{eq:EquationMotion_i_mod} is identical to Eq.~\eqref{eq:EquationMotion_i}, suggesting that the numerical characteristics of the resulting time integration scheme are likely similar. However, it is important to acknowledge that substituting the expressions changes the numerical outcomes\footnote{Remark: When considering the $\alpha$-stabilization technique, it is worth noting that the system matrices can also be divided into a standard matrix and a stabilization matrix. The standard matrix is computed using $\alpha(\mathbf{x}_\mathrm{phys})\,{=}\,1$ in the physical domain and $\alpha(\mathbf{x}_\mathrm{fict})\,{=}\,0$ in the fictitious domain. On the other hand, the stabilization matrix is computed using $\alpha(\mathbf{x}_\mathrm{fict})\,{=}\,\alpha_0$ in the fictitious domain and $\alpha(\mathbf{x}_\mathrm{phys})\,{=}\,0$ in the physical domain. Hence, the same error analysis is applicable to the $\alpha$-method as well. Further exploration of this perspective will be presented in an upcoming publication that addresses a method aimed at mitigating the loss of accuracy introduced by stabilization.}. The magnitude of the error resulting from this adjustment depends on the choice of the stabilization parameter $\epsilon_{\mathrm{S}}$. It is crucial to strike a balance between accuracy and effective stabilization, as these goals are often conflicting with each other.

When utilizing the stabilized matrices, the resulting expression to calculate the displacement field at time step $i\,{+}\,1$ can be expressed as follows:
\begin{equation}
	\left(\cfrac{1}{\Delta t^2}\textcolor{red}{\mathbf{M}^\mathrm{Mod}} + \cfrac{1}{2\Delta t}\textcolor{red}{\mathbf{C}^\mathrm{Mod}}\right)\mathbf{U}_{i+1} =\; \mathbf{F}^i_\mathrm{ex} - \left[\left(\textcolor{red}{\mathbf{K}^\mathrm{Mod}} - \cfrac{2}{\Delta t^2}\textcolor{red}{\mathbf{M}^\mathrm{Mod}}\right)\mathbf{U}_i + \left(\cfrac{1}{\Delta t^2}\textcolor{red}{\mathbf{M}^\mathrm{Mod}} - \cfrac{1}{2\Delta t}\textcolor{red}{\mathbf{C}^\mathrm{Mod}}\right)\mathbf{U}_{i-1}\right] \,.
	\label{eq:CDM_stab}
\end{equation}
Equation~\eqref{eq:CDM_stab} can be again reformulated as a static equilibrium equation, resulting in
\begin{equation}
	\tilde{\mathbf{K}}^\mathrm{Mod}\mathbf{U}_{i+1} = \tilde{\mathbf{F}}^\mathrm{Mod}_i\,.
	\label{eq:StaticEquilibrium_mod_i}
\end{equation}

It is important to note that depending on the chosen variant of the EVS-technique (refer to Table~\ref{tab:Variants}), not all matrices need to be stabilized. Since the damping matrix is simply a linear combination of the mass and stiffness matrices, its stabilization depends on these two matrices and is not discussed separately.

For applications in explicit dynamics, improving the condition number through the EVS-technique is merely a side-effect and does not significantly impact the process. Since no system of equations is solved in each time step (assuming mass and damping matrices are diagonal), the conditioning of the stiffness matrix is not of utmost importance. Therefore, stabilizing $\mathbf{K}$ is not essentially required for explicit dynamics.

On the other hand, when considering implicit time integration schemes, an improved condition number of the system matrices becomes crucial for both direct and iterative solvers, which are used in each time step. As a result, a difference in the application of the EVS-technique for explicit and implicit dynamics may arise.
\subsection{Amplification matrix and load-operator}
To assess the implications of using the EVS-technique for explicit time stepping, it is necessary to re-evaluate the properties of the selected time integrator. The stability analysis of the CDM follows the methodology initially introduced in Ref.~\cite{ArticleBathe1972}. In this analysis, a single degree of freedom (SDOF) system is examined\footnote{Remark: For linear elastodynamics the superposition principle holds and therefore, it is possible to decouple the semi-discrete equations of motion by means of a modal decomposition, which renders the stiffness and mass matrices diagonal. This means that the movement of a structure is essentially governed by a weighted superposition of all mode shapes.}, and the amplification matrix $\mathbf{A}$, as well as the load-operator $\mathbf{L}$, are derived. The derivation of $\mathbf{A}$ and $\mathbf{L}$ is based on different versions of the equation of motion discussed in the previous subsections. For clarity, we will present the required expressions in a scalar format. The equation of motion is conventionally expressed as follows:
\begin{equation}
	m\ddot{u}_i + c\dot{u}_i + ku_i = f^i_\mathrm{ex}\,.
	\label{eq:EquationMotion_i_scalar}
\end{equation}
In contrast, the stabilized version of the equation of motion produces the following expression:
\begin{equation}
	\textcolor{red}{m^\mathrm{Mod}}\ddot{u}_i + \textcolor{red}{c^\mathrm{Mod}}\dot{u}_i + \textcolor{red}{k^\mathrm{Mod}}u_i = f^i_\mathrm{ex}\,.
	\label{eq:EquationMotion_i_stab_scalar}
\end{equation}
It is observed that changing from the unstabilized to the stabilized version of the equation of motion, one only needs to exchange the mass, damping, and stiffness parameters with their stabilized counterparts. Therefore, it is sufficient to derive $\mathbf{A}$ and $\mathbf{L}$ for the unstabilized scheme.
The central difference expressions are defined by
\begin{alignat}{2}
	\dot{u}_i & = && \cfrac{1}{2\Delta t} (u_{i+1} - u_{i-1})\,, \label{eq:CD_v_i_scalar}\\
	\shortintertext{and} & && \nonumber \\[-18pt]
	\ddot{u}_i & = && \cfrac{1}{\Delta t^2} (u_{i+1} - 2u_{i} + u_{i-1})\,. \label{eq:CD_a_i_scalar}
\end{alignat}
Equations~\eqref{eq:CD_v_i_scalar} and \eqref{eq:CD_a_i_scalar} are substituted into Eq.~\eqref{eq:EquationMotion_i_scalar}, which yields
\begin{equation}
	\cfrac{1}{\Delta t^2} \,m\, (u_{i+1} - 2u_{i} + u_{i-1}) + \cfrac{1}{2\Delta t} \,c\, (u_{i+1} - u_{i-1}) + ku_{i} = f^i_\mathrm{ex}\,.
\end{equation}
By separating all terms related to $u_{i-1}$, $u_i$, and $u_{i+1}$, we can express the time integration scheme in a recursive relation of the form:
\begin{equation}
	\mathbf{u}_{i+1} = \mathbf{A} \mathbf{u}_i + \mathbf{L} \label{eq:RecursiveTimeInt}
\end{equation}
where $\mathbf{u}_{i+1}$ and $\mathbf{u}_i$ represent the vectors storing the solution quantities, which can be displacements, velocities, or accelerations depending on the selected time integrator. In the case of the standard CDM implementation, Eq.~\eqref{eq:RecursiveTimeInt} takes the following form:
\begin{equation}
	\begin{Bmatrix}
		u_{i+1}\\
		u_{i}
	\end{Bmatrix} = 
	\mathbf{A}
	\begin{Bmatrix}
		u_{i}\\
		u_{i-1}
	\end{Bmatrix} + 
	\mathbf{L}\,. \label{eq:Recurr2Terms}
\end{equation}
The structure of the amplification matrix $\mathbf{A}$ can be described as follows:
\begin{equation}
	\mathbf{A} = \hat{k}^{-1}
	\begin{bmatrix}
		A_{11} & A_{12}\\
		\hat{k} & 0
	\end{bmatrix}
\end{equation}
with
\begin{alignat}{1}
	A_{11} & =  2m - \Delta t^2 k\,, \\
	\shortintertext{and} & \nonumber \\[-18pt]
	A_{12} & = \cfrac{\Delta t}{2} \,c - m = \hat{k} - 2m\,.
\end{alignat}
The load-operator, on the other hand, takes the following form:
\begin{equation}
	\mathbf{L} =\hat{k}^{-1}
	\begin{bmatrix}
		L_1 \\
		0
	\end{bmatrix},
	\label{eq:L_CDM}
\end{equation}
with
\begin{equation}
	L_1 = \Delta t^2 f^i_\mathrm{ex}\,.
\end{equation}
For the sake of a compact notation, the auxiliary variable $\hat{k}$ is defined as
\begin{equation}
	\hat{k} = m + \cfrac{\Delta t}{2} \,c\,.
\end{equation}

To assess the stability of a direct time integration scheme, the spectral radius of the amplification matrix $\mathbf{A}$ is examined by solving the eigenvalue problem
\begin{equation}
	\mathbf{A}\mathbf{x} = \lambda\mathbf{x}\,,
\end{equation}
which yields non-trivial solutions only for
\begin{equation}
	\det(\mathbf{A} - \lambda\mathbf{I}) = 0\,. \label{eq:DetA}
\end{equation}
In order for a time integration scheme to be stable, it is necessary for the absolute value of the spectral radius to satisfy the condition $|\varrho(\mathbf{A})|\,{\le}\,1$.

The characteristic equation (polynomial) $p(\lambda)$ of a general $2\times2$ matrix $\mathbf{A}$ can be expressed as
\begin{equation}
p_2(\lambda) = \lambda^2 - A_1\lambda + A_0 \qquad \text{with} \quad A_1 = \mathrm{tr}(\mathbf{A}) \quad \text{and} \quad A_0 = \det(\mathbf{A})\,.
\label{eq:CharacteristicEq_p=2}
\end{equation}
Given the specific structure of the amplification matrix in the CDM, the two constants correspond to $A_1\,{=}\,\hat{k}^{-1}A_{11}$ and $A_0\,{=}\,-\hat{k}^{-1}A_{12}$, respectively. Consequently, the characteristic polynomial can be expressed as:
\begin{equation}
p_2(\lambda) = \lambda^2 - \hat{k}^{-1}A_{11}\lambda - \hat{k}^{-1}A_{12} = \hat{k}^{-1} \left( \hat{k}\lambda^2 - A_{11}\lambda - A_{12} \right)\,.
\label{eq:CharacteristicEq_p=2_A}
\end{equation}
To establish the range of $\Delta t$ in which the CDM remains stable, the following approach can be employed: Initially, we substitute $\lambda\,{=}\,\pm 1$ into Eq.~\eqref{eq:CharacteristicEq_p=2_A} and subsequently solve the resulting quadratic equation to obtain the corresponding values of $\Delta t$. Evaluating the equation for $\lambda\,{=}\,+1$ leads us to the conclusion that either $k\,{=}\,0$ or $\Delta t\,{=}\,0$, which is not physically reasonable for investigating structural dynamics problems. However, from the condition $\lambda\,{=}\,- 1$, the following expression can be derived:
\begin{equation}
	0 = 4m - \Delta t^2 k
\end{equation}
assuming that none of the quantities $m$, $c$, $k$, or $\Delta t$ can be zero or negative. From this expression we infer that the critical time step size for the CDM is:
\begin{equation}
	\Delta t_\mathrm{cr} = \pm\sqrt{\cfrac{4}{m^{-1}k}} = \pm\cfrac{2}{\omega}\,,
	\label{eq:CDM_dt_crit_1}
\end{equation}
where the result with the negative sign is discarded due to physical reasons. Interestingly, the derived value for the critical time step size is independent of the introduced physical damping\footnote{Remark: In \emph{truly explicit} time integration schemes, the effective stiffness matrix is determined solely by the mass $m$ and remains unaffected by the damping parameter $c$. Consequently, the critical time step size exhibits a dependence on the value of $c$. However, for the CDM, the effective stiffness matrix is dependent on both $m$ and $c$, leading to the critical time step size $\Delta t_\mathrm{cr}$ being exclusively determined by $k$ and $m$, while being independent of $c$.}. Consequently, the critical time step size should be selected to satisfy the condition $\Delta t_\mathrm{cr} \in\;\; ]0, \sfrac{2}{\omega}]$, where $\omega$ represents the natural frequency of the (undamped) system.

When examining the stability limit for the time stepping scheme based on Eq.~\eqref{eq:EquationMotion_i_stab_scalar}, the analysis is analogous to the one previously discussed. Hence, we can proceed by substituting the original matrices with their stabilized counterparts, leading to the following expression:
\begin{equation}
	\Delta t^\star_\mathrm{cr} = \sqrt{\cfrac{4}{\left( m^\star \right)^{-1}k^\star}} = \cfrac{2}{\omega^\star}\,.
	\label{eq:CDM_dt_crit_2}
\end{equation}
At this point, we want to stress that this expression -- Eq.~\eqref{eq:CDM_dt_crit_2} -- is utilized to obtain the numerical results for the critical time step size listed in Sect.~\ref{sec:EVST_dtcr}. Depending on the chosen variant of the EVS-technique, all quantities denoted by a star $\square^\star$ represent either the original quantity $\square$ or its modified version $\square^\mathrm{Mod}$.
\section{Critical time step size}
\label{sec:EVST_dtcr}
\begin{figure}[t!]
	\centering
	\subfloat[Cut element]{\includegraphics[scale=1.0]{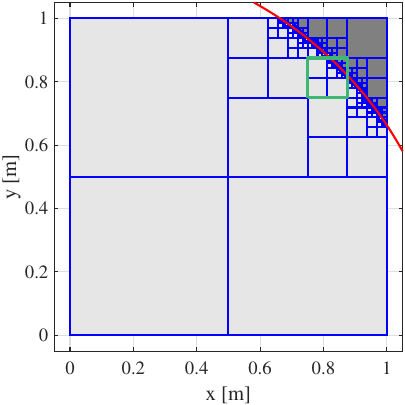}}
	\hspace*{1.5cm}
	\subfloat[Detail view (seagreen boundary)]{\includegraphics[scale=1.0]{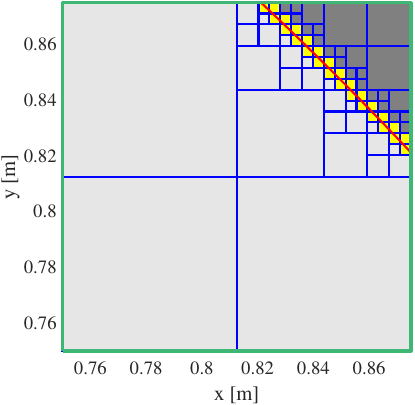}}
	\caption{Model of a rectangular plate intersected by a circular domain boundary (Discretization: $1\times1$ finite elements, Subdivision depth: $k\,{=}\,8$; Red solid line: physical boundary of the circular hole). \emph{Color coding for the numerical integration} -- Dark gray: integration domains located in the physical domain; Light gray: integration domains located in the fictitious domain; Yellow: cut integration domains (leaf cells of the tree data structure). \label{fig:CircularVoid}}
\end{figure}

In order to evaluate the impact of the EVS-technique on increasing the critical time step size $\Delta t_\mathrm{cr}$, while disregarding its effect on the accuracy of the numerical simulation, a simplified model comprising a single cut element is analyzed. To this end, we consider the following two-dimensional setup: Figure~\ref{fig:CircularVoid} illustrates a single finite element intersected by a circular void region. The circle's origin coincides with the lower left corner vertex of the quadrilateral element. In the following investigations, the radius of the circle $r_{\circ}$ is set to $1.2$, while the element's side lengths remain fixed at $1$, resulting in a volume fraction $\chi$ (with respect to the physical domain) of $4.9\%$. As described in Section~\ref{sec:FCM}, a composite numerical integration scheme based on a spacetree decomposition of the integration domain is employed, with the subdivision level set to $k\,{=}\,8$. For the simulations, a plane stress state is assumed, and the material properties chosen are those of steel: Young's modulus $E\,{=}\,210\,$GPa, Poisson's ratio $\nu\,{=}\,0.3$, and mass density $\rho\,{=}\,7850\,\sfrac{\mathrm{kg}}{\mathrm{m}^3}$.

Using this example, the different variants of the EVS-technique (see Table~\ref{tab:Variants}) are evaluated. For the initial investigation, the two parameters of the stabilization technique, $\epsilon_\lambda$ and $\epsilon_\mathrm{S}$, are chosen as $10^{-4}$ and $10^{-4}$, respectively. Moreover, the numerical results obtained using the $\alpha$-method are also included as a reference for different values of $\alpha_0 \in \{0,\,10^{-12},\,10^{-5}\}$, denoted as \emph{variant 0}.

The numerical results for the critical time step size are listed in Table~\ref{tab:dtcr_1}, where all values are normalized with respect to the critical time increment obtained for \textcolor{Matlab6}{\emph{variant 0e}}. For different polynomial orders $p\in[1,8]$, the reference values are $[27.1141,\, 16.4569,\, 12.7724,\, 5.60615,\, 4.33158,\, 3.13476,\, 2.44573,\, 1.82912]\,\mu$s. Note that these results correspond to the standard $\alpha$-stabilization with $\alpha_0\,{=}\,10^{-5}$ and a lumped mass matrix. The lumping procedure is discussed in detail in Ref.~\cite{ArticleJoulaian2014a}, where it is shown that only HRZ-lumping ensures the positive-definiteness of the mass matrix for cut elements\footnote{Remark: In the context of the spectral cell method (SCM), it is important to note that mass lumping of uncut elements (standard spectral elements) is accomplished through a nodal quadrature technique, while cut elements are diagonalized using the HRZ-method.}. In order to be useful, the proposed EVS-technique should yield an improved performance regarding higher critical time step sizes, while maintaining accuracy.
\begin{table}[t!]
	\caption{Normalized critical time step size $\Delta t^\mathrm{norm}_\mathrm{cr}$ for different variants of the EVS-technique and polynomial degrees $p$ of the shape functions. Parameters: $\epsilon_\lambda\,{=}\,10^{-4}$, ${}_3\epsilon^\lambda_{\mathrm{S}}$ with $\epsilon_\mathrm{S}\,{=}\,10^{-4}$ (\emph{variant 0}: $\alpha$-method; \emph{variants 1 -- 3}: EVS-technique $\rightarrow$ $\alpha_0\,{=}\,0$).\label{tab:dtcr_1}}
	\centering
	\begin{tabular}{cccccccccccc}
		\toprule[2pt]
		\phantom{x} & \diagbox{No.}{$p$} & \phantom{xxx} & $1$ & $2$ & $3$ & $4$ & $5$ & $6$ & $7$ & $8$ & \phantom{x}\\
		\midrule[1pt]
		& 0a && 0.9135 & 0.9821 & 0.9738 & 0.8813 & 0.9082 & 0.8262 & 0.8286 & 0.7900 \\ 
		& 0b && 0.4792 & 0.3322 & 0.2328 & 0.3283 & 0.3007 & 0.3523 & 0.3424 & 0.4131 \\ 
		& 0c && 0.9135 & 0.9821 & 0.9738 & 0.8813 & 0.9082 & 0.8262 & 0.8286 & 0.7900 \\ 
		& 0d && 0.4792 & 0.3324 & 0.2742 & 0.4646 & 0.4647 & 0.5283 & 0.5749 & 0.6469 \\ 
		\rowcolor{Grey2}
		& 0e && 1      & 1      & 1      & 1      & 1      & 1      & 1      & 1      \\ 
		& 0f && 0.5978 & 0.6070 & 0.5915 & 0.9845 & 1.0379 & 1.1036 & 1.1800 & 1.3030 \\\hline\\[-1.75ex]
		& 1a && 0.9135 & 0.9530 & 0.5804 & 0.3278 & 0.3418 & 0.2660 & 0.2500 & 0.2103 \\ 
		& 1b && 0.4792 & 0.0056 & $\dagger$ & $\dagger$ & $\dagger$ & $\dagger$ & $\dagger$ & $\dagger$ \\\hline\\[-1.75ex]
		& 2a && 1.3679 & 1.1338 & 1.1331 & 1.4304 & 1.3587 & 1.4469 & 1.2487 & 1.2914 \\ 
		\rowcolor{Grey1}
		& 2b && 1.3125 & 1.1639 & 1.1758 & 1.6458 & 1.4768 & 1.6851 & 1.3060 & 1.3487 \\ 
		\rowcolor{Grey1}
		& 2c && 1.3420 & 1.1516 & 1.1686 & 1.6552 & 1.4635 & 1.6893 & 1.2900 & 1.3336 \\ 
		& 2d && 0.7439 & 0.7950 & 0.5756 & 1.0265 & 1.1329 & 1.1655 & 1.2632 & 1.4538 \\ 
		& 2e && 0.7409 & 0.7972 & 0.6709 & 1.3317 & 1.3062 & 1.4989 & 1.3422 & 1.5101 \\ 
		& 2f && 0.7434 & 0.8066 & 0.7242 & 1.3665 & 1.3580 & 1.5479 & 1.3589 & 1.5316 \\\hdashline\\[-1.75ex]
		& 2g && 0.9135 & 0.9821 & 1.0198 & 1.3227 & 1.2629 & 1.2577 & 1.1643 & 1.1961 \\ 
		& 2h && 0.9135 & 0.9821 & 1.0198 & 1.3227 & 1.2629 & 1.2577 & 1.1643 & 1.1961 \\ 
		& 2i && 0.9135 & 0.9821 & 1.0198 & 1.3227 & 1.2629 & 1.2577 & 1.1643 & 1.1961 \\ 
		& 2j && 0.4792 & 0.3322 & 0.2752 & 0.5669 & 0.5359 & 0.4735 & 0.4675 & 0.5239 \\ 
		& 2k && 0.4792 & 0.3322 & 0.2752 & 0.5669 & 0.5359 & 0.4735 & 0.4675 & 0.5239 \\ 
		& 2l && 0.4792 & 0.3322 & 0.2752 & 0.5669 & 0.5359 & 0.4735 & 0.4675 & 0.5239 \\\hline\\[-1.75ex]
		& 3a && 1.3679 & 1.1291 & 1.1304 & 1.4284 & 1.3577 & 1.4463 & 1.2462 & 1.2902 \\ 
		& 3b && 1.3125 & 1.1597 & 1.1739 & 1.6443 & 1.4762 & 1.6848 & 1.3038 & 1.3477 \\ 
		& 3c && 1.3420 & 1.1469 & 1.1662 & 1.6532 & 1.4626 & 1.6889 & 1.2876 & 1.3323 \\ 
		& 3d && 0.7439 & 0.7949 & 0.5748 & 1.0260 & 1.1325 & 1.1651 & 1.2631 & 1.4535 \\ 
		& 3e && 0.7409 & 0.7972 & 0.6705 & 1.3313 & 1.3059 & 1.4988 & 1.3419 & 1.5097 \\ 
		& 3f && 0.7434 & 0.8066 & 0.7240 & 1.3663 & 1.3578 & 1.5479 & 1.3587 & 1.5312 \\\hdashline\\[-1.75ex] 
		& 3g && 0.9135 & 0.9530 & 1.0008 & 1.3194 & 1.2612 & 1.2542 & 1.1613 & 1.1945 \\ 
		& 3h && 0.9135 & 0.9530 & 1.0008 & 1.3194 & 1.2612 & 1.2542 & 1.1613 & 1.1945 \\ 
		& 3i && 0.9135 & 0.9530 & 1.0008 & 1.3194 & 1.2612 & 1.2542 & 1.1613 & 1.1945 \\ 
		& 3j && 0.4792 & 0.0056 & 0.0015 & 0.0066 & 0.0010 & 0.0001 & $\dagger$ & $\dagger$ \\ 
		& 3k && 0.4792 & 0.0056 & 0.0015 & 0.0066 & 0.0010 & 0.0001 & $\dagger$ & $\dagger$ \\ 
		& 3l && 0.4792 & 0.0056 & 0.0015 & 0.0066 & 0.0010 & 0.0001 & $\dagger$ & $\dagger$ \\ 
		\bottomrule[2pt]
	\end{tabular}
	\caption*{\raggedright{The symbol $\dagger$ denotes values that are either below a numerical threshold of $10^{-4}$ or complex.}}
\end{table}%

When choosing a suitable variant of the proposed stabilization scheme, it is essential to consider different application scenarios. In implicit time integration methods, such as the trapezoidal rule of the Newmark family, the ill-conditioning (with respect to matrix inversion) of the stiffness and mass matrices becomes a major concern. This is because implicit schemes require solving a linear system of equations in each time step, where the effective stiffness matrix consists of a linear combination of both stiffness and mass matrices. In such cases, ill-conditioning can result in inaccurate solutions when using direct equation solvers, or significantly increased iteration counts when employing iterative equation solvers. Hence, it is advisable to stabilize both the stiffness and mass matrices in the context of implicit methods. Consequently, appropriate options should be taken from group \emph{3}.

On the other hand, in explicit time integration schemes like the CDM, the ill-conditioning (with respect to matrix inversion) of the system matrices is of lesser importance, i.e., lowering the condition number for inversion is not a priority. This is because only simple matrix--vector products are required to compute the internal force vector in linear problems, and no system of equations needs to be solved as long as the mass matrix for truly explicit schemes and the mass and damping matrices for explicit schemes are diagonal. Thus, a different notion of ill-conditioning (with respect to matrix-vector products) becomes important and must be considered for explicit methods. Considering that ill-conditioning in matrix--vector products is less likely to pose a problem in numerical methods, it appears reasonable to focus on stabilizing only the mass matrix for explicit dynamics. Consequently, appropriate options should be taken from group \emph{2}.

In the subsequent paragraphs, we will analyze the numerical results presented in Table~\ref{tab:dtcr_1}. One important observation is that the numerical time increments obtained for \emph{variants 2g}, \emph{2h}, and \emph{2i} are identical. In these variants of the EVS-scheme, only the mass matrix is stabilized, and the eigenvalue decomposition is based on the lumped mass matrix. Therefore, it is essential to recall that the eigenvectors of a diagonal matrix contain only one non-zero value each; otherwise, off-diagonal elements would arise from the spectral decomposition. Consequently, the computed mass stabilization matrix is diagonal as well. In essence, the stabilization entails adding additional mass to those nodes within the fictitious domain, whose shape functions exhibit little support in the physical domain. Since the mass stabilization matrix is already diagonal, the application of any mass lumping technique becomes inconsequential. This explains why the resulting time step sizes for all three variants are identical. In our specific example (refer to Fig.~\ref{fig:CircularVoid}), it is evident that all nodes located near the origin of the global coordinate system require stabilization. Keep in mind that the size of the stabilization zone varies depending on the chosen value of $\epsilon_\lambda$. Consequently, in order to offer reliable recommendations regarding the most suitable variant of the EVS-technique to employ, it becomes imperative to thoroughly study the effects of both $\epsilon_\mathrm{S}$ and $\epsilon_\lambda$.

Upon analyzing this initial example, it becomes apparent that the performances of \emph{variants 2a} to \emph{2c} and \emph{3a} to \emph{3c} are markedly better compared to all other options. To narrow down the number of possible approaches, all variants from group \emph{3} are excluded from further considerations as we want to keep the complexity of the proposed stabilization approach as low as possible. Hence, it is preferred to only stabilize the mass matrix, while the stiffness matrix remains unstabilized (use options from group \emph{2}). Furthermore, it is observed that by stabilizing the stiffness matrix no additional advantages are gained. Despite performing reasonably well, \emph{variant 2a} is also excluded from the set of suitable options. This is due to the fact that the mass stabilization matrix is not lumped and therefore, severe performance penalties for explicit time integration schemes are incurred.

To reach a definitive conclusion regarding the most suitable variant of the EVS-technique for explicit time integration schemes, we conduct a thorough investigation by varying the stabilization parameter $\epsilon_\mathrm{S}$ and the threshold value $\epsilon_\lambda$. The results of these parameter studies for \emph{variants 2b} and \emph{2c} are provided in Tables~\ref{tab:dtcr_2b} and \ref{tab:dtcr_2c}, respectively.

\begin{table}[b!]
	\caption{Normalized critical time step size $\Delta t^\mathrm{norm}_\mathrm{cr}$ for \emph{variant 2b} of the EVS-technique and different polynomial degrees $p$ of the shape functions. \label{tab:dtcr_2b}}
	\centering
	\begin{tabular}{cccccccclll}
		\toprule[2pt]
		$p\,{=}\,1$ & $2$ & $3$ & $4$ & $5$ & $6$ & $7$ & $8$ & \multicolumn{3}{c}{Parameters} \\
		\midrule[1pt]
		3.3801	& 2.4723 	& 1.8662 	& 2.8171 	& 2.5957 	& 2.6823 	& 2.6667 	& 2.8452 &  & SE-GLL  & CMM\\
		5.9659	& 3.9090 	& 2.8107 	& 3.9757 	& 3.4586 	& 3.4278 	& 3.3017 	& 3.4368 &  & SE-GLL  & LMM\\\hline\\[-1.75ex]
		1      	& 1     	& 1     	& 1     	& 1     	& 1         & 1     	& 1      & $\epsilon_\lambda\,{=}\,0$   	& $\epsilon_\mathrm{S}\,{=}\,0$			& $\alpha_0\,{=}\,10^{-5}$\\\hline\\[-1.75ex]
		0.9141 	& 0.9829 	& 0.9773 	& 0.9057 	& 0.9366 	& 0.8846 	& 0.8446 	& 0.8142 & $\epsilon_\lambda\,{=}\,10^{-4}$ & $\epsilon_\mathrm{S}\,{=}\,10^{-7}$, 	& $\alpha_0\,{=}\,0$\\
		0.9194 	& 0.9893 	& 0.9984 	& 0.9855 	& 1.0094 	& 1.0105 	& 0.8948 	& 0.8793 & $\epsilon_\lambda\,{=}\,10^{-4}$ & $\epsilon_\mathrm{S}\,{=}\,10^{-6}$, 	& $\alpha_0\,{=}\,0$\\
		0.9699 	& 1.0337 	& 1.0601 	& 1.2013 	& 1.1722 	& 1.2661 	& 1.0281 	& 1.0355 & $\epsilon_\lambda\,{=}\,10^{-4}$ & $\epsilon_\mathrm{S}\,{=}\,10^{-5}$, 	& $\alpha_0\,{=}\,0$\\
		1.3125  & 1.1639    & 1.1758    & 1.6458    & 1.4768    & 1.6851    & 1.3060    & 1.3487 & $\epsilon_\lambda\,{=}\,10^{-4}$ & $\epsilon_\mathrm{S}\,{=}\,10^{-4}$, 	& $\alpha_0\,{=}\,0$\\
		1.8367 	& 1.4861 	& 1.3806 	& 2.3930 	& 1.9463 	& 2.3321 	& 1.7337 	& 1.8771 & $\epsilon_\lambda\,{=}\,10^{-4}$ & $\epsilon_\mathrm{S}\,{=}\,10^{-3}$, 	& $\alpha_0\,{=}\,0$\\
		1.9099 	& 1.8156 	& 1.6060 	& 3.2533 	& 2.5503 	& 3.0144 	& 2.3236 	& 2.5995 & $\epsilon_\lambda\,{=}\,10^{-4}$ & $\epsilon_\mathrm{S}\,{=}\,10^{-2}$, 	& $\alpha_0\,{=}\,0$\\
		1.9357 	& 2.4191 	& 1.9734 	& 3.8364 	& 2.9883 	& 3.3577 	& 2.9900 	& 3.2254 & $\epsilon_\lambda\,{=}\,10^{-4}$ & $\epsilon_\mathrm{S}\,{=}\,10^{-1}$, 	& $\alpha_0\,{=}\,0$\\
		2.1208 	& 3.5499 	& 2.2196 	& 4.0026 	& 3.4376 	& 3.5552 	& 3.3999 	& 3.4234 & $\epsilon_\lambda\,{=}\,10^{-4}$ & $\epsilon_\mathrm{S}\,{=}\,1$, 		& $\alpha_0\,{=}\,0$\\\hdashline\\[-1.75ex]
		1.0005 	& 1.0008 	& 1.0023 	& 1.0114 	& 1.0124 	& 1.0214 	& 1.0038 	& 1.0046 & $\epsilon_\lambda\,{=}\,10^{-4}$ & $\epsilon_\mathrm{S}\,{=}\,10^{-7}$, 	& $\alpha_0\,{=}\,10^{-5}$\\
		1.0050 	& 1.0075 	& 1.0175 	& 1.0650 	& 1.0591 	& 1.0975 	& 1.0242 	& 1.0279 & $\epsilon_\lambda\,{=}\,10^{-4}$ & $\epsilon_\mathrm{S}\,{=}\,10^{-6}$, 	& $\alpha_0\,{=}\,10^{-5}$\\
		1.0479 	& 1.0529 	& 1.0698 	& 1.2488 	& 1.1948 	& 1.2993 	& 1.1086 	& 1.1200 & $\epsilon_\lambda\,{=}\,10^{-4}$ & $\epsilon_\mathrm{S}\,{=}\,10^{-5}$, 	& $\alpha_0\,{=}\,10^{-5}$\\
		1.3496 	& 1.1887 	& 1.1792 	& 1.6622 	& 1.4834 	& 1.6918 	& 1.3341 	& 1.3735 & $\epsilon_\lambda\,{=}\,10^{-4}$ & $\epsilon_\mathrm{S}\,{=}\,10^{-4}$, 	& $\alpha_0\,{=}\,10^{-5}$\\
		1.8380 	& 1.5184 	& 1.3778 	& 2.3746 	& 1.9463 	& 2.3297 	& 1.7365 	& 1.8756 & $\epsilon_\lambda\,{=}\,10^{-4}$ & $\epsilon_\mathrm{S}\,{=}\,10^{-3}$, 	& $\alpha_0\,{=}\,10^{-5}$\\
		1.9114 	& 1.8370 	& 1.5754 	& 3.2505 	& 2.5446 	& 3.0072 	& 2.3254 	& 2.5882 & $\epsilon_\lambda\,{=}\,10^{-4}$ & $\epsilon_\mathrm{S}\,{=}\,10^{-2}$, 	& $\alpha_0\,{=}\,10^{-5}$\\
		1.9393 	& 2.4999 	& 1.9320 	& 3.8348 	& 2.9780 	& 3.3541 	& 2.9991 	& 3.1999 & $\epsilon_\lambda\,{=}\,10^{-4}$ & $\epsilon_\mathrm{S}\,{=}\,10^{-1}$, 	& $\alpha_0\,{=}\,10^{-5}$\\
		2.1434 	& 3.4157 	& 2.1609 	& 4.0029 	& 3.4168 	& 3.5552 	& 3.4322 	& 3.4121 & $\epsilon_\lambda\,{=}\,10^{-4}$ & $\epsilon_\mathrm{S}\,{=}\,1$, 		& $\alpha_0\,{=}\,10^{-5}$\\\hline\\[-1.75ex]
		0.9135 	& 1.4388 	& 1.4968 	& 2.7622 	& 2.2891 	& 2.7606 	& 2.1583 	& 2.3967 & $\epsilon_\lambda\,{=}\,10^{-7}$ & $\epsilon_\mathrm{S}\,{=}\,10^{-2}$, 	& $\alpha_0\,{=}\,0$\\
		0.9135 	& 1.4388 	& 1.5767 	& 2.7956 	& 2.4464 	& 2.7758 	& 2.2874 	& 2.4359 & $\epsilon_\lambda\,{=}\,10^{-6}$ & $\epsilon_\mathrm{S}\,{=}\,10^{-2}$, 	& $\alpha_0\,{=}\,0$\\
		0.9135 	& 1.7219 	& 1.6060 	& 2.9812 	& 2.5423 	& 2.9404 	& 2.3000 	& 2.5236 & $\epsilon_\lambda\,{=}\,10^{-5}$ & $\epsilon_\mathrm{S}\,{=}\,10^{-2}$, 	& $\alpha_0\,{=}\,0$\\
		1.9099 	& 1.8156 	& 1.6060 	& 3.2533 	& 2.5503 	& 3.0144 	& 2.3236 	& 2.5995 & $\epsilon_\lambda\,{=}\,10^{-4}$ & $\epsilon_\mathrm{S}\,{=}\,10^{-2}$, 	& $\alpha_0\,{=}\,0$\\
		1.9099 	& 1.8156 	& 1.7216 	& 3.2673 	& 2.6999 	& 3.1005 	& 2.4031 	& 2.6316 & $\epsilon_\lambda\,{=}\,10^{-3}$ & $\epsilon_\mathrm{S}\,{=}\,10^{-2}$, 	& $\alpha_0\,{=}\,0$\\
		2.4771 	& 2.2155 	& 1.8489 	& 3.4627 	& 2.7596 	& 3.2902 	& 2.4375 	& 2.7129 & $\epsilon_\lambda\,{=}\,10^{-2}$ & $\epsilon_\mathrm{S}\,{=}\,10^{-2}$, 	& $\alpha_0\,{=}\,0$\\\hdashline\\[-1.75ex]
		1      	& 1      	& 1      	& 1      	& 1      	& 1      	& 1      	& 1      & $\epsilon_\lambda\,{=}\,10^{-7}$ & $\epsilon_\mathrm{S}\,{=}\,10^{-2}$, 	& $\alpha_0\,{=}\,10^{-5}$\\
		1      	& 1      	& 1      	& 1.3172 	& 1.2144 	& 1.1888 	& 1.1082 	& 1.1538 & $\epsilon_\lambda\,{=}\,10^{-6}$ & $\epsilon_\mathrm{S}\,{=}\,10^{-2}$, 	& $\alpha_0\,{=}\,10^{-5}$\\
		1      	& 1.1450 	& 1.2414 	& 1.6544 	& 1.5682 	& 1.7337 	& 1.4607 	& 1.5582 & $\epsilon_\lambda\,{=}\,10^{-5}$ & $\epsilon_\mathrm{S}\,{=}\,10^{-2}$, 	& $\alpha_0\,{=}\,10^{-5}$\\
		1.9114 	& 1.8370 	& 1.5754 	& 3.2505 	& 2.5446 	& 3.0072 	& 2.3254 	& 2.5882 & $\epsilon_\lambda\,{=}\,10^{-4}$ & $\epsilon_\mathrm{S}\,{=}\,10^{-2}$, 	& $\alpha_0\,{=}\,10^{-5}$\\
		1.9114 	& 1.8160 	& 1.7203 	& 3.2670 	& 2.6997 	& 3.1004 	& 2.4036 	& 2.6321 & $\epsilon_\lambda\,{=}\,10^{-3}$ & $\epsilon_\mathrm{S}\,{=}\,10^{-2}$, 	& $\alpha_0\,{=}\,10^{-5}$\\
		2.4780 	& 2.2167 	& 1.8492 	& 3.4628 	& 2.7597 	& 3.2903 	& 2.4381 	& 2.7134 & $\epsilon_\lambda\,{=}\,10^{-2}$ & $\epsilon_\mathrm{S}\,{=}\,10^{-2}$, 	& $\alpha_0\,{=}\,10^{-5}$\\\hline\\[-1.75ex]
		9.5581 	& 4.2002 	& 3.0151 	& 4.4653 	& 4.3341 	& 4.8341 	& 4.1968 	& 4.1005 & $\epsilon_\lambda\,{=}\,10^{-2}$ & $\epsilon_\mathrm{S}\,{=}\,1$, 		& $\alpha_0\,{=}\,0$\\
		9.5582 	& 4.2012 	& 3.0149 	& 4.4665 	& 4.3340 	& 4.8342 	& 4.1968 	& 4.1006 & $\epsilon_\lambda\,{=}\,10^{-2}$ & $\epsilon_\mathrm{S}\,{=}\,1$, 		& $\alpha_0\,{=}\,10^{-5}$\\
		\bottomrule[2pt]
	\end{tabular}
\end{table}%
\begin{table}[t!]
	\caption{Normalized critical time step size $\Delta t^\mathrm{norm}_\mathrm{cr}$ for \emph{variant 2c} of the EVS-technique and different polynomial degrees $p$ of the shape functions. \label{tab:dtcr_2c}}
	\centering
	\begin{tabular}{cccccccclll}
		\toprule[2pt]
		$p\,{=}\,1$ & $2$ & $3$ & $4$ & $5$ & $6$ & $7$ & $8$ & \multicolumn{3}{c}{Parameters} \\
		\midrule[1pt]
		3.3801	& 2.4723 	& 1.8662 	& 2.8171 	& 2.5957 	& 2.6823 	& 2.6667 	& 2.8452 &  & SE-GLL  & CMM\\
		5.9659	& 3.9090 	& 2.8107 	& 3.9757 	& 3.4586 	& 3.4278 	& 3.3017 	& 3.4368 &  & SE-GLL  & LMM\\\hline\\[-1.75ex]
		1      	& 1     	& 1     	& 1     	& 1     	& 1         & 1     	& 1      & $\epsilon_\lambda\,{=}\,0$   	& $\epsilon_\mathrm{S}\,{=}\,0$			& $\alpha_0\,{=}\,10^{-5}$\\\hline\\[-1.75ex]
		0.9142 	& 0.9827 	& 0.9764 	& 0.9025 	& 0.9334 	& 0.8798 	& 0.8431 	& 0.8124 & $\epsilon_\lambda\,{=}\,10^{-4}$ & $\epsilon_\mathrm{S}\,{=}\,10^{-7}$, 	& $\alpha_0\,{=}\,0$\\
		0.9201 	& 0.9877 	& 0.9936 	& 0.9759 	& 1.0010 	& 0.9998 	& 0.8902 	& 0.8744 & $\epsilon_\lambda\,{=}\,10^{-4}$ & $\epsilon_\mathrm{S}\,{=}\,10^{-6}$, 	& $\alpha_0\,{=}\,0$\\
		0.9759 	& 1.0255 	& 1.0492 	& 1.1827 	& 1.1575 	& 1.2518 	& 1.0178 	& 1.0248 & $\epsilon_\lambda\,{=}\,10^{-4}$ & $\epsilon_\mathrm{S}\,{=}\,10^{-5}$, 	& $\alpha_0\,{=}\,0$\\
		1.3420 	& 1.1516 	& 1.1686 	& 1.6552 	& 1.4635 	& 1.6893 	& 1.2900 	& 1.3336 & $\epsilon_\lambda\,{=}\,10^{-4}$ & $\epsilon_\mathrm{S}\,{=}\,10^{-4}$, 	& $\alpha_0\,{=}\,0$\\
		1.8409 	& 1.4834 	& 1.4010 	& 2.4966 	& 2.0092 	& 2.4241 	& 1.7361 	& 1.8825 & $\epsilon_\lambda\,{=}\,10^{-4}$ & $\epsilon_\mathrm{S}\,{=}\,10^{-3}$, 	& $\alpha_0\,{=}\,0$\\
		1.9099 	& 1.8156 	& 1.6060 	& 3.2533 	& 2.5503 	& 3.0144 	& 2.3236 	& 2.5995 & $\epsilon_\lambda\,{=}\,10^{-4}$ & $\epsilon_\mathrm{S}\,{=}\,10^{-2}$, 	& $\alpha_0\,{=}\,0$\\
		1.3349 	& 3.1825 	& 2.1520 	& 3.6006 	& 3.3236 	& 3.6439 	& 3.0739 	& 3.1754 & $\epsilon_\lambda\,{=}\,10^{-4}$ & $\epsilon_\mathrm{S}\,{=}\,10^{-1}$, 	& $\alpha_0\,{=}\,0$\\
		2.1208 	& 3.5499 	& 2.2196 	& 4.0026 	& 3.4376 	& 3.5552 	& 3.3999 	& 3.4234 & $\epsilon_\lambda\,{=}\,10^{-4}$ & $\epsilon_\mathrm{S}\,{=}\,1$, 		& $\alpha_0\,{=}\,0$\\\hdashline\\[-1.75ex]
		1.0006 	& 1.0005 	& 1.0018 	& 1.0097 	& 1.0107 	& 1.0192 	& 1.0034 	& 1.0041 & $\epsilon_\lambda\,{=}\,10^{-4}$ & $\epsilon_\mathrm{S}\,{=}\,10^{-7}$, 	& $\alpha_0\,{=}\,10^{-5}$\\
		1.0056 	& 1.0048 	& 1.0141 	& 1.0580 	& 1.0532 	& 1.0902 	& 1.0221 	& 1.0257 & $\epsilon_\lambda\,{=}\,10^{-4}$ & $\epsilon_\mathrm{S}\,{=}\,10^{-6}$, 	& $\alpha_0\,{=}\,10^{-5}$\\
		1.0534 	& 1.0379 	& 1.0608 	& 1.2326 	& 1.1817 	& 1.2871 	& 1.1015 	& 1.1134 & $\epsilon_\lambda\,{=}\,10^{-4}$ & $\epsilon_\mathrm{S}\,{=}\,10^{-5}$, 	& $\alpha_0\,{=}\,10^{-5}$\\
		1.3778 	& 1.1584 	& 1.1723 	& 1.6728 	& 1.4710 	& 1.6967 	& 1.3210 	& 1.3625 & $\epsilon_\lambda\,{=}\,10^{-4}$ & $\epsilon_\mathrm{S}\,{=}\,10^{-4}$, 	& $\alpha_0\,{=}\,10^{-5}$\\
		1.8424 	& 1.4862 	& 1.3964 	& 2.4952 	& 2.0096 	& 2.4233 	& 1.7403 	& 1.8912 & $\epsilon_\lambda\,{=}\,10^{-4}$ & $\epsilon_\mathrm{S}\,{=}\,10^{-3}$, 	& $\alpha_0\,{=}\,10^{-5}$\\
		1.8581 	& 2.0462 	& 1.7307 	& 3.4732 	& 2.6919 	& 3.1800 	& 2.3912 	& 2.6722 & $\epsilon_\lambda\,{=}\,10^{-4}$ & $\epsilon_\mathrm{S}\,{=}\,10^{-2}$, 	& $\alpha_0\,{=}\,10^{-5}$\\
		1.2984 	& 3.1793 	& 2.0990 	& 3.5973 	& 3.2985 	& 3.6414 	& 3.0696 	& 3.1501 & $\epsilon_\lambda\,{=}\,10^{-4}$ & $\epsilon_\mathrm{S}\,{=}\,10^{-1}$, 	& $\alpha_0\,{=}\,10^{-5}$\\
		$\dagger$ 	& 4.9810 	& 2.4551 	& 4.5803 	& 4.1526 	& 3.4432 	& 2.4809 	& $\dagger$ & $\epsilon_\lambda\,{=}\,10^{-4}$ & $\epsilon_\mathrm{S}\,{=}\,1$, 		& $\alpha_0\,{=}\,10^{-5}$\\\hline\\[-1.75ex]
		0.9135 	& 1.6279 	& 1.6687 	& 3.0879 	& 2.3765 	& 2.9968 	& 2.2573 	& 2.5301 & $\epsilon_\lambda\,{=}\,10^{-7}$ & $\epsilon_\mathrm{S}\,{=}\,10^{-2}$, 	& $\alpha_0\,{=}\,0$\\
		0.9135 	& 1.6279 	& 1.7202 	& 3.0879 	& 2.6305 	& 2.9968 	& 2.3726 	& 2.5718 & $\epsilon_\lambda\,{=}\,10^{-6}$ & $\epsilon_\mathrm{S}\,{=}\,10^{-2}$, 	& $\alpha_0\,{=}\,0$\\
		0.9135 	& 1.7760 	& 1.7574 	& 3.1594 	& 2.6975 	& 2.9306 	& 2.3918 	& 2.5504 & $\epsilon_\lambda\,{=}\,10^{-5}$ & $\epsilon_\mathrm{S}\,{=}\,10^{-2}$, 	& $\alpha_0\,{=}\,0$\\
		1.8595 	& 2.0483 	& 1.7574 	& 3.4761 	& 2.6975 	& 3.1837 	& 2.3913 	& 2.6897 & $\epsilon_\lambda\,{=}\,10^{-4}$ & $\epsilon_\mathrm{S}\,{=}\,10^{-2}$, 	& $\alpha_0\,{=}\,0$\\
		1.8595 	& 2.0483 	& 1.7978 	& 3.4956 	& 2.7816 	& 3.2421 	& 2.4593 	& 2.6871 & $\epsilon_\lambda\,{=}\,10^{-3}$ & $\epsilon_\mathrm{S}\,{=}\,10^{-2}$, 	& $\alpha_0\,{=}\,0$\\
		2.6691 	& 2.1029 	& 1.8868 	& 3.4972 	& 2.7986 	& 2.3969 	& 2.4645 	& 2.7558 & $\epsilon_\lambda\,{=}\,10^{-2}$ & $\epsilon_\mathrm{S}\,{=}\,10^{-2}$, 	& $\alpha_0\,{=}\,0$\\\hdashline\\[-1.75ex]
		1      	& 1      	& 1      	& 1      	& 1      	& 1      	& 1      	& 1      & $\epsilon_\lambda\,{=}\,10^{-7}$ & $\epsilon_\mathrm{S}\,{=}\,10^{-2}$, 	& $\alpha_0\,{=}\,10^{-5}$\\
		1      	& 1      	& 1      	& 1.5205 	& 1.6118 	& 1.5275 	& 1.2161 	& 1.8329 & $\epsilon_\lambda\,{=}\,10^{-6}$ & $\epsilon_\mathrm{S}\,{=}\,10^{-2}$, 	& $\alpha_0\,{=}\,10^{-5}$\\
		1      	& 0.0000 	& 1.5940 	& 2.2516 	& 2.2067 	& 0.0000 	& 1.2459 	& 1.5564 & $\epsilon_\lambda\,{=}\,10^{-5}$ & $\epsilon_\mathrm{S}\,{=}\,10^{-2}$, 	& $\alpha_0\,{=}\,10^{-5}$\\
		1.8581 	& 2.0462 	& 1.7307 	& 3.4732 	& 2.6919 	& 3.1800 	& 2.3912 	& 2.6722 & $\epsilon_\lambda\,{=}\,10^{-4}$ & $\epsilon_\mathrm{S}\,{=}\,10^{-2}$, 	& $\alpha_0\,{=}\,10^{-5}$\\
		1.8581 	& 2.0462 	& 1.7964 	& 3.4956 	& 2.7816 	& 3.2418 	& 2.4599 	& 2.6877 & $\epsilon_\lambda\,{=}\,10^{-3}$ & $\epsilon_\mathrm{S}\,{=}\,10^{-2}$, 	& $\alpha_0\,{=}\,10^{-5}$\\
		2.6700 	& 2.1044 	& 1.8871 	& 3.4972 	& 2.7987 	& 2.3967 	& 2.4651 	& 2.7563 & $\epsilon_\lambda\,{=}\,10^{-2}$ & $\epsilon_\mathrm{S}\,{=}\,10^{-2}$, 	& $\alpha_0\,{=}\,10^{-5}$\\\hline\\[-1.75ex]
		6.5835 	& 4.8270 	& 3.7534 	& 5.1328 	& 5.1561 	& 2.8917 	& 3.1875 	& 4.1762 & $\epsilon_\lambda\,{=}\,10^{-2}$ & $\epsilon_\mathrm{S}\,{=}\,1$, 		& $\alpha_0\,{=}\,0$\\
		6.5832 	& 4.8260 	& 3.7535 	& 5.1314 	& 5.1557 	& 2.8905 	& 3.1876 	& 4.1761 & $\epsilon_\lambda\,{=}\,10^{-2}$ & $\epsilon_\mathrm{S}\,{=}\,1$, 		& $\alpha_0\,{=}\,10^{-5}$\\
		\bottomrule[2pt]
	\end{tabular}
\end{table}%

After conducting a comprehensive assessment of these parametric studies, it becomes evident that \textcolor{red}{\emph{variant 2b}} stands out as the most favorable method, showcasing a significant increase in the critical time increment. In contrast to \emph{variant 2c}, no problems with respect to the method's stability when using large values for $\epsilon_\mathrm{S}$ and $\epsilon_\lambda$ are observed. The difference between both variants lies in the fact that in \emph{2b} the mass stabilization matrix is diagonalized by means of the HRZ-method, while in \emph{2c} row-summing is used. Since row-summing can lead to negative components in the lumped matrix, the robustness of the proposed method cannot be guaranteed for all possible cases. Another noteworthy observation is that using a combined EVS--$\alpha$ stabilization scheme does not results in a larger critical time step size. Therefore, it is recommended to employ the EVS-technique only. As for the values of the parameters  $\epsilon_\mathrm{S}$ and $\epsilon_\lambda$, it seems advisable to select a threshold value $\epsilon_\lambda$ of at least $10^{-4}$, while the stabilization constant $\epsilon_\mathrm{S}$ should be at least $10^{-3}$. Keep in mind that these consideration are merely based on the attainable critical time step size and accuracy considerations have not yet entered the list of evaluation criteria (see Sect.~\ref{sec:Examples}).

\begin{figure}[p!]
	\centering
	\subfloat[$p\,{=}\,2$, $\alpha_0\,{=}\,0$]{\includegraphics[scale=1.0]{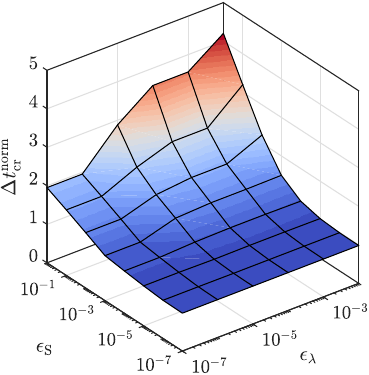}}
	\hspace*{1.5cm}
	\subfloat[$p\,{=}\,2$, $\alpha_0\,{=}\,10^{-5}$]{\includegraphics[scale=1.0]{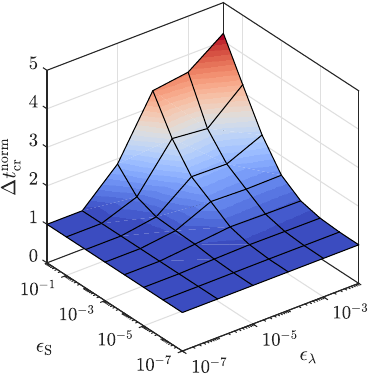}}\\
	\subfloat[$p\,{=}\,5$, $\alpha_0\,{=}\,0$]{\includegraphics[scale=1.0]{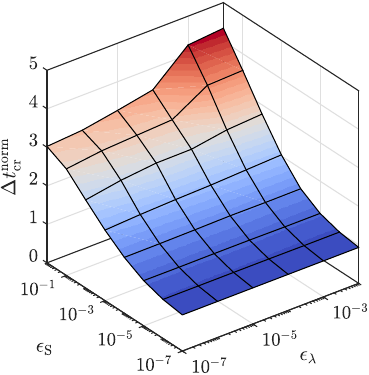}}
	\hspace*{1.5cm}
	\subfloat[$p\,{=}\,5$, $\alpha_0\,{=}\,10^{-5}$]{\includegraphics[scale=1.0]{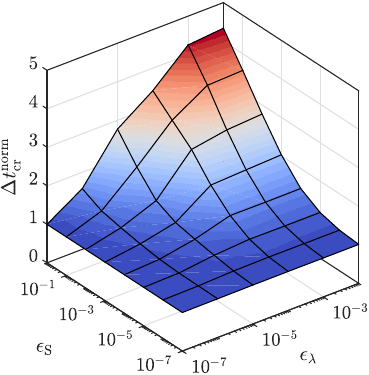}}\\
	\subfloat[$p\,{=}\,8$, $\alpha_0\,{=}\,0$]{\includegraphics[scale=1.0]{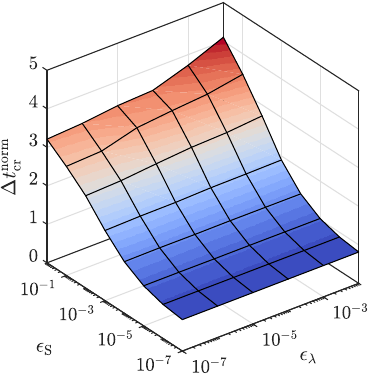}}
	\hspace*{1.5cm}
	\subfloat[$p\,{=}\,8$, $\alpha_0\,{=}\,10^{-5}$]{\includegraphics[scale=1.0]{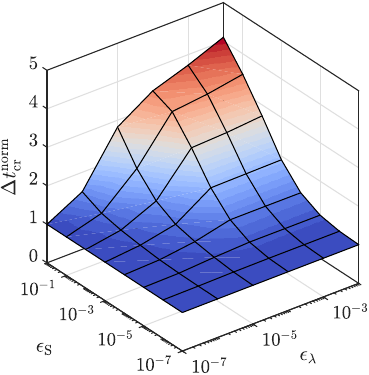}}\\
	\caption{Normalized critical time step size $\Delta t^\mathrm{norm}_\mathrm{cr}$ for \emph{variant 2b}  considering different values of $\epsilon_\mathrm{S}$ and $\epsilon_\lambda$. \label{fig:Variant_10_dtcr}}
\end{figure}

To ensure a comprehensive and conclusive analysis of the default parameters, we have undertaken a parametric study focusing on both $\epsilon_\lambda$ and $\epsilon_\mathrm{S}$. The outcomes of this extensive investigation are visually depicted in Fig.~\ref{fig:Variant_10_dtcr}. This supplementary figure offers a deeper understanding of the relationship between $\epsilon_\mathrm{S}$ and $\epsilon_\lambda$, providing valuable insights to aid in the selection of the most effective stabilization strategy.

In Fig.~\ref{fig:Variant_10_dtcr}, several key observations can be made. For the pure EVS-scheme, depicted in the left column, the critical factor affecting the stability limit is the stabilization parameter $\epsilon_\mathrm{S}$. Remarkably, variations in $\epsilon_\lambda$ within a wide range have minimal influence on the overall outcome, except when $\epsilon_\mathrm{S}$ assumes large values. In such instances, we observe a dependence of the critical time step size on the threshold value $\epsilon_\lambda$.

The behavior is markedly different for the combined EVS--$\alpha$-stabilization. Here, even for smaller values of $\epsilon_\mathrm{S}$, we observe a notable correlation between $\Delta t_\mathrm{cr}$ and $\epsilon_\lambda$. This result aligns with expectations, as the $\alpha$-method stabilizes the entire fictitious domain, causing an increase in eigenvalues. Therefore, a change in the threshold value will also result in the stabilization of more modes. This phenomenon is not encountered in the EVS-technique across a reasonable range of values of $\epsilon_\lambda$. This intriguing result highlights the fact that in the EVS-technique, the stabilized modes are, indeed, closely associated with (almost) zero eigenvalues, thus emphasizing the efficacy of the stabilization approach.

Moreover, it is worth highlighting that the achievable (relative) improvements in critical time step sizes become slightly more pronounced when employing cut elements with high-order shape functions. This suggests that the utilization of the EVS technique is particularly effective with high-order elements, thereby enhancing the efficiency and accuracy of transient simulations.

The numerical findings presented in Fig.~\ref{fig:Variant_10_dtcr} substantiate our initial recommendation to solely employ the EVS technique. The use of the combined EVS--$\alpha$ method does not yield significant advantages. On the contrary, it adversely affects performance, and consequently, we exclude the combined scheme from further assessments.

\section{Condition number}
\label{sec:Kappa}
In the pertinent literature, the term ``ill-conditioning'' is often not well-defined and used rather loosely. Therefore, it is crucial to clarify the specific definition of the condition number used, particularly concerning what aspect of ill-conditioning is being considered. Typically, the condition number regarding matrix inversion, denoted as $\kappa_\mathrm{inv}$, is addressed, which is relevant when solving a linear system of equations (important for implicit time integration schemes). However, it is essential to remember that there is also a condition number related to matrix-vector multiplications, denoted as $\kappa_\mathrm{mvp}$ (important for linear explicit time integration schemes either at the global or element level depending on implementational aspects). Both condition numbers are significant in different scenarios and can be defined for the generic system of equations $\mathbf{A}\mathbf{x} = \mathbf{b}$ as follows:
\begin{alignat}{2}
	\kappa_\mathrm{inv}(\mathbf{A}) & = && \|\mathbf{A}\| \|\mathbf{A}^{-1}\| 
	\shortintertext{and} & && \nonumber \\[-18pt]
	\kappa_\mathrm{mvp}(\mathbf{A}) & = && \cfrac{\|\mathbf{A}\| \|\mathbf{x}\|}{\|\mathbf{A} \mathbf{x}\|} = \cfrac{\|\mathbf{A}\| \|\mathbf{x}\|}{\|\mathbf{b}\|}\,.
\end{alignat}
It is essential to emphasize that ill-conditioning in matrix-vector products is less common, as these operations tend to be more robust and less sensitive to numerical issues compared to matrix inversion.

In this section, the condition numbers of the system matrices related to a square-shaped plate with a circular hole is investigated. The plate's dimensions are $200\,$mm$\times200\,$mm, and the circular void is centered at ($x_\mathrm{o}\,{=}\,100\,$mm, $y_\mathrm{o}\,{=}\,0\,$mm) with a radius of $r\,{=}\,70\,$mm. For the analysis, we assume plane stress conditions and use material properties of steel: $E\,{=}\,210\,$GPa, $\nu\,{=}\,0.3$, and $\rho\,{=}\,7850\,\sfrac{\mathrm{kg}}{\mathrm{m}^3}$.

\begin{figure}[b!]
	\centering
	\subfloat[Mesh\label{fig:PlateHoleMesh}]{\includegraphics[scale=1.0]{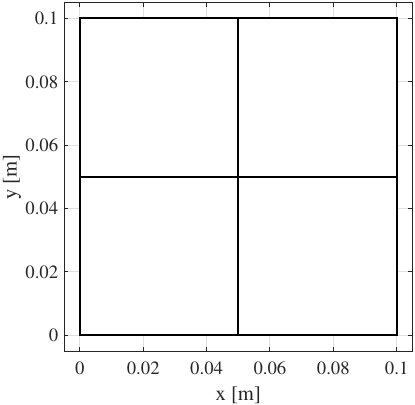}}
	\hspace*{30pt}
	\subfloat[Integration domains\label{fig:PlateHoleInt}]{\includegraphics[scale=1.0]{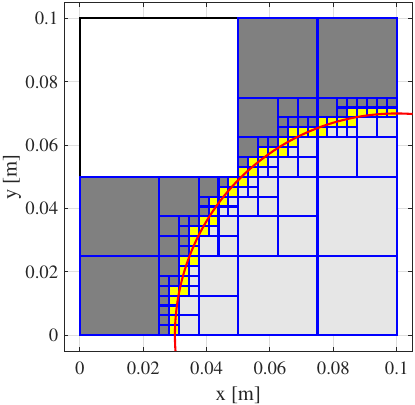}}
	\caption{Model of a rectangular plate with a circular hole (Discretization: $2\times2$ finite element, Subdivision depth: $k\,{=}\,4$; Red solid line: physical boundary of the circular hole). \emph{Color coding for the numerical integration} -- White: conventional finite elements; Dark gray: integration domains located in the physical domain; Light gray: integration domains located in the fictitious domain; Yellow: cut integration domains (leaf cells). \label{fig:PlateHole}}
\end{figure}

To exploit the structure's symmetry, we model only a quarter of the plate and apply symmetry boundary conditions at the right and bottom edges, fixing the displacements in the normal direction. Figure~\ref{fig:PlateHole} illustrates the geometry of the model, including its dimensions, with a spatial discretization of $2\times2$ finite elements. To enhance visualization, we utilize color-coding: The integration domains of cut elements are represented in various shades of gray, while white denotes uncut elements (i.e., conventional spectral elements). Dark gray subdomains correspond to regions within the physical domain, while light gray regions belong to the fictitious domain. Additionally, we use yellow to mark all integration subdomains intersected by the physical boundary (red solid line).

As illustrated in Fig.~\ref{fig:PlateHoleInt}, the mesh for this benchmark problem consists of one spectral element (white) and three cut elements (gray) with varying volume fractions. Notably, the critical cut element is located in the bottom-right corner. Among the cut elements, two exhibit volume fractions of $\chi\,{=}\,73.02\%$, while the last one is barely part of the physical domain, resulting in a tiny volume fraction of $\chi\,{=}\,0.02\%$\footnote{Remark: The reported volume fraction values $\chi$ are computed with a subcell division level of $k\,{=}\,12$. This ensures a highly accurate numerical integration of the element matrices.}. Dealing with such a small volume fraction ($\chi$) poses significant challenges in the numerical analysis, necessitating careful attention when employing such discretizations.

In this context, it is paramount to recognize that the envisioned simulation framework aims for complete automation. Models featuring poorly cut elements are inevitable without expert human interventions. Hence, it is imperative for any immersed boundary method to exhibit reliability and robustness in dealing with such situations. Consequently, the use of specialized and tailor-made stabilization techniques becomes indispensable to ensure the method's efficacy and accuracy even under challenging scenarios with extremely small volume fractions.

Before delving into the numerical results (see Sect.~\ref{sec:Examples}) and assessing the convergence properties, we want to confirm that the proposed stabilization approach effectively resolves the ill-conditioning issues (related to matrix inversion) as previously discussed in Refs.~\cite{ArticleLoehnert2014, ArticleGarhuom2022} for static cases. The results of our tests are compiled in Tables~\ref{tab:KappaPlateHoleK} and \ref{tab:KappaPlateHoleM} for the stiffness and mass matrices, respectively.
\begin{table}[b!]
	\caption{Condition number $\kappa^\mathrm{K}_\mathrm{inv}$ of the stiffness matrix $\mathbf{K}$ -- plate with a circular hole ($2\times2$ finite elements; $k\,{=}\,4$).\label{tab:KappaPlateHoleK}}
	\centering
	\begin{tabular}{ccccccc}
		\toprule[2pt]
		Method: $\rightarrow$ & Unstabilized & \multicolumn{2}{c}{$\alpha$-Method} & EVS & \multicolumn{2}{c}{Combined EVS--$\alpha$-method} \\
		$p$: $\downarrow$ & $\alpha_0\,{=}\,0$ & $\alpha_0\,{=}\,10^{-12}$ & $\alpha_0\,{=}\,10^{-5}$ & $\epsilon_\mathrm{S}\,{=}\,10^{-2}$, $\epsilon_\lambda\,{=}\,10^{-3}$ & $\alpha_0\,{=}\, 10^{-12}$ & $\alpha_0\,{=}\, 10^{-5}$\\
		\midrule[1pt]
		1  & 2.64E+02 & 2.64E+02 & 2.64E+02 & 2.64E+02 & 2.64E+02 & 2.64E+02 \\
		2  & 3.29E+23 & 1.08E+13 & 1.08E+06 & 1.51E+03 & 1.65E+03 & 1.08E+04 \\
		3  & 4.12E+26 & 2.77E+13 & 2.74E+06 & 3.21E+03 & 3.35E+03 & 2.51E+04 \\
		4  & 1.25E+25 & 6.96E+13 & 6.96E+06 & 1.10E+04 & 9.47E+03 & 1.53E+04 \\
		5  & 1.87E+25 & 1.22E+14 & 1.27E+07 & 2.45E+04 & 1.81E+04 & 2.59E+04 \\
		6  & 5.05E+24 & 2.09E+14 & 2.59E+07 & 1.94E+04 & 2.48E+04 & 1.93E+04 \\
		7  & 6.47E+25 & 3.17E+14 & 4.27E+07 & 2.57E+04 & 2.77E+04 & 2.66E+04 \\
		8  & 2.81E+26 & 4.64E+14 & 6.34E+07 & 3.10E+04 & 3.52E+04 & 3.68E+04 \\
		9  & 1.53E+27 & 6.42E+14 & 9.35E+07 & 4.23E+04 & 4.27E+04 & 4.05E+04 \\
		10 & 1.12E+26 & 8.68E+14 & 1.32E+08 & 5.09E+04 & 4.83E+04 & 5.12E+04 \\
		\bottomrule[2pt]
	\end{tabular}
	%
	\vspace*{3pt}
	%
	\caption{Condition number $\kappa^\mathrm{M}_\mathrm{inv}$ of the (lumped) mass matrix $\mathbf{M}$ -- plate with a circular hole ($2\times2$ finite elements; $k\,{=}\,4$).\label{tab:KappaPlateHoleM}}
	\centering
	\begin{tabular}{ccccccc}
		\toprule[2pt]
		Method:  $\rightarrow$ & Unstabilized & \multicolumn{2}{c}{$\alpha$-Method} & EVS & \multicolumn{2}{c}{Combined EVS--$\alpha$-method} \\
		$p$: $\downarrow$ & $\alpha_0\,{=}\,0$ & $\alpha_0\,{=}\,10^{-12}$ & $\alpha_0\,{=}\,10^{-5}$ & $\epsilon_\mathrm{S}\,{=}\,10^{-2}$, $\epsilon_\lambda\,{=}\,10^{-3}$ & $\alpha_0\,{=}\, 10^{-12}$ & $\alpha_0\,{=}\, 10^{-5}$\\
		\midrule[1pt]
		1  & 8.43E+00 & 8.43E+00 & 8.43E+00 & 8.43E+00 & 8.43E+00 & 8.43E+00 \\
		2  & 3.59E+10 & 3.57E+10 & 5.69E+05 & 8.82E+01 & 8.82E+01 & 5.45E+04 \\
		3  & 4.26E+11 & 3.99E+11 & 6.06E+05 & 3.43E+02 & 3.43E+02 & 1.04E+05 \\
		4  & 4.19E+10 & 4.17E+10 & 9.05E+05 & 2.38E+02 & 2.38E+02 & 1.28E+04 \\
		5  & 1.88E+11 & 1.85E+11 & 1.12E+06 & 1.47E+02 & 1.47E+02 & 1.03E+04 \\
		6  & 2.64E+10 & 2.63E+10 & 1.39E+06 & 1.15E+02 & 1.15E+02 & 1.51E+04 \\
		7  & 8.29E+10 & 8.25E+10 & 1.70E+06 & 1.33E+02 & 1.33E+02 & 1.40E+04 \\
		8  & 2.71E+10 & 2.71E+10 & 1.98E+06 & 2.22E+02 & 2.22E+02 & 1.09E+04 \\
		9  & 2.24E+10 & 2.23E+10 & 2.32E+06 & 2.18E+02 & 2.18E+02 & 1.55E+03 \\
		10 & 3.91E+10 & 3.91E+10 & 2.88E+06 & 2.73E+02 & 2.73E+02 & 1.66E+03 \\
		\bottomrule[2pt]
	\end{tabular}
\end{table}%

Let us begin by discussing the results obtained using the $\alpha$-method. It is evident that by employing larger values for the stabilization parameter $\alpha_0$, both condition numbers $\kappa^\mathrm{K}_\mathrm{inv}$ and $\kappa^\mathrm{M}_\mathrm{inv}$ experience a significant reduction. However, it is crucial to be mindful that this stabilization approach affects the entire fictitious domain contribution of a cut element and not just specific modes. If $\alpha_0$ is chosen too large, a substantial amount of stiffness and mass is added to the original system. Despite this, the $\alpha$-method achieves a significant reduction in the condition numbers, reaching six to seven orders of magnitude improvement when compared to the hardly stabilized case with $\alpha_0\,{=}\,10^{-12}$, which is used as a reference. It is noteworthy that the mass matrix $\kappa^\mathrm{M}_\mathrm{inv}$ exhibits reasonably low condition numbers even without any stabilization, owing to the application of a nodal quadrature lumping technique for spectral elements and the HRZ-method for cut elements.

Next, we consider the results obtained by employing only the EVS-technique, with a stabilization parameter $\epsilon_\mathrm{S}\,{=}\,10^{-2}$ and an eigenvalue stabilization threshold $\epsilon_\lambda\,{=}\,10^{-3}$. Here, we observe a significantly increased performance in reducing the condition numbers. An additional three to four orders of magnitude are achieved. Note again that the EVS-technique targets specific modes and results in a less intrusive adaptation of the original system, as previously mentioned.

As discussed before, the combined EVS--$\alpha$-method does not lead to improved results. In contrast, for large values of $\alpha_0$ the effectiveness of the pure EVS-scheme is disturbed and thus, the conditioning of the system matrices worsens. This might initially seem counterintuitive, but it can be explained by the fact that a larger $\alpha_0$ results in a reasonably stable element formulation and thus, the EVS-technique has less impact. To conclude, in terms of the conditioning of the system of equations, it is advisable to employ the pure EVS-technique.

In examining the condition number of the stiffness/mass matrix with respect to matrix-vector products ($\kappa^\square_\mathrm{mvp}$), which is particularly relevant in explicit dynamics, we observe that it remains unaffected by the stabilization scheme and consistently stays around a value of $10^2$ for this specific test case. This suggests that we can expect only a negligible loss in accuracy when considering matrix-vector products. These findings reinforce the notion that stabilizing the stiffness matrix is unnecessary in the context of explicit dynamics as no numerical penalties are incurred.

The observed low values of $\kappa^\square_\mathrm{mvp}$ suggest that the proposed stabilization schemes have minimal impact on the efficiency and accuracy of matrix--vector operations. Therefore, we can confidently assert that in explicit dynamics simulations, concentrating on stabilizing the mass matrix is adequate, and stabilizing the stiffness matrix is not meaningful in this specific context. Instead, our attention should be directed towards other aspects, such as the attainable accuracy.
\section{Numerical examples}
\label{sec:Examples}
In this section, we thoroughly examine two benchmark examples to evaluate the effectiveness of the proposed stabilization schemes concerning both computational efficiency and achievable numerical accuracy. Throughout the analysis (explicit dynamics), we employ \emph{variant 2b}. For the sake of completeness, we list the chosen default settings:
\begin{enumerate}
	\item The mass matrix is stabilized.
	\item The stiffness matrix is \textbf{not} stabilized.
	\item The consistent mass matrix (CMM) is used for the eigenvalue decomposition.
	\item Both the mass matrix $\mathbf{M}_c$ and the mass stabilization matrix $\mathbf{M}^\mathrm{S}_c$ are lumped by means of the HRZ-method.
	\item The threshold for stabilizing mode shapes is $\epsilon_\lambda\,{\approx}\,10^{-3}$ with $\epsilon_\lambda \,\in\, [10^{-4},\, 10^{-2}]$.
	\item The stabilization parameter is $\epsilon_\mathrm{S}\,{\approx}\,10^{-3}$ with $\epsilon_\mathrm{S} \,\in\, [10^{-5},\, 10^{-2}]$.
\end{enumerate}
By adopting these standardized settings, our primary goal is to maintain consistency and enable a fair comparison of results across various problems in this article. Note that a combination of the EVS-technique with the $\alpha$-method has been ruled out since no advantages in terms of the critical time step sizes are gained.

Additionally, our focus in this study is exclusively on transient problems, which we will investigate through carefully chosen benchmark examples. It is worth noting that implementing the iterative correction scheme proposed for static problems in Ref.~\cite{ArticleGarhuom2022} is currently not a feasible option. Doing so would necessitate a substantial re-writing effort for existing time integration schemes, a task that falls beyond the scope of this paper. However, we anticipate addressing this aspect in forthcoming publications dedicated to this subject.

Furthermore, it is essential to acknowledge the inherent complexities associated with solving eigenvalue problems within the context of immersed boundary methods. These challenges manifest in the form of a significant number of spurious modes observed in structures discretized with only a few elements. Such issues arise due to constraints related to degrees of freedom with limited support \cite{ArticleEisentraeger2022a} and ill-conditioning problems. The extent to which these effects are prevalent in the solution largely depends on the employed stabilization technique.

On the other hand, in cases involving structures containing only a few cut elements, where the majority of elements is entirely located within the physical domain, problems are not expected when computing only the first few eigenvalues. Nevertheless, it is important to emphasize that addressing eigenvalue problems, such as modal or buckling analyses, requires dedicated research endeavors and the formulation of tailored solution strategies. The utilization of standard eigenvalue solvers presents challenges in distinguishing spurious modes from physically meaningful ones. Consequently, despite the significance of modal analyses in dynamics, this contribution does not delve into this specific topic.
\subsection{Wave propagation in a simple waveguide---rectangular plate}
\label{sec:TransientRect}
\begin{figure}[t!]
	\centering
	\begin{overpic}[scale=1.0]{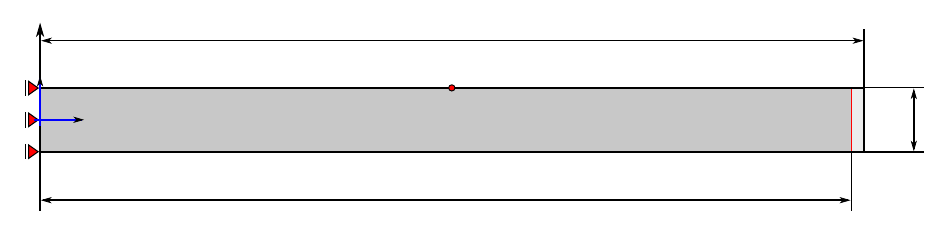}
		\put(3,23){$F(t)$}
		\put(49.5,4.75){$l_\mathrm{phys}$}
		\put(50,21.5){$l_\mathrm{d}$}
		\put(95.3,11.6){\rotatebox{90}{$h$}}
		\put(9.25,11.7){$x$}
		\put(4.75,17){$y$}
		\put(47,16.65){$P_1$}
		\put(93.4,6){$\Omega_\mathrm{fict}$}
		\put(93.15,7.6){\color{blue}\line(-1,2.5){1.75}}
		\put(91.4,11.975){\color{blue}\circle*{1.5pt}}
		\put(70,11.6){$\Omega_\mathrm{phys}$}
	\end{overpic}
	\caption{Model of the rectangular waveguide.\label{fig:RectangularPlateModel}}
\end{figure}
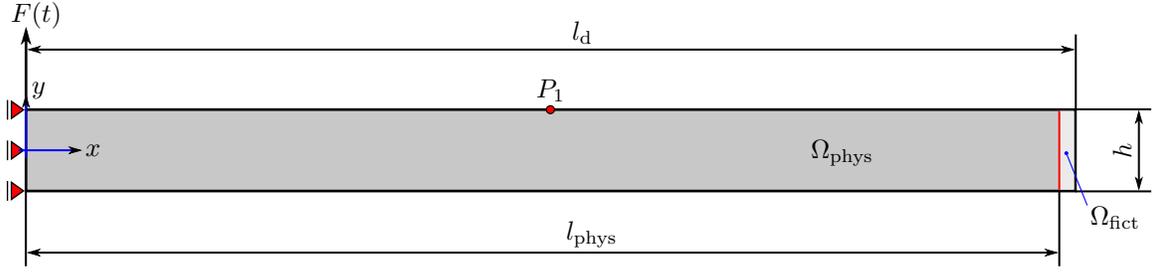

In the first examples, a simple rectangular domain is chosen in order to study the propagation of guided waves. A sketch of the waveguide along with its dimensions and Dirichlet as well as Neumann boundary conditions is depicted in Fig.~\ref{fig:RectangularPlateModel}. For this example, it is assumed that the waveguide consists of aluminum and therefore, the material properties are Young's modulus $E\,{=}\,70\,$GPa, Poisson's ratio $\nu\,{=}\,0.3$, and the mass density $\rho\,{=}\,2700\,\sfrac{\mathrm{kg}}{\mathrm{m}^3}$. Note that plane strain conditions are utilized to reduce the dimensionality of the problem from 3 to 2. This also means that the plate has an infinite length in $z$-direction. 

Considering the geometrical dimensions, the (extended) length of the plate (including the fictitious domain) is set to $l_\mathrm{d}\,{=}\,200\,$mm, while a thickness of $h\,{=}\,2\,$mm is selected. Hence, the length of the fictitious part 
\begin{equation}
	l_\mathrm{fict} = l_\mathrm{d} - l_\mathrm{phys}
\end{equation}
is directly related to the volume fraction $\chi$ of the cut elements
\begin{equation}
	\chi = 1 - \cfrac{l_\mathrm{fict}}{h_\mathrm{l}}\,,
\end{equation}
determining the severity of the cut. In this example, the volume fraction is selected as $\chi\,{=}\,0.05$, meaning that only 5\% of the cut elements should be located in the physical domain. Since we divide the plate into $n^\mathrm{x}_\mathrm{Elem}\,{=}\,200$ finite elements along the longitudinal direction and $n^\mathrm{y}_\mathrm{Elem}\,{=}\,2$ finite elements in the thickness direction, the element size is $h_\mathrm{l}\,{=}\,1\,$mm. This also means that the fictitious part of the plate is of length $l_\mathrm{fict}\,{=}\,0.95\,$mm and consequently, $l_\mathrm{phys}\,{=}\,199.05\,$mm

Note that an observation point, denoted as $P_1$ (see Fig.~\ref{fig:RectangularPlateModel}), is introduced on the top surface of the plate at ($x\,{=}\,100\,$mm, $y\,{=}\,\sfrac{h}{2}\,{=}\,1\,$mm) to evaluate the wave signals. To make sure that a node is always located at the observation point, the geometric model is partitioned at  $x\,{=}\,100\,$mm.

Based on the spatial discretization described above, a \emph{p}-refinement study is conducted in this section, where the numerical integration over the cut elements is performed by means of the quadtree-based composed integration scheme. Here, the number of subdivision steps for the quadtree-algorithm is set to $k\,{=}\,8$.

\begin{figure}[t!]
	\centering
	\subfloat[Time-domain signal: $p(t)$]{\includegraphics[scale=1.0]{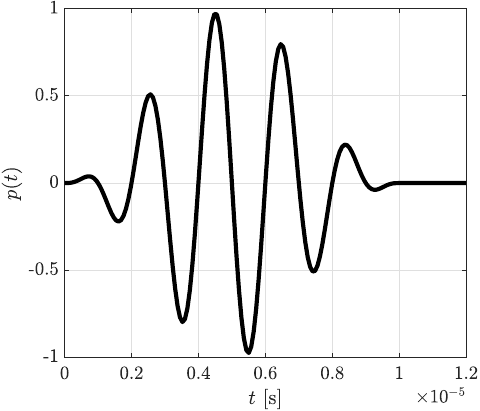}}\hfill
	\subfloat[Frequency spectrum: $P(f)$]{\includegraphics[scale=1.0]{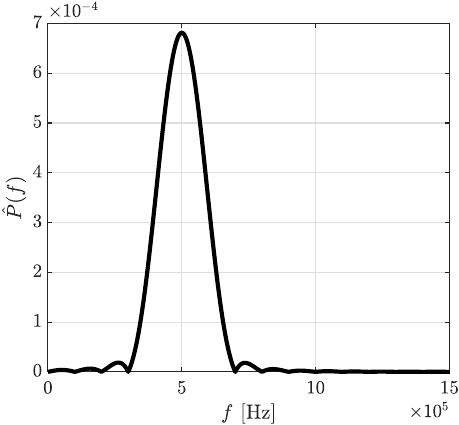}}
	\caption{Excitation signal -- Sine burst modulated by a Hann-window. \label{fig:ExciteSignal}}
\end{figure}

The waves are excited at the top surface of the structure by means of a time-dependent (concentrated) force $F(t)$ acting at one node only, which is defined as
\begin{equation}
	F(t) = \bar{F} p(t) \quad\text{with}\quad p(t) = \sin(\omega t) \sin^2\left(\cfrac{\omega t}{2n}\right) \quad\forall\,t \in \left[0\,,\cfrac{n}{f_\mathrm{c}}\right] \,,\label{eq:Hann}
\end{equation}
where $\bar{F}$ denotes the amplitude of the external excitation force, $\omega\,{=}\,2\pi f_\mathrm{c}$ is the (central) angular excitation frequency, and $n$ stands for the number of periods. In the current example, the following values have been selected: $\bar{F}\,{=}\,10^6\,$N, $f_\mathrm{c}\,{=}\,500\,$kHz, and $n\,{=}\,5$. The time-domain signal and its frequency spectrum (according to these values) are plotted in Fig.~\ref{fig:ExciteSignal}. Due to the comparably narrow-banded excitation, the spatial discretization can be easily set up according to the guidelines developed in Ref.~\cite{ArticleWillberg2012} based on the value of $f_\mathrm{c}$. The wavelengths of the traveling waves need to be derived from dispersion diagrams for Lamb waves \cite{BookGiurgiutiu2008}, featuring the wave velocity as a function of the excitation frequency. For a simpler approach to designing an analysis-suitable mesh, the bulk wave velocities
\begin{alignat}{2}
	c_\mathrm{P} & = && \sqrt{\cfrac{\lambda + 2\mu}{\rho}} \label{eq:CP} \\
	\shortintertext{and} & && \nonumber \\[-18pt]
	c_\mathrm{S} & = && \sqrt{\cfrac{\mu}{\rho}} \label{eq:CS}
\end{alignat}
are used instead of the values for guided waves (Lamb waves). Here, $c_\mathrm{P}$ denotes the pressure wave and $c_\mathrm{S}$ is the shear wave velocity. Only the two Lam\'e constants $\lambda$ and $\mu$ (shear modulus) need to be known, which can be computed from the engineering constants $E$ and $\nu$ by
\begin{alignat}{2}
	\lambda & = && \cfrac{E\nu}{(1+\nu)(1-2\nu)} \label{eq:Lame1} \\
	\shortintertext{and} & && \nonumber \\[-18pt]
	\mu & = && \cfrac{E}{2(1+\nu)}\,. \label{eq:Lame2}
\end{alignat}

In the general case, the shear wave velocity $c_\mathrm{S}$ is taken to determine the element size as a multiple of the shear wavelength $\lambda_\mathrm{S}$, which is defined as
\begin{equation}
	\lambda_\mathrm{S} = \cfrac{c_\mathrm{S}}{f_\mathrm{c}}\,.
\end{equation}
Depending on the polynomial degree of the shape functions, at least 3 to 20 nodes per wavelength should be employed \cite{ArticleWillberg2012}; the higher the polynomial degree, the fewer nodes per wavelength are required. Considering aluminum as the plate's material, the pressure and shear wave velocities are $5907.6\,\sfrac{\mathrm{m}}{\mathrm{s}}$ and $3157.8\,\sfrac{\mathrm{m}}{\mathrm{s}}$, respectively. Consequently, the wavelength of a shear wave at $500\,$kHz is $6.32\,$mm. As mentioned above, the rectangular waveguide is discretized by 200 spectral elements/cells in $x$-direction and therefore, the element size is $1\,$mm, resulting in slightly more than six elements per shear wavelength. This discretization is indeed notably too coarse for low-order elements, but sufficiently accurate results are expected for $p\,{\ge}\,3$ \cite{ArticleWillberg2012}.

The results are evaluated at the observation point $P_1$, where the error between a suitable reference solution of the displacement history in $x$-direction, denoted as $u^{P_1}_\mathrm{x}$, and its numerical approximation, denoted as $\tilde{u}^{P_1}_\mathrm{x}$, is computed in the $L_2$-norm as
\begin{equation}
	e_{L_2} = \sqrt{\cfrac{\displaystyle\int\limits_{t=0}^{t_\mathrm{sim}}\Bigl(u^{P_1}_\mathrm{x}(t) - \tilde{u}^{P_1}_\mathrm{x}(t)\Bigr)^2 \mathrm{d}t}{\displaystyle\int\limits_{t=0}^{t_\mathrm{sim}}\Bigl(u^{P_1}_\mathrm{x}(t)\Bigr)^2 \mathrm{d}t}} \times 100 [\%]\,.
	\label{eq:Error_L2_ux_P1}
\end{equation}

All explicit simulations are performed using a lumped mass matrix formulation and the CDM discussed in Sect.~\ref{sec:StabTimeInt}. To ensure that the temporal discretization does not interfere with the spatial one, a time step size of $\Delta t\,{=}\,3\times10^{-9}\,$s is chosen. This value is below the critical one for all simulations and corresponds to 667 sampling points per central frequency $f_\mathrm{c}$ guaranteeing a highly-accurate time integration. The overall simulation time is set to $t_\mathrm{sim}\,{=}\,120\,\mu$s and therefore, $40{,}000$ time steps are computed.

\begin{figure}[b!]
	\centering
	\subfloat[$u_\mathrm{x}$ at $P_1$]{\begin{overpic}[scale=1.0]{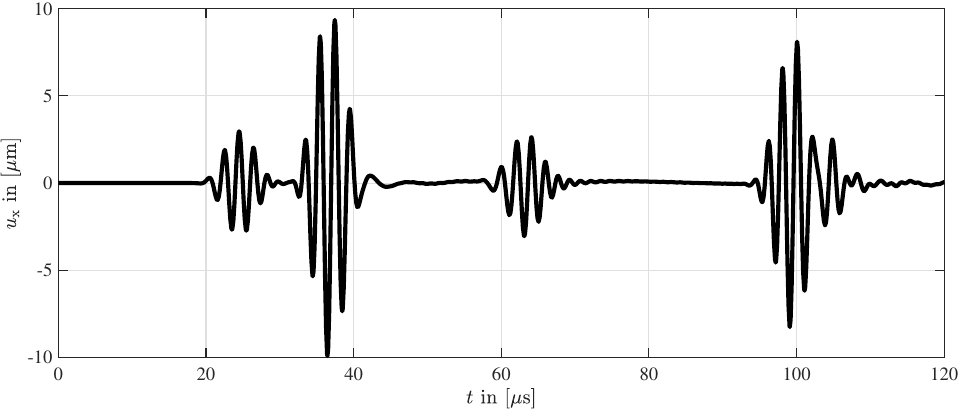}
	\put(22,31){$S_0$-mode}
	\put(52,30){\textcolor{blue}{$S_0$-mode}}
	\put(37,8){$A_0$-mode}
	\put(86,9){\textcolor{blue}{$A_0$-mode}}
	\end{overpic}}\\
	\subfloat[$u_\mathrm{y}$ at $P_1$]{\includegraphics[scale=1.0]{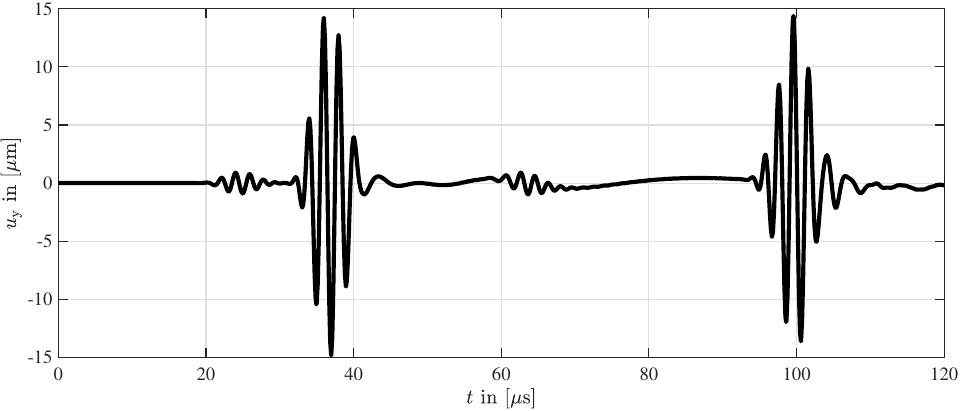}}
	\caption{Time-history of the displacement components at the observation point $P_1$ -- rectangular waveguide, $\chi\,{=}\,5\%$. \label{fig:TimeHistP1_rectWG}}
\end{figure}

The attained error will be evaluated against a numerical \emph{overkill} solution computed using a geometry-conforming mesh consisting of $199\times2$ spectral elements of order $p\,{=}\,10$ ($n_\mathrm{DOF}\,{=}\,83{,}622$), whereas the time integration has been performed using a high-order accurate Pad\'e-based (implicit) time integration scheme \cite{ArticleSong2022a, ArticleSong2022b}. The spectral elements are lumped by means of nodal quadrature. However, row-summing and diagonal scaling (HRZ-method) would yield identical results for Cartesian meshes \cite{ArticleDuczek2019b}. The parameters of the Pad\'e-scheme are selected such that no numerical (algorithmic) damping ($\rho_\infty\,{=}\,1$) is introduced during the simulation and that a time integrator of order $p_\mathrm{Int}\,{=}\,10$ is employed. The prescribed time step size, is also fixed to $\Delta t\,{=}\,3\times10^{-9}\,$s. The reference histories of the displacements in $x$- and $y$-directions for point $P_1$ are depicted in Fig.~\ref{fig:TimeHistP1_rectWG} for comparison purposes. Especially from the displacement history in $x$-direction, we can clearly distinguish the incident wave packets and their reflections at the right end of the rectangular waveguide. The first two wave packets correspond to the incident $S_0$- and $A_0$-modes, while the next wave packets are their reflections at the free boundary. 

\begin{figure}[t!]
	\centering
	\subfloat[Error $e_{L_2}$]{\includegraphics[scale=1.0]{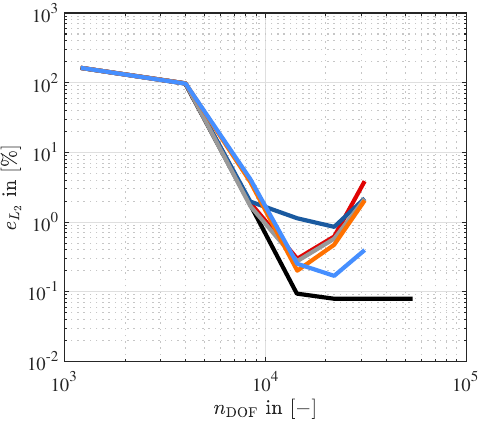}}\hfill
	\subfloat[Critical time step $\Delta t_\mathrm{cr}$]{\includegraphics[scale=1.0]{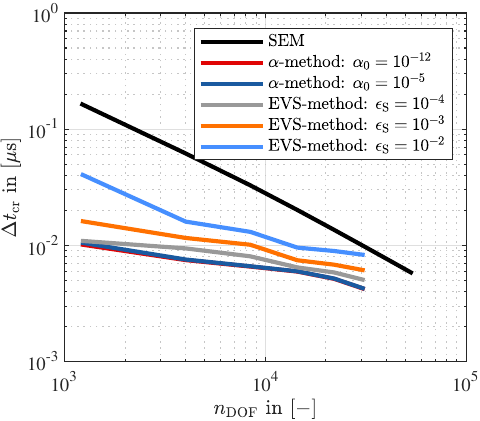}}
	\caption{Evolution of the error and the critical time step size for different stabilization parameters $\alpha_0$ and $\epsilon_\mathrm{S}$ on a fixed mesh ($200\times2$ spectral elements) under \emph{p}-refinement with $p\in[1,6]$ -- rectangular waveguide, $\chi\,{=}\,5\%$. In the EVS-technique, the stabilization threshold is set to $\epsilon_\lambda\,{=}\,10^{-3}$. \label{fig:RectangularWaveguide}}
%
\vspace*{6pt}
%
	\centering
	\includegraphics[scale=1.0]{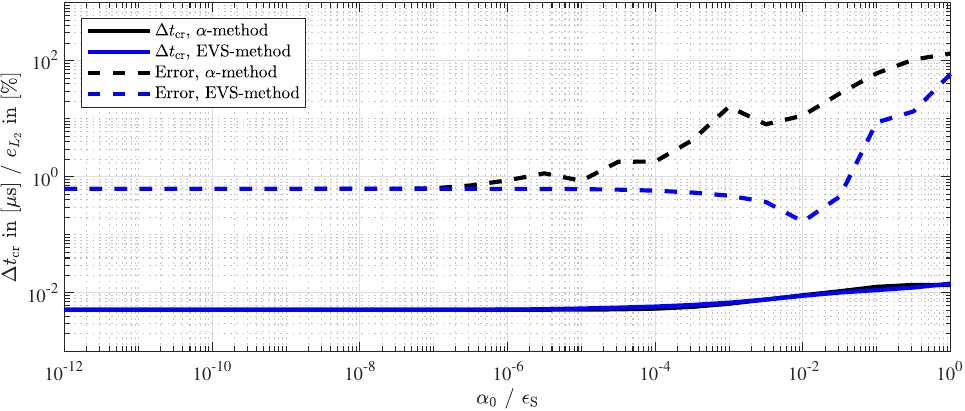}
	\caption{Evolution of the error and the critical time step size for different stabilization parameters $\alpha_0$ and $\epsilon_\mathrm{S}$ on a fixed mesh ($200\times2$ spectral elements with $p\,{=}\,5$, $n_\mathrm{DOF}\,{=}\,22{,}022$) -- rectangular waveguide, $\chi\,{=}\,5\%$. In the EVS-technique, the stabilization threshold is set to $\epsilon_\lambda\,{=}\,10^{-3}$. \label{fig:RectangularWaveguide_p5}}
\end{figure}

In Figs.~\ref{fig:RectangularWaveguide} and \ref{fig:RectangularWaveguide_p5}, the numerical results in terms of the attainable error and critical time step size are depicted for the rectangular waveguide example. To facilitate the comparison of the different stabilization methods, we also included the results of SEM simulations as a reference.

\begin{figure}[b!]
	\centering
	\subfloat[Error $e_{L_2}$]{\includegraphics[scale=1.0]{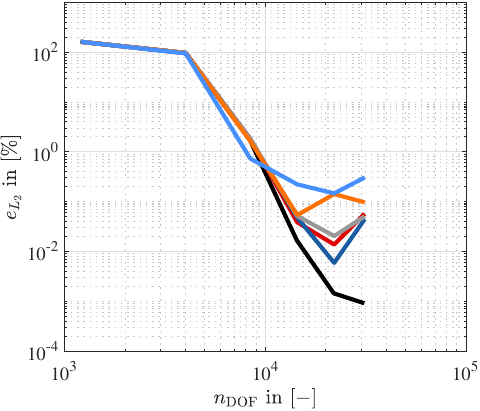}}\hfill
	\subfloat[Critical time step $\Delta t_\mathrm{cr}$]{\includegraphics[scale=1.0]{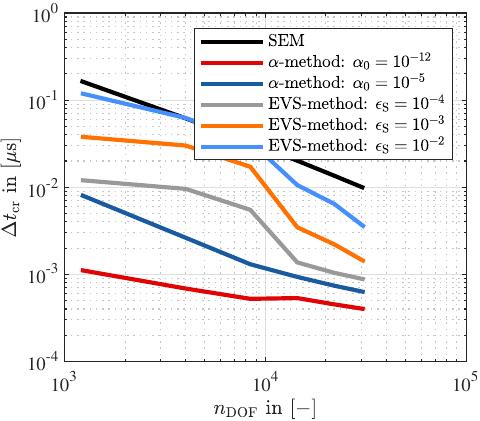}}
	\caption{Evolution of the error and the critical time step size for different stabilization parameters $\alpha_0$ and $\epsilon_\mathrm{S}$ on a fixed mesh ($200\times2$ spectral elements) under \emph{p}-refinement with $p\in[1,6]$ -- rectangular waveguide, $\chi\,{=}\,0.5\%$. In the EVS-technique, the stabilization threshold is set to $\epsilon_\lambda\,{=}\,10^{-3}$. \label{fig:RectangularWaveguide2}}
%
\vspace*{6pt}
%
	\centering
	\includegraphics[scale=1.0]{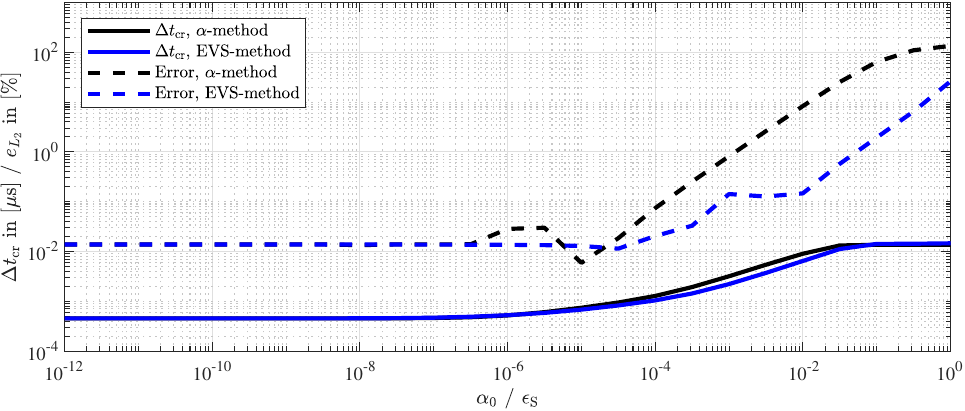}
	\caption{Evolution of the error and the critical time step size for different stabilization parameters $\alpha_0$ and $\epsilon_\mathrm{S}$ on a fixed mesh ($200\times2$ spectral elements with $p\,{=}\,5$, $n_\mathrm{DOF}\,{=}\,22{,}022$) -- rectangular waveguide, $\chi\,{=}\,0.5\%$. In the EVS-technique, the stabilization threshold is set to $\epsilon_\lambda\,{=}\,10^{-3}$. \label{fig:RectangularWaveguide2_p5}}
\end{figure}

It can be observed that the accuracy of the SEM results is limited by the accuracy of the CDM, which already highlights one problem when using high-order spatial discretizations. We also require a matching high-order (explicit) time integration scheme to balance the spatial and temporal discretization errors. Despite the selection of an exceedingly small time increment for the simulations, the second-order accurate time integration scheme struggles to match the precision of the high-order SEM. By comparing the achieved error plateaus in Figs.~\ref{fig:RectangularWaveguide} and \ref{fig:RectangularWaveguide2}, we clearly observe the effect of decreasing the used time step size by a factor of $10$. This drastic decrease becomes necessary to ensure that all simulations can be run with the same time step size even for very critical cuts, where the elements hardly belongs to the physical domain anymore. However, the difference in accuracy is striking: while in the former case, the attainable error remains confined to approximately $0.1\%$, employing a time increment that is ten times smaller reduces the error level to a remarkable $0.001\%$, representing an improvement of nearly two orders of magnitude in terms of accuracy.

In our analysis of the numerical results obtained using immersed techniques and various stabilization schemes, we observe that the achievable accuracy is similar independent of whether the $\alpha$-method or the EVS-technique is employed during the simulations. However, when we shift our focus to the critical time step size, a noteworthy difference is noted. The EVS-technique exhibits superior performance compared to the $\alpha$-method in this regard. This advantage lies in the flexibility it offers with respect to the stabilization parameter $\epsilon_\mathrm{S}$. Unlike the $\alpha_0$ parameter, larger values of $\epsilon_\mathrm{S}$ can be chosen without compromising the accuracy of the results. Depending on the polynomial order used for spatial discretization, we observe a substantial improvement in the critical time step size. Specifically, for this example, the critical time step size increases by a factor of approximately $2.6$. This finding underscores the practical advantages of the EVS-technique, particularly in scenarios where managing the critical time step size is a critical consideration for time-dependent simulations.

Additional insights can be gained form Fig.~\ref{fig:RectangularWaveguide_p5}, where the spatial discretization is fixed ($200\times2$ spectral elements with $p\,{=}\,5$, $n_\mathrm{DOF}\,{=}\,22{,}022$), while systematically varying the stabilization parameters $\alpha_0$ and $\epsilon_{\mathrm{S}}$ over a wide range. The important point to notice is that the behavior of the critical time step size as it evolves with increasing stabilization parameters remains strikingly similar for both schemes. However, a noteworthy distinction is that the EVS-method consistently exhibits superior accuracy at higher stabilization values. This distinction is evident in the error curve for the EVS-technique, which is shifted towards larger values of $\epsilon_\mathrm{S}$.  Therefore, we recommend $\epsilon_\mathrm{S}\,{>}\,\alpha_0$ and still achieve a similar level of accuracy. This advantage translates into a larger time step size and, consequently, shorter simulation times. Therefore, the numerical findings presented in this section corroborate the observations discussed in Sect.~\ref{sec:EVST_dtcr}, where it was previously noted that the EVS-scheme excels in increasing the critical time step size.

In the remainder of this section, we revisit our analysis, this time focusing on a more challenging cut scenario where the volume fraction is further reduced to $\chi\,{=}\,0.5\%$. In practical terms, this reduction translates to a decrease in the length of the physical domain, now measuring $l_\mathrm{phys}\,{=}\,199.005\,$mm.

The displacement histories in both the $x$- and $y$-directions are not presented for this particular example. The reason being that they closely resemble those depicted in Fig.~\ref{fig:TimeHistP1_rectWG} for the slightly longer waveguide. However, it is important to note that in the case of such an extreme cut, we anticipate witnessing an even more pronounced demonstration of the EVS-technique's effectiveness. This improved performance can be attributed to the targeted nature of the EVS-approach.

In order to address this rather severe cut, it becomes imperative to reduce the time step size by an order of magnitude, bringing it down to $\Delta t\,{=}\,3\times10^{-10}\,$s. This adjustment ensures that all simulations,  regardless of the polynomial order ($p$), can be executed using the same time increment. The numerical results obtained with the decreased time step size are presented in Figs.~\ref{fig:RectangularWaveguide2} and \ref{fig:RectangularWaveguide2_p5}. In these figures, we observe a behavior akin to what was previously documented in Figs.~\ref{fig:RectangularWaveguide} and \ref{fig:RectangularWaveguide_p5}, with one notable distinction. The range within which the critical time step size can be enhanced expands significantly. While in the initial example with a volume fraction of $\chi\,{=}\,5\%$, we observed a critical time step size increase of approximately $2.6$ times, this figure improves to $32.3$ for the severely cut case with a volume fraction of $\chi\,{=}\,0.5\%$. This outcome reaffirms the notion that stabilization schemes deliver their best performance when confronted with critically cut elements.
\subsection{Wave propagation in a porous waveguide---perforated plate}
\label{sec:TransientCircHoles}
\begin{figure}[b!]
	\centering
	\begin{overpic}[scale=1.0]{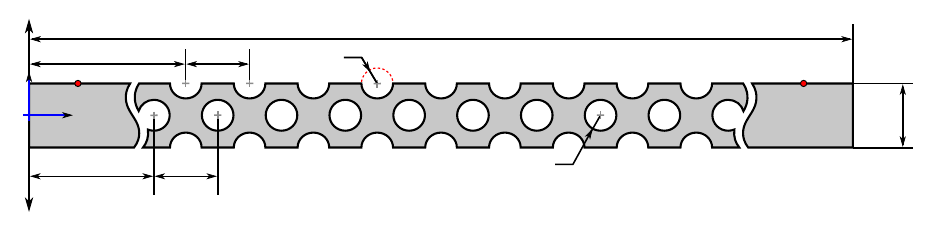}
		\put(2,23.5){$F(t)$}
		\put(2,0.25){$F(t)$}
		\put(9.6,6.9){$l_\mathrm{c}$}
		\put(18.25,6.9){$\Delta l_\mathrm{c}$}
		\put(11.95,18.75){$l_\mathrm{u}$}
		\put(21.65,18.75){$\Delta l_\mathrm{u}$}
		\put(50,21.5){$l_\mathrm{d}$}
		\put(94.35,12.1){\rotatebox{90}{$h$}}
		\put(8.25,12.15){$x$}
		\put(1.25,17.25){$y$}
		\put(59.25,8.0){$r_\circ$}
		\put(36.5,19.5){$r_\cup$}
		\put(6,13.75){$P_1$}
		\put(85.5,13.75){$P_2$}
	\end{overpic}
	\caption{Model of the perforated waveguide. \label{fig:PerforatedPlateModel}}
\end{figure}
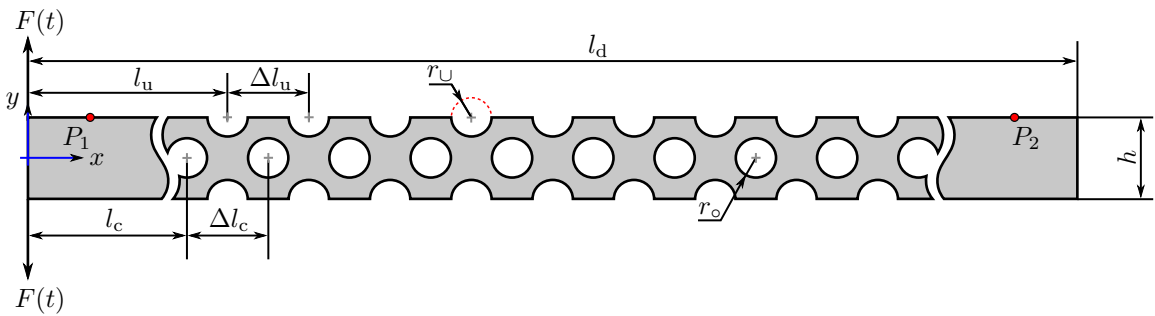

The second example is taken from Ref.~\cite{ArticleJoulaian2014a} and features a waveguide with 37 cutouts (see Fig.~\ref{fig:PerforatedPlateModel}). This is a highly complex geometry, where the generation of a geometry-conforming mesh leads to a large number of finite elements in order to accurately resolve the geometrical details.  It is particularly in scenarios of this complexity that immersed boundary methods emerge as invaluable tools, promising not only simplification but also the automation of the simulation workflow.

The complex waveguide model, as depicted in Fig.~\ref{fig:PerforatedPlateModel} alongside the applied boundary conditions, is characterized by several parameters. First, note that only Neumann boundary conditions are applied, whereas Dirichlet boundary conditions are neglected. For this example, two concentrated forces are applied at the top and bottom surfaces acting in opposite directions. Their time-dependent behavior is given by Eq.~\eqref{eq:Hann}, while the excitation parameters for this simulation are defined as follows: $\bar{F}\,{=}\,10^{8}\,$N, $f_\mathrm{c}\,{=}\,200\,$kHz, and $n\,{=}\,5$. The geometry of the plate is defined by its length $l_\mathrm{d}\,{=}\,600\,$mm and its thickness $h\,{=}\,5\,$mm. Moreover, the semi-circular cutouts, located at the bottom ($y\,{=}\,-\sfrac{h}{2}\,{=}\,-2.5\,$mm) and top ($y\,{=}\,\sfrac{h}{2}\,{=}\,2.5\,$mm) of the plate, are defined by the coordinates of their centroids. The first of these cutouts is positioned at $x\,{=}\,l_\mathrm{u}\,{=}\,154\,$mm. The radius is selected as $r\,{=}\,1\,$mm for all cutouts. Note that the semi-circular cutouts are repeated 11 times, leading to a total of 24 cutouts (12 on the top and 12 on the bottom of the plate). The distance between the centroids of individual cutouts on the same height-level is $\Delta l_\mathrm{u}\,{=}\,\Delta l_\mathrm{c}\,{=}\,4\,$mm. The first fully circular cutout is located at ($x\,{=}\,l_\mathrm{c}\,{=}\,152\,$mm, $y\,{=}\,0\,$mm) and is repeated 12 times, bringing the total count of circular cutouts to 13. Regarding material properties, the waveguide is assumed to be composed of aluminum, characterized by its Young's modulus $E\,{=}\,70\,$GPa, Poisson's ratio $\nu\,{=}\,0.33$, and the mass density $\rho\,{=}\,2700\,\sfrac{\mathrm{kg}}{\mathrm{m}^3}$. The transient simulations are executed under plane strain conditions.

The symmetry of both the geometry and applied boundary conditions relative to the mid-plane of the plate is evident. This implies that the external loading scenario exclusively excites symmetric Lamb wave modes. Additionally, given the symmetry of all obstacles, there will be no mode conversion to antisymmetric modes \cite{ArticleAhmad2012}. This characteristic simplifies the interpretation of the received wave signal, particularly in the context of structural health monitoring applications.

Because of this inherent symmetry in the waveguide, encompassing both its geometry and boundary conditions, it is sufficient to model only half of the plate. This can be achieved by cutting it along the midplane and applying symmetry boundary conditions. While not utilized in this paper, this approach yields a substantial reduction in the number of degrees of freedom and thus, in the computational costs and simulation run time.

\begin{figure}[b!]
	\centering
	\subfloat[Finite element mesh: $3{,}244$ elements \label{fig:PorousPlateFEM_mesh}]{\includegraphics[width=1.0\textwidth]{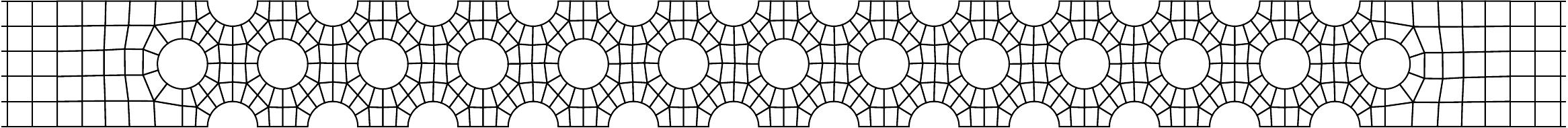}}\\
	\subfloat[Spectral cell mesh: $1{,}920$ elements of which $130$ are cut elements \label{fig:PorousPlateFCM_mesh}]{\includegraphics[width=1.0\textwidth]{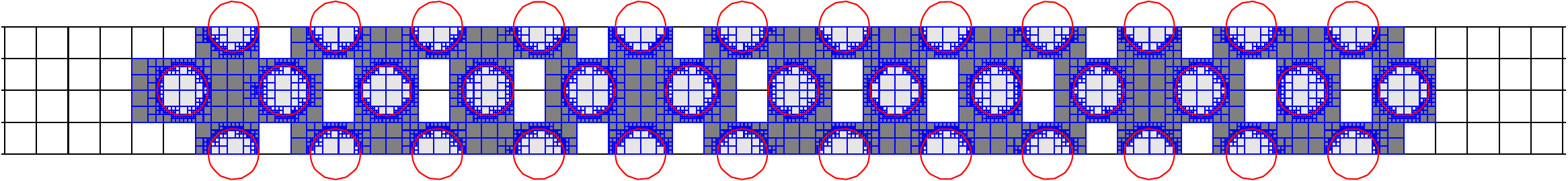}}
	\caption{Spatial discretization of the porous waveguide -- $r_\circ\,{=}\,r_\cup\,{=}\,r\,{=}\,1\,$mm -- using geometry-conforming and immersed boundary methods. The subdivision level for the numerical integration is set to $k\,{=}\,4$. \emph{Color coding for the numerical integration} -- White: conventional finite elements; Dark gray: integration domains located in the physical domain; Light gray: integration domains located in the fictitious domain; Yellow: cut integration domains (leaf cells $\rightarrow$ zoom in 64x). \label{fig:PorousPlate_mesh}}
\end{figure}

The spatial discretization employed for the SEM reference solution consists of $3{,}244$ elements, each having an average size of $h_\mathrm{l}\,{=}\,1\,$mm, as illustrated in Fig.~\ref{fig:PorousPlateFEM_mesh}. Considering the excitation frequency of $200\,$kHz, the shear wavelength for this particular example is $15.79\,$mm. Consequently, there are approximately $15.79$ elements per wavelength available. To generate a precise numerical reference solution, a polynomial degree of $p\,{=}\,8$ is selected, resulting in $n_\mathrm{DOF}\,{=}428{,}040\,$ degrees of freedom.

\begin{figure}[t!]
	\centering
	\subfloat[$u_\mathrm{x}$ at $P_1$]{\includegraphics[scale=1.0]{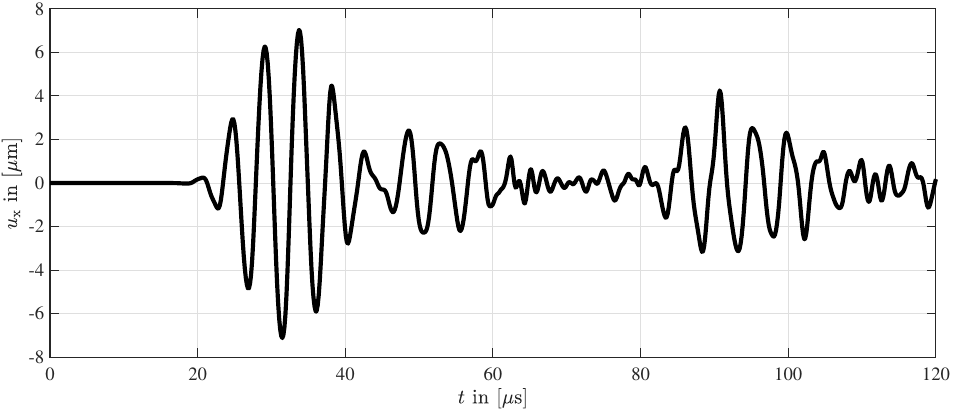}}\\
	\subfloat[$u_\mathrm{y}$ at $P_1$]{\includegraphics[scale=1.0]{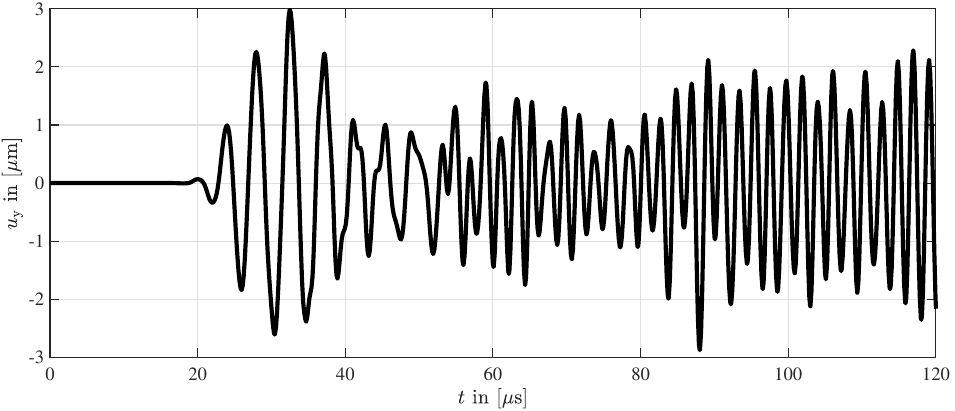}}
	\caption{Time-history of the displacement components at the observation point $P_1$ -- porous waveguide, $r_\circ\,{=}\,r_\cup\,{=}\,r\,{=}\,1\,$mm. \label{fig:TimeHistP1_porousWG}}
\end{figure}

\begin{figure}[t!]
	\centering
	\subfloat[$u_\mathrm{x}$ at $P_2$]{\includegraphics[scale=1.0]{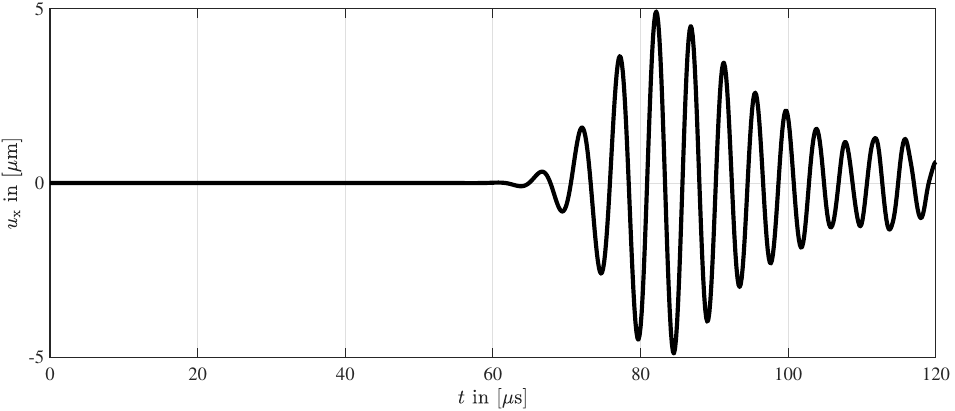}}\\
	\subfloat[$u_\mathrm{y}$ at $P_2$]{\includegraphics[scale=1.0]{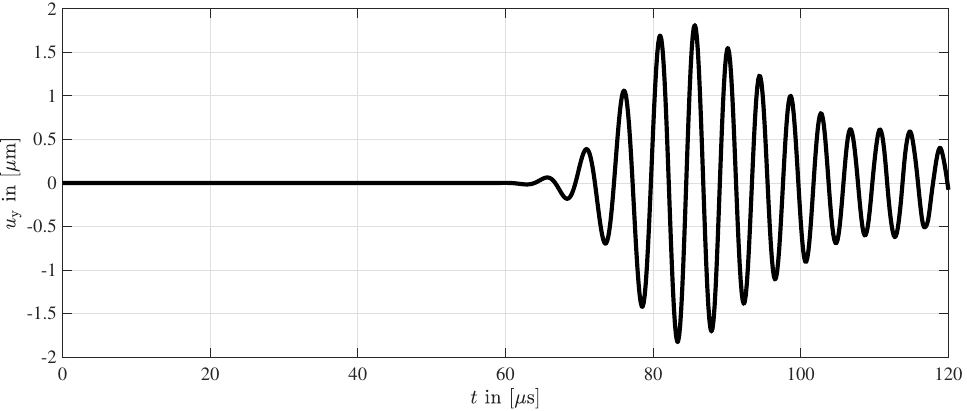}}
	\caption{Time-history of the displacement components at the observation point $P_2$ -- porous waveguide, $r_\circ\,{=}\,r_\cup\,{=}\,r\,{=}\,1\,$mm. \label{fig:TimeHistP2_porousWG}}
\end{figure}

For time integration, a 4\textsuperscript{th}-order accurate Pad\'e-based (implicit) time integration scheme is employed \cite{ArticleSong2022a, ArticleSong2022b}, with a time step size of $\Delta t\,{=}\,10^{-9}\,$s. The reference solutions of the displacement history at two designated observation points $P_1$ and $P_2$ (see Fig.~\ref{fig:PerforatedPlateModel}), located on the top surface of the plate at ($x\,{=}\,100\,$mm, $y\,{=}\,2.5\,$mm) and ($x\,{=}\,300\,$mm, $y\,{=}\,2.5\,$mm), respectively, are are depicted in Figs.~\ref{fig:TimeHistP1_porousWG} and \ref{fig:TimeHistP2_porousWG}. These figures vividly illustrate the complexity of the received wave signal.

In terms of geometry mapping within our in-house code, it is important to note that we employ a subparametric element formulation. This means that we approximate the geometry of the cutouts using 8-noded quadrilateral finite elements of Serendipity type. While it is worth mentioning that more precise descriptions of boundaries can be achieved through techniques like the blending function method (BFM) \cite{BookSzabo1991} or a quasi-regional mapping approach \cite{ArticleKiralyfalvi1997}, these methods are more complex to implement and in case of the BFM often result in superparametric geometry approximations. This introduces complications, particularly in ensuring the correct description of rigid body modes (RBMs). To avoid these complexities, we have opted to stick with the proposed subparametric approach.

In a recent study (Ref.~\cite{ArticleHildebrandt2022}), it was demonstrated that, especially in the context of high-order FEMs, iso- and/or subparametric formulations tend to offer greater robustness and improved accuracy. While the mentioned study favored the isoparametric element concept, we recommend adhering to the proposed subparametric approach due to its simplicity and compatibility with existing mesh generation software, which is predominantly based on low-order finite elements.

In this context, it is important to acknowledge that the subparametric approach introduces a slight loss in accuracy due to an inexact geometry approximation. In practice, this means that a small geometry error has to be accepted and therefore, the discretization error and the geometry error might not be balanced.

Regarding the spatial discretization set-up in the SCM, we employ a grid consisting of $480 \times 4$ elements, each with a dimension of $1.25\,$mm, as illustrated in Fig.~\ref{fig:PorousPlateFCM_mesh}. Specifically, the plate's thickness is discretized using four elements, while in the direction of wave propagation, we employ $480$ elements. This results in an allocation of approximately $12.63$ elements per wavelength within the mesh.

In total, the mesh comprises $1{,}790$ uncut elements and $130$ elements intersected by the physical boundary of the domain of interest. The volume fractions range from $26.4\%$ to $92.7\%$, indicating that we do not anticipate encountering the issue of small cut elements. It is important to note that this outcome represents a fortunate circumstance, as it is generally challenging to avoid generating poorly cut elements for intricate structures.

The time integration is accomplished using the CDM  with a fixed time step size of $\Delta t\,{=}\,10^{-9}\,$s. This value is deliberately chosen to remain below the critical threshold for all simulations, ensuring stability and accuracy in the time-stepping process.

\begin{figure}[t!]
	\centering
	\subfloat[Error $e_{L_2}$]{\includegraphics[scale=1.0]{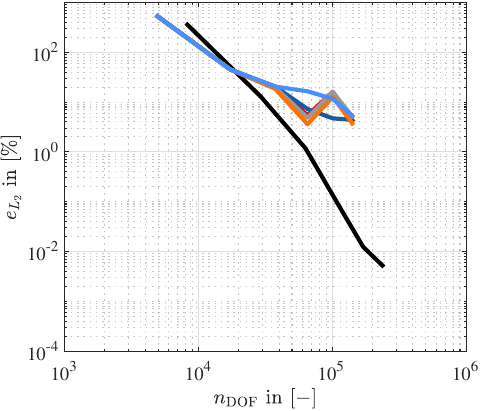}}\hfill
	\subfloat[Critical time step $\Delta t_\mathrm{cr}$]{\includegraphics[scale=1.0]{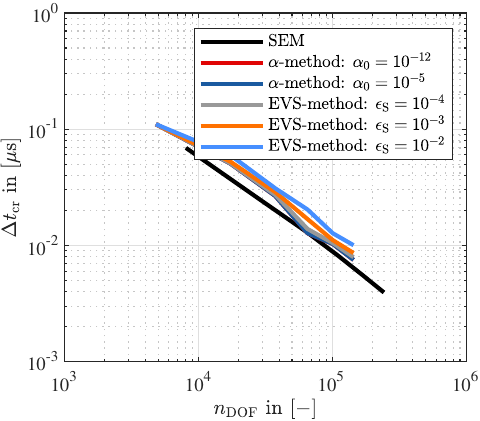}}
	\caption{Evolution of the error and the critical time step size for different stabilization parameters $\alpha_0$ and $\epsilon_\mathrm{S}$ on a fixed mesh ($480\times4$ spectral elements) under \emph{p}-refinement with $p\in[1,6]$ -- porous waveguide, $r_\circ\,{=}\,r_\cup\,{=}\,r\,{=}\,1\,$mm. In the EVS-technique, the stabilization threshold is set to $\epsilon_\lambda\,{=}\,10^{-3}$. \label{fig:PorousWaveguide}}
\end{figure}

In the upcoming analysis, we conduct a \emph{p}-refinement study ($p\in[1,6]$) to observe the convergence characteristics of the SCM across various stabilization techniques and parameters.  To this end, the results are evaluated at the observation point $P_1$, where the error between a reference solution of the displacement history in $x$-direction, denoted as $u^{P_1}_\mathrm{x}$, and its numerical approximation, denoted as $\tilde{u}^{P_1}_\mathrm{x}$, is computed in the $L_2$-norm using Eq.~\eqref{eq:Error_L2_ux_P1}. The outcomes are depicted in Fig.~\ref{fig:PorousWaveguide}.

Upon examining the convergence curves, a noteworthy observation emerges: both stabilization methods exhibit nearly identical performance. There is minimal distinction in terms of the attainable error and the critical time step size. This outcome aligns with our expectations, as the presence of critically cut elements, which typically necessitates stabilization schemes, is absent in this scenario. Consequently, the stabilization techniques do not significantly impact the solution of the equations of motion.

Furthermore, it is worth noting that the critical time step size for the immersed methods, despite their complexity, is actually larger when compared to a geometry-conforming mesh. This phenomenon can be attributed to the intricate shape of the specimen, which prompts the generation of small, distorted elements to accurately represent geometric details.

The reduced accuracy of the immersed methods in this particular example can be attributed to three distinct sources, listed in order of their significance\footnote{Preliminary investigations validating this assumption have been carried out, and their findings are briefly presented in Appendix~\ref{App:Comp_CMM_LMM}.}:
\begin{enumerate}
	\item Lumping of cut elements by means of the HRZ-method.
	\item Insufficient accuracy of the numerical integration of cut elements ($k$ should be increased beyond a value of $4$).
	\item Stabilization of the fictitious domain.
\end{enumerate}

\begin{figure}[t!]
	\centering
	\subfloat[Finite element mesh: $3{,}940$ elements \label{fig:PorousPlate2FEM_mesh}]{\includegraphics[width=1.0\textwidth]{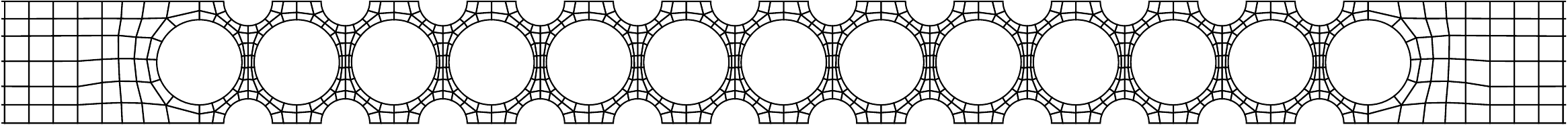}}\\
	\subfloat[Spectral cell mesh: $1{,}920$ elements of which $146$ are cut and $20$ fictitious elements \label{fig:PorousPlate2FCM_mesh}]{\includegraphics[width=1.0\textwidth]{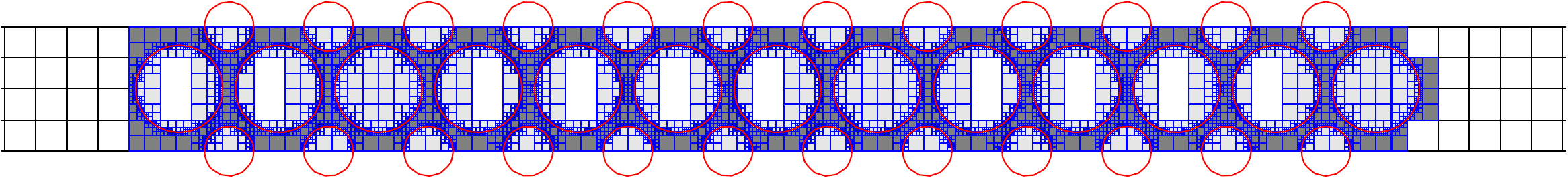}}
	\caption{Spatial discretization of the porous waveguide -- $r_\circ\,{=}\,1.75\,$mm, $r_\cup\,{=}\,1\,$mm -- using geometry-conforming and immersed boundary methods. The subdivision level for the numerical integration is set to $k\,{=}\,4$. \emph{Color coding for the numerical integration} -- White: conventional finite elements; Dark gray: integration domains located in the physical domain; Light gray: integration domains located in the fictitious domain; Yellow: cut integration domains (leaf cells $\rightarrow$ zoom in 64x). \label{fig:PorousPlate2_mesh}}
\end{figure}

In the next step, adjustments are made to the radius of the fully circular cutouts, serving the dual purpose of fine-tuning the volume fraction of cut elements and generating critical cuts. To this end, a radius of $r_\circ\,{=}\,1.75\,$mm is selected for the circular cutouts, while the radius of $r_\cup\,{=}\,1\,$mm remains unchanged for the semi-circular cutouts. The spatial discretizations corresponding to both the geometry-aligned approach (SEM) and the immersed boundary technique (SCM) are presented in Fig.~\ref{fig:PorousPlate2_mesh}.
\begin{figure}[t!]
	\centering
	\subfloat[$u_\mathrm{x}$ at $P_1$]{\includegraphics[scale=1.0]{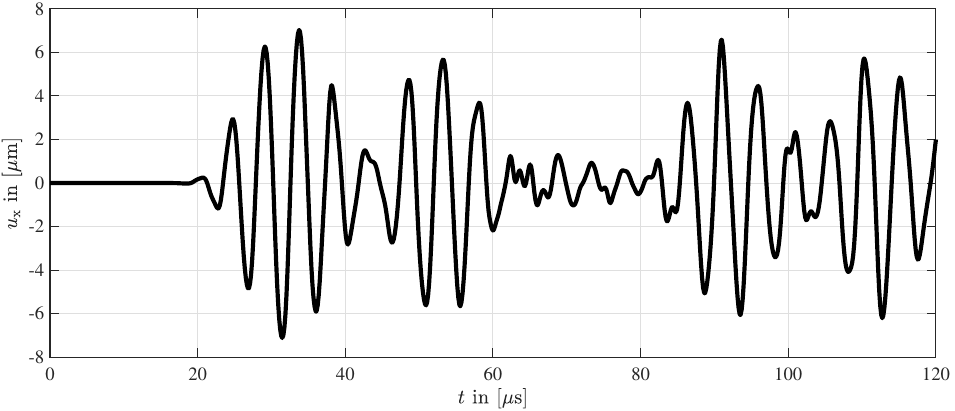}}\\
	\subfloat[$u_\mathrm{y}$ at $P_1$]{\includegraphics[scale=1.0]{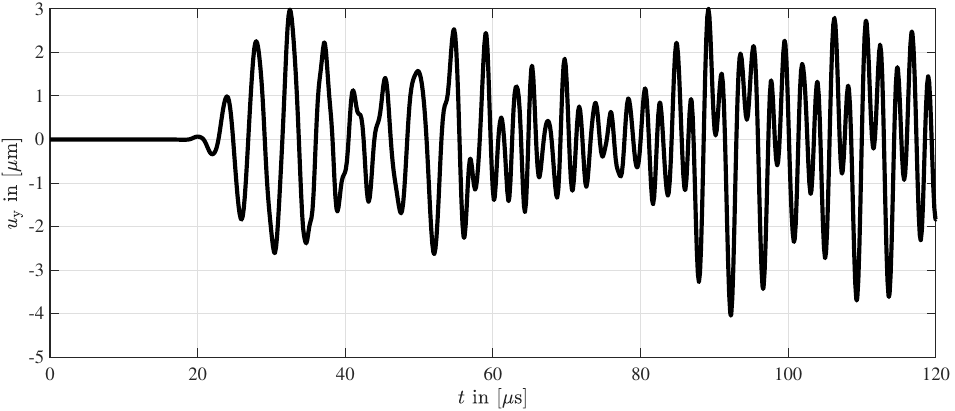}}
	\caption{Time-history of the displacement components at the observation point $P_1$ -- porous waveguide, $r_\circ\,{=}\,1.75\,$mm, $r_\cup\,{=}\,1\,$mm. \label{fig:TimeHistP1_porousWG2}}
\end{figure}%
\begin{figure}[t!]
	\centering
	\subfloat[$u_\mathrm{x}$ at $P_2$]{\includegraphics[scale=1.0]{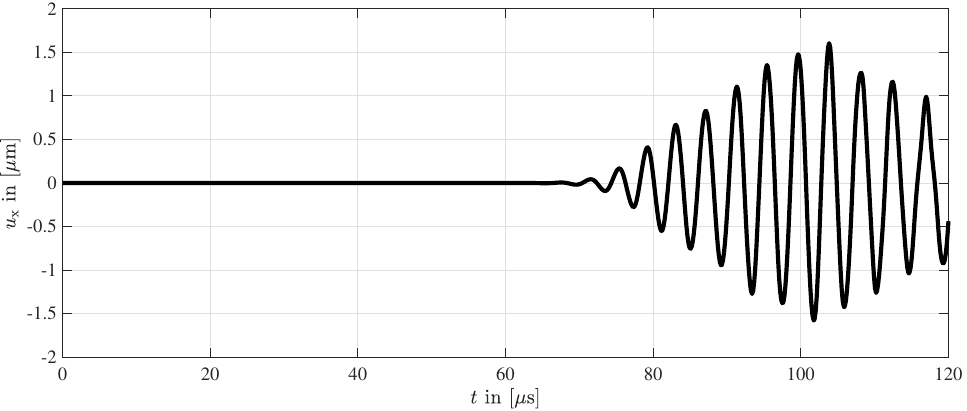}}\\
	\subfloat[$u_\mathrm{y}$ at $P_2$]{\includegraphics[scale=1.0]{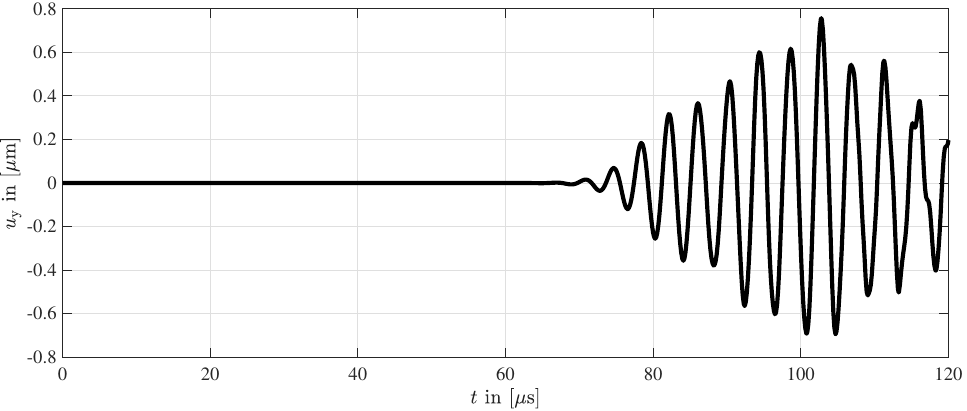}}
	\caption{Time-history of the displacement components at the observation point $P_2$ -- porous waveguide, $r_\circ\,{=}\,1.75\,$mm, $r_\cup\,{=}\,1\,$mm. \label{fig:TimeHistP2_porousWG2}}
\end{figure}%

In the context of the SEM reference solution, the spatial discretization encompasses a total of $3{,}940$ elements, each of order $p\,{=}\,8$, which are characterized by an average element size of $h_\mathrm{l}\,{=}\,1\,$mm, as illustrated in Fig.~\ref{fig:PorousPlate2FEM_mesh}. This configuration yields $n_\mathrm{DOF}\,{=}517{,}736\,$ degrees of freedom for the entire model.

As for the time integration, a time step size of $\Delta t\,{=}\,2\times10^{-10}\,$s is selected to ensure an accurate solution. The reference results for the displacement history at the two designated observation points, as denoted in Fig.~\ref{fig:PerforatedPlateModel} are depicted in Figs.~\ref{fig:TimeHistP1_porousWG2} and \ref{fig:TimeHistP2_porousWG2}.

In our immersed boundary approach, the mesh comprises a total of $1{,}920$ elements, out of which $1{,}754$ remain uncut, and $146$ are intersected by the physical boundary of the domain of interest. Additionally, $20$ elements have been excluded from the analysis as they are fully located within the fictitious domain. The volume fractions within the mesh exhibit a broad spectrum, ranging from as low as $0.02\%$ to as high as $94.4\%$.

This range in volume fractions poses a challenge, particularly when employing explicit methods within an immersed setting, as it imposes restrictions on the critical time step size. Consequently, we decided for a time-stepping approach utilizing the CDM  with a fixed time step size of $\Delta t\,{=}\,2\times10^{-10}\,$s, which is below the critical value for all simulations.

It is worth noting that the time increment has been significantly reduced due to the more pronounced cuts present in the current configuration. This adjustment aligns with the rest of the set-up, which mirrors the previous example where all cutouts were of the same size.

In light of these considerations, a pertinent question arises regarding the wisdom of retaining cells with such small volume fractions. From a robustness standpoint, it would be advisable to eliminate these cells, as it is expected that their removal would have minimal impact on the overall accuracy of the simulation.

\begin{figure}[t!]
	\centering
	\subfloat[Error $e_{L_2}$]{\includegraphics[scale=1.0]{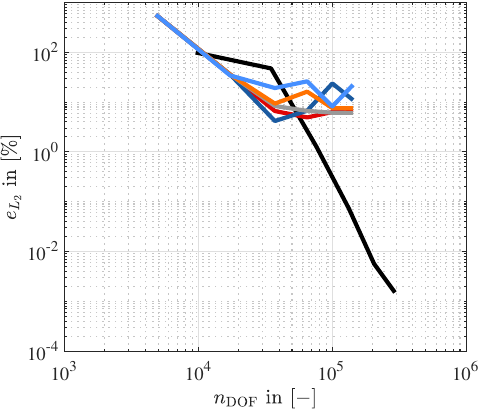}}\hfill
	\subfloat[Critical time step $\Delta t_\mathrm{cr}$]{\includegraphics[scale=1.0]{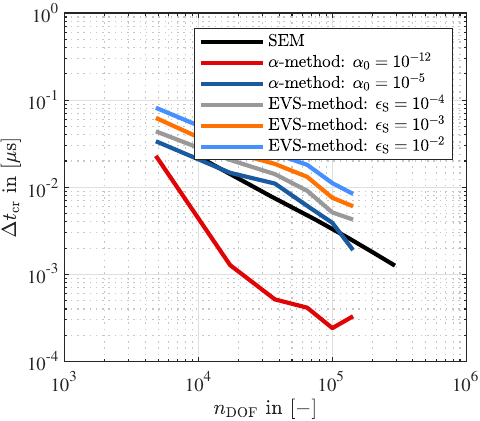}}
	\caption{Evolution of the error and the critical time step size for different stabilization parameters $\alpha_0$ and $\epsilon_\mathrm{S}$ on a fixed mesh ($480\times4$ spectral elements) under \emph{p}-refinement with $p\in[1,6]$ -- porous waveguide, $r_\circ\,{=}\,1.75\,$mm, $r_\cup\,{=}\,1\,$mm. In the EVS-technique, the stabilization threshold is set to $\epsilon_\lambda\,{=}\,10^{-3}$. \label{fig:PorousWaveguide2}}
\end{figure}

In this particular example, we conduct a final comparison between the conventional $\alpha$-stabilization method and the EVS-technique. To facilitate this comparison, we once again employ a \emph{p}-extension approach, wherein the mesh remains fixed while we systematically elevate the polynomial degree of the shape functions. The results of this numerical convergence study, illustrated in Fig.~\ref{fig:PorousWaveguide2}, demonstrate that both stabilization methods yield similar performance concerning the attainable error.

However, it is important to note that in this specific example, which includes severely cut elements, we observe more pronounced effects on the critical time step size. Notably, the $\alpha$-method, when employed with a small stabilization parameter of $\alpha_0\,{=}\,10^{-12}$, proves insufficient to effectively stabilize the fictitious domain. Consequently, the achievable time step sizes are significantly lower compared to cases with geometry-conforming discretization. Interestingly, increasing the value of $\alpha_0$ to $10^{-5}$ brings the time step size closer to that achieved with the SEM.

In contrast, the EVS-technique demonstrates its capability to substantially increase the time step size, surpassing the SEM values. This highlights the advantage of employing a targeted stabilization method over one that stabilizes the entire fictitious domain.

Overall, it is reassuring to observe that the performance of the EVS-technique, in terms of achievable error, closely aligns with that of the established $\alpha$-method across all examples. This translates into an increase in the time step size while maintaining a similar error threshold.

It is worth noting that the reasons for the higher errors compared to geometry-aligned approaches have been previously explained and need not be reiterated here.
\section{Summary and conclusions}
\label{sec:Summary}

In this research, we introduced and implemented an eigenvalue stabilization technique, referred to as the EVS-technique, within the context of dynamic problems. Our findings demonstrate the outstanding performance of this method, particularly in elevating the critical time step size for explicit simulations. Additionally, it effectively addresses ill-conditioning issues, even though these issues are not critical for our specific applications.

Compared to the $\alpha$-method, the EVS-technique offers a more targeted approach, effectively eliminating the source of the observed stability problems in immersed boundary methods. Our results further indicate that the combination of the EVS-technique with the $\alpha$-method does not yield significant benefits in the context of linear explicit dynamics.

Before discussing future research directions, we want to summarize the main \emph{highlights} of the EVS-technique that have been observed in both static (as seen in previous works) and dynamic applications (as explored in the current study):
\begin{enumerate}
	\item Excellent reduction of the condition numbers of both stiffness and mass matrices.
	\item Enhanced robustness in nonlinear static computations, allowing for larger displacements under finite strains.
	\item Highly accurate in nonlinear static computations through the implementation of an iterative force correction scheme allowing for larger stabilization parameters.
	\item Increased critical time step size in (linear) explicit dynamics by exclusively stabilizing the mass matrix.
	\item Targeted approach to stabilizing immersed boundary methods, avoiding  an unnecessary overuse of stabilization in the fictitious domain.
\end{enumerate}
Considering the points above, it is generally recommended to apply the EVS-technique instead of the commonly employed $\alpha$-stabilization.

Based on the numerical findings presented in this study, we highly recommend to set the (mode-independent) stabilization parameter, denoted as $\epsilon_\mathrm{S}$, to a value of $10^{-3}$. Additionally, we advise to set the stabilization threshold, denoted as $\epsilon_\lambda$, to a value of $10^{-3}$. These parameter settings have proven to be robust and reliable in our analysis.

A promising avenue for future research lies in the development of force correction techniques tailored for transient simulations. These techniques are designed to alleviate the detrimental effects of introducing additional stabilization terms to the system, terms that are not inherent to the original mathematical problem. As a result, the loss in accuracy can be partially mitigated.

In the realm of non-linear static analyses, a viable iterative force correction technique was previously proposed in Ref.~\cite{ArticleGarhuom2022}, demonstrating a noteworthy enhancement in accuracy. However, for transient simulations, especially within the context of explicit time-stepping methods, iterative approaches may not be as practical. This is because the potential gains in accuracy through iterations could offset the advantages of achieving higher critical time step sizes. The computational costs per time step are directly proportional to the number of iterations required to enhance accuracy.

Consequently, the need arises for the development of more sophisticated force correction procedures tailored specifically for dynamic simulations. This research endeavor holds significant promise for improving the accuracy and efficiency of transient simulations in complex engineering and scientific applications.

Implementing the EVS-technique, along with force correction terms, in dynamic simulations introduces an additional layer of complexity. This complexity stems from the fact that the force correction terms are dependent on the displacement, velocity, and acceleration fields. Incorporating these functions into the semi-discrete equations of motion necessitates the development and thorough analysis of tailored time integration schemes to ensure their numerical stability and accuracy.

In essence, unlike their application in (nonlinear) statics, it is not merely a matter of adding force correction terms in dynamic scenarios. Instead, it requires a more nuanced approach that involves devising specialized time integration methods capable of handling the intricate interactions between these correction terms and the evolving dynamic fields. This challenge underscores the need for dedicated research efforts in the realm of dynamic simulations when implementing the EVS-technique with force correction.

It is worth noting that the semi-discrete equations of motion are conventionally solved using direct time-stepping methods. Consequently, the introduction of force correction terms brings about alterations in the expressions for the amplification matrix and the load operator within the chosen time integration scheme. This implies that the implementation of a correction technique can potentially transform a time integrator that was previously unconditionally stable into one that is no longer unconditionally stable but unconditionally unstable. Even in cases where stability is maintained, there is a need to establish new estimates for the critical time step size, especially in the context of explicit dynamics.

In this context, there is a tangible risk that the introduction of force correction terms may actually lead to a reduction in the stability limit. In such instances, the pursuit of improved accuracy might be outweighed by the restrictions imposed on the time step size. Consequently, the question of how to effectively approximate or compute the correction terms becomes a delicate and intricate challenge, warranting dedicated research efforts, which the authors are currently engaging in.
\paragraph{Acknowledgment} The authors (SE, AD, and LR) gratefully acknowledge the financial support provided by the German Research Foundation (Deutsche Forschungsgemeinschaft -- DFG) through the research grants EI 1188/3-1 (project number: 497531141) and DU 405/20-1 (project number: 503865803).
%
\appendix
\renewcommand\thefigure{\Roman{figure}}  
\setcounter{figure}{0}
\renewcommand\theequation{\Roman{equation}}  
\setcounter{equation}{0}  
\section{Stabilization parameter}
\label{App:Eps_S}
In this part of the appendix, we briefly introduce different formulations of the stabilization parameter taken from recent literature \cite{ArticleLoehnert2014, ArticleGarhuom2022, PhDGarhuom2023}. These different approaches are included for the sake of completeness and to provide the reader with an overview of available options. However, a detailed analysis of the pros and cons of each approach are out of the scope of this contribution.

The first approach has been suggested by L\"ohnert in the original publication on the EVS-technique \cite{ArticleLoehnert2014}. Here, a functional dependence of the actual stabilization parameter and the magnitude of the eigenvalue to be stabilized is introduced:
\begin{equation}
	_1\epsilon^\lambda_{\mathrm{S},i} = \left(\epsilon_\mathrm{S} \lambda_\mathrm{max} - \hat{\lambda}_{\mathrm{u},i}\right)f^+\left(\epsilon_\mathrm{S} \lambda_\mathrm{max} - \hat{\lambda}_{\mathrm{u},i}\right)\qquad \forall\; i \in [1,n_\mathrm{u}]
	\qquad\text{with}\quad f^+(x) = 
	\begin{cases}
		0 & \forall\; x < 0 \\
		1 & \forall\; x \ge 0
	\end{cases}\,.
	\label{eq:StabFac2}
\end{equation}%

\begin{figure}[b!]
	\centering
	\includegraphics[scale=1.0]{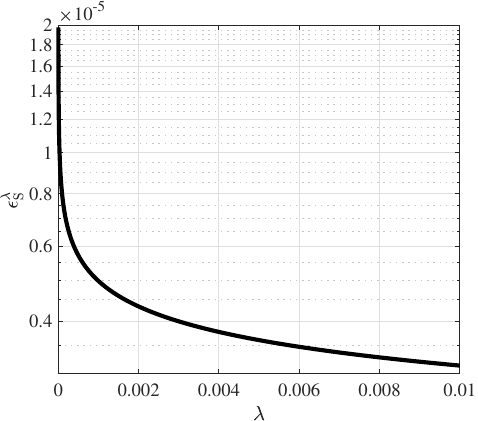}
	\caption{Example of the functional dependence of the stabilization factor ${}_2\epsilon^\lambda_\mathrm{S}$ on the eigenvalue $\lambda$ for $\epsilon_\mathrm{S}\,{=}\,1 \times 10^{-4}$, $\beta\,{=}\,5$, and $n_{\epsilon_\mathrm{S}}\,{=}\,80$.\label{fig:StabFac3}}
\end{figure}

Another approach, proposed by Garhuom et al. in Ref.~\cite{ArticleGarhuom2022}, suggests the following empirically determined relation:
\begin{equation}
	_2\epsilon^\lambda_{\mathrm{S},i} = \cfrac{\epsilon_\mathrm{S}}{n_{\epsilon_\mathrm{S}}\sqrt[\beta]{\hat{\lambda}_{\mathrm{u},i}}}\qquad \forall\; i \in [1,n_\mathrm{u}]\,.
	\label{eq:StabFac3}
\end{equation}%
The additional factors $\beta$ and $n_{\epsilon_\mathrm{S}}$ are introduced to tune the stabilization factor to different applications. In Ref.~\cite{ArticleGarhuom2022}, a good performance has been demonstrated for $\beta\,{=}\,5$ and $n_{\epsilon_\mathrm{S}}\,{=}\,80$, see Fig.~\ref{fig:StabFac3}. Although the magnitude of the eigenvalues is taken into account by  Eq.~\eqref{eq:StabFac3}, there is still a potential issue when different systems of units are utilized. Keep in mind that the mode shapes are invariant to a change in the system of units, while the eigenvalues will be different. Due to the nonlinear characteristics of the stabilization factor introduced in Eq.~\eqref{eq:StabFac3}, one would obtain different relative stabilization matrices despite that fact that the physics of the structure has not changed. To address this issue and ensure a wider range of applicability, it is also recommended to introduce the scaling approach proposed in Sect.~\ref{subsec:Scaling}.

A third idea, put forward in Ref.~\cite{PhDGarhuom2023}, consists in establishing a functional dependence of the stabilization factor on the material parameters, e.g., the first Lam\'{e} constant $\lambda_\mathrm{L}$, and the volume fraction of a cut element $\chi$. This results in a constant stabilization factor
\begin{equation}
	_3\epsilon_{\mathrm{S},i} = \epsilon_\mathrm{S}\, \lambda_\mathrm{L} \, \left( 1 - \chi \right)^\beta\,.
	\label{fig:StabFac4}
\end{equation}
The additional parameter $\beta$ is again introduced to tune the properties of the proposed scheme. In Ref.~\cite{PhDGarhuom2023}, it is suggested to use $\beta\,{=}\,1$.

At this point, we want to note that our approach, which consists in using a constant value for the stabilization factor $\epsilon_\mathrm{S}$ in conjunction with the scaling method described in Sect.~\ref{subsec:Scaling}, is closely related to the third approach given by Eq.~\eqref{fig:StabFac4}. By means of the scaling relation we implicitly take the material properties into account. While we do not consider the volume fraction of the cut element in our formulation, one can easily infer that for badly cut elements, which is the primary area of application of the stabilization scheme, the term $(1 - \chi)$ quickly goes to 1. Therefore, it is expected to observe a similar performance.
\section{Alternative expression of the stabilized/modified matrix}
\label{App:AlternativeExp}
In this part of the appendix, an alternative expression for the system stabilization matrices is give.We recall that our objective is to stabilize the system matrices for cut elements in immersed boundary methods. This is typically accomplished by introducing an additional stabilization term $\mathbf{A}^\mathrm{S}$ to the generic system matrix $\mathbf{A}$. Therefore, our stabilized or modified matrix $\mathbf{A}^\mathrm{Mod}$ is defined as
\begin{equation}
	\mathbf{A}^\mathrm{Mod} = \mathbf{A} + \mathbf{A}^\mathrm{S}\,.
\end{equation}

When employing the EVS-technique, the system matrix is initially transformed into its spectral form
\begin{equation}
	\mathbf{A} = \fatgreek{\Phi}\, \fatgreek{\Lambda}\, \fatgreek{\Phi}^\mathrm{T}\,,
\end{equation}
where $\fatgreek{\Phi}$ and $\fatgreek{\Lambda}$ represent the matrix of mode shapes (stored column-wise) and the diagonal matrix of eigenvalues, respectively. The stabilization matrix $\mathbf{A}^\mathrm{S}$ can be expressed in a similar manner as
\begin{equation}
	\mathbf{A}^\mathrm{S} = \fatgreek{\Phi}\, \fatgreek{\epsilon}^\mathrm{S}\, \fatgreek{\Phi}^\mathrm{T}\,,
\end{equation}
with $\fatgreek{\epsilon}^\mathrm{S}$ representing the diagonal matrix of stabilization parameters. The stabilization parameter for mode $i$ is denoted as $\epsilon_i^\mathrm{S}$, and it can differ for each mode. Thus, the matrix of stabilization parameters takes the following form
\begin{equation}
	\fatgreek{\epsilon}^\mathrm{S} = n^\mathrm{a}_\alpha
	\begin{bmatrix}
		\epsilon_1^\mathrm{S} & & & & \\
		& \epsilon_2^\mathrm{S} & & & \\
		& & \epsilon_3^\mathrm{S} & & \\
		& & & \ddots & \\
		& & & & \epsilon_{n_\mathrm{dim}}^\mathrm{S}
	\end{bmatrix},
\end{equation}
where $n_\mathrm{dim}$ denotes the dimension of the square matrix. It is important to recognize that the scaling parameter $n^\mathrm{a}_\alpha$, as defined in Sect.~\ref{subsec:Scaling}, remains a constant value. For modes that do not require stabilization, $\epsilon_i^\mathrm{S}$ takes a value of $0$. Otherwise, it might be a function of the corresponding eigenvalue $\lambda_i$, i.e., $\epsilon_k^\mathrm{S}\,{\ne}\,\epsilon_l^\mathrm{S}$, or a constant value, i.e., $\epsilon_k^\mathrm{S}\,{=}\,\epsilon_l^\mathrm{S}\,{=}\,\epsilon^\mathrm{S}$. Therefore, we can express the modified matrix as follows
\begin{equation}
	\mathbf{A}^\mathrm{Mod} = \fatgreek{\Phi} \left(\fatgreek{\Lambda} + \fatgreek{\epsilon}^\mathrm{S} \right) \fatgreek{\Phi}^\mathrm{T}\,.
	\label{eq:Amod}
\end{equation}
The result is evidently equivalent to the derivation presented in Sect.~\ref{sec:EVST}, albeit expressed differently. Nevertheless, by examining Eq.~\eqref{eq:Amod}, it becomes apparent that the EVS technique effectively enhances the stability of the system matrices by increasing the eigenvalues associated with the physically meaningless zero eigenspace.

Numerical analyses demonstrate that in the presence of a uniform stabilization parameter denoted as $\epsilon_\mathrm{S}$, the stabilization exclusively impacts the main diagonal, leaving the off-diagonal elements unaffected. This situation is markedly different in scenarios involving varying stabilization parameters $\epsilon_i^\mathrm{S}$ or incomplete mode stabilization, where the stabilization also extends its influence to the off-diagonal components within the system matrix.
\section{Influence of mass lumping applied to cut element on the attainable accuracy in dynamics}
\label{App:Comp_CMM_LMM}
In this part of the appendix, we briefly investigate the reasons for the impaired convergence behavior observed in the numerical examples. While it is well-documented that mass lumping in the context of spectral elements leads to optimal convergence, the observed deterioration in accuracy is likely connected to the treatment of cut elements in the spatial discretization. Therefore, we conduct a more focused investigation on the porous plate case, as discussed in Sect.~\ref{sec:TransientCircHoles}, with the critical setup parameters: $r_\circ\,{=}\,1.75\,$mm, $r_\cup\,{=}\,1\,$mm. To this end, we select a specific spatial discretization consisting of $480 \times 4$ elements with a polynomial degree of $p\,{=}\,6$, resulting in $n_\mathrm{DOF}\,{=}142{,}950\,$ degrees of freedom. We also make adjustments to the mass matrix formulation and the subdivision level for numerical integration ($k$). For stabilization, we employ the EVS-technique (variant \emph{2b}) with parameters $\alpha_0\,{=}\,0$, $\epsilon_\mathrm{S}\,{=}\,10^{-2}$, and $\epsilon_\lambda\,{=}\,10^{-3}$.

\begin{figure}[t!]
	\centering
	\includegraphics[width=1.0\textwidth]{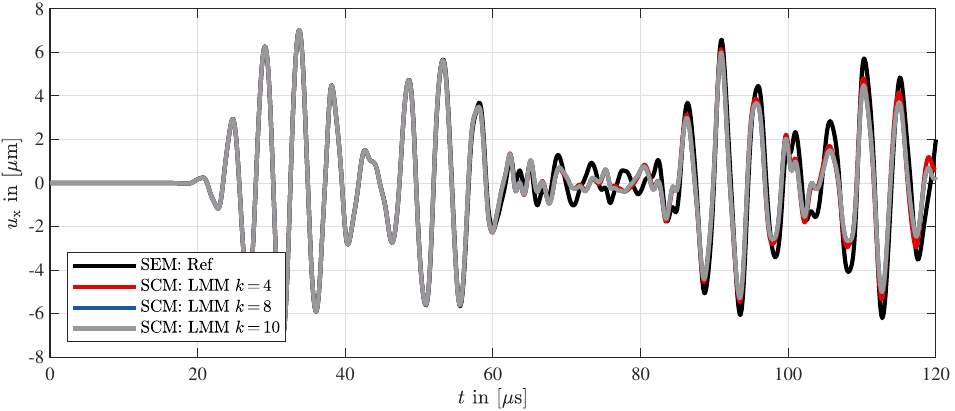}
	\caption{Influence of the numerical integration accuracy on the displacement signal. \label{fig:Influence_k}}
\end{figure}

First, we investigate the impact of numerical integration accuracy, particularly by increasing the value of $k$ from $4$ to $10$. The results are depicted in Fig.~\ref{fig:Influence_k}, demonstrating that a subdivision level of $k\,{=}4$ is insufficient for our needs. While some improvement is noticeable with higher values of $k$, it is evident that the overall trend of the signal remains inaccurately approximated. Thus, despite the fact that enhancing numerical integration accuracy is beneficial, it alone does not resolve the issue of poor convergence and is also not the main contributor.

\begin{figure}[b!]
	\centering
	\includegraphics[width=1.0\textwidth]{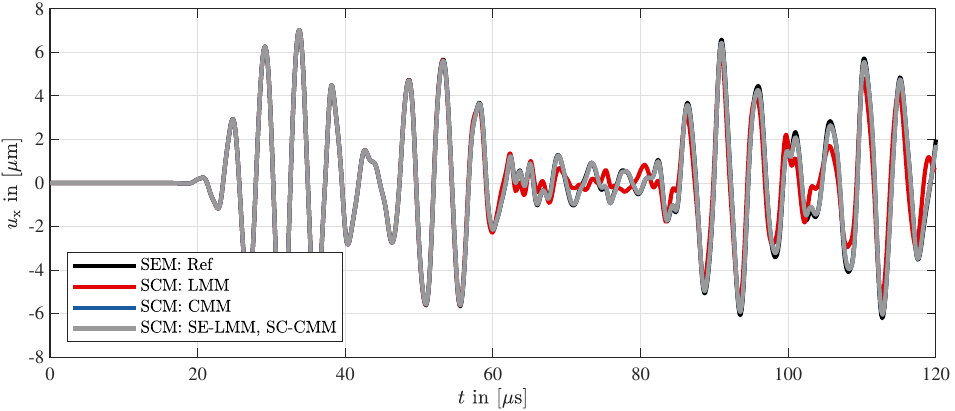}
	\caption{Influence of the employed mass matrix formulation on the displacement signal. \label{fig:Influence_M}}
\end{figure}

As a second attempt to explain the suboptimal convergence rate, we adjust the mass matrix formulation, while a subdivision level of $k\,{=}\,4$ is employed. To this end, we consider the following settings:
\begin{enumerate}
	\item Use HRZ-lumping for cut elements and nodal quadrature for spectral elements (standard approach) $\rightarrow$ SCM: LMM.
	\item Use the consistent mass matrix for all elements $\rightarrow$ SCM: CMM.
	\item Use the consistent mass matrix for cut elements and nodal quadrature for spectral elements  $\rightarrow$ SCM: SE-LMM, SC-CMM.
\end{enumerate}
The results, depicted in Fig.~\ref{fig:Influence_M}, reveal that the consistent mass matrix formulation produces outcomes nearly indistinguishable from the reference solution, with slight discrepancies attributed to the finer spatial and temporal resolutions of the reference and numerical integration of cut elements. Even when employing a consistent mass matrix formulation for cut cells only, we achieve a remarkable agreement with the reference. Only when cut elements are lumped using the HRZ-scheme do we observe a deterioration in the results.

\begin{figure}[t!]
	\centering
	\includegraphics[width=1.0\textwidth]{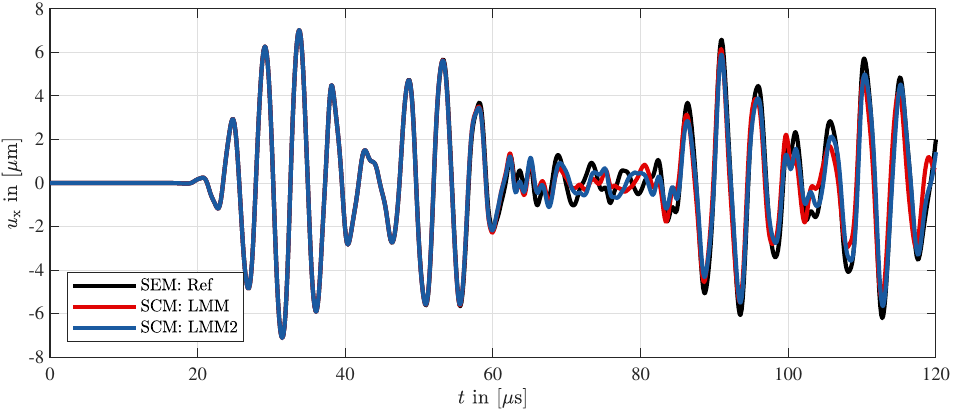}
	\caption{Influence of the HRZ-lumping algorithm. \label{fig:Influence_LMM}}
\end{figure}

Moreover, it has been even observed that the results vary on the specific HRZ-lumping procedure employed (see Fig.~\ref{fig:Influence_LMM}). The red curve illustrates results derived by first computing the consistent mass matrix of cut elements and subsequently applying the HRZ-scheme. In contrast, the blue curve represents outcomes achieved through the direct application of the HRZ-technique to individual subdomain contributions, followed by the summation of the pre-diagonalized matrices to obtain the lumped mass matrix of the cut element.  This, once again, emphasizes the imperative need for dedicated research efforts aimed at addressing the unresolved question of mass lumping in the context of immersed boundary methods.
\bibliographystyle{ieeetr}
\bibliography{./bib/references}
\end{document}